\documentclass{amsart}

\usepackage{amsbsy}
\usepackage{amsfonts}
\usepackage{latexsym}
\usepackage[dvips]{graphics}
\usepackage{epsfig}
\input{psfig.sty}

\usepackage{amsmath}
\usepackage{amssymb}
\usepackage{amscd}
\usepackage{amsthm}

\numberwithin{equation}{subsection}
 
\theoremstyle{plain}

\newtheorem{Le}[equation]{Lemma}

\newtheorem{N}{Notation}

\theoremstyle{remark}

\theoremstyle{definition}

\newtheorem{Def}[equation]{Definition}

\newcommand{\Z}{\mathbb Z}
\newcommand{\R}{\mathbb R}
\renewcommand{\H}{{\mathbb H}}
\renewcommand{\N}{{\mathbb N}}

 \def\Area{\text{Area}}
\def\supp{\text{supp}}

\def\Link{\text{Link}}
\def\length{\text{length}}
\def\sup{\text{sup}}
\def\liminf{\text{liminf}}
\def\interior{\text{interior}}
\def\image{\text{image}}
\def\CAT{\text{CAT}}

\def\id{\text{id}}

\def\De{\Delta}

\def\ra{\rightarrow}

\def\e{\emph}
\def\i{\infty}
\def\p{\partial}
\def\b{\begin}

\begin{document}

\title{    \flushleft{\bf{Mostow rigidity for Fuchsian buildings}}      }

\maketitle


\noindent 
Xiangdong Xie\newline
Department of Mathematics,
Washington University,
St.Louis, MO 63130.\newline
Email:   {  xxie@math.wustl.edu}

\noindent
Current address: Department of Mathematics,
University of Cincinnati,  Cincinnati, Ohio 45221.



\pagestyle{myheadings}

\markboth{{\upshape Xiangdong Xie}}{{\upshape   Mostow rigidity for Fuchsian buildings}}

      
\vspace{3mm}

\noindent
{\small {\bf Abstract.}
We show that if a  homeomorphism between the ideal boundaries of two Fuchsian buildings preserves the combinatorial cross ratio 
almost everywhere, then it extends to  an isomorphism between the Fuchsian buildings.  
It  follows  that    Mostow rigidity  holds  for Fuchsian
    buildings:  if a group acts properly and cocompactly on two 
 Fuchsian buildings $X$ and $Y$, then $X$ and $Y$  are equivariantly isomorphic.  }


\vspace{3mm}



\vspace{3mm}

\tableofcontents

\setcounter{subsection}{0}
\setcounter{subsubsection}{0}


\subsection{Introduction}




The classical Mostow rigidity theorem states that if $X$ and $Y$ are two irreducible symmetric spaces with nonpositive sectional curvature  and dimension $>2$, and  $G$ is a group acting properly and cocompactly
 on both $X$ and $Y$, then $X$ and $Y$ are equivariantly isometric (after rescaling of the metrics).
 In this paper we consider the same question  
for a class of singular spaces --- Fuchsian buildings.

Let $R$ be a convex compact polygon in the hyperbolic plane  whose angles are of the form $\pi/m$, $m\in \Z$, $m\ge 2$.
The group generated by the reflections about the edges of $R$ is a Coxeter group $W$, 
 and $W(R)$ is a tessellation  of the hyperbolic plane.  The hyperbolic plane  equipped with such a 
 tessellation is the Coxeter complex of $W$.  A two dimensional hyperbolic building  $\De$  with chamber  $R$ (see Section \ref{pre}
 for the definition) is a piecewise hyperbolic 2-complex which has nice local structure and contains 
 plenty of subcomplexes isomorphic to the Coxeter complex. Each  $2$-cell  of $\De$ is isometric to $R$ and is called 
 a chamber. Each subcomplex of $\De$ isomorphic to the Coxeter complex is called an apartment.
 In particular, any two chambers of $\De$ are contained in a common apartment.  
            $\De$  with  the path metric   is   a 
 $\CAT(-1)$ space. Each geodesic in $\De$ is also contained in an apartment. 
 The vertex links of  $\De$  
are generalized polygons (see \cite{R} or  Section  \ref{pre}).  There is a combinatorial map $f:\De\ra R$ 
 such that the restriction of $f$ to each chamber is an isomorphism.  A two dimensional hyperbolic building  $\De$
 is a   \e{Fuchsian building} if  for each edge $e$ of $R$ there is an integer $q_e\ge 2$ such that 
 any edge of $\De$ that is mapped to $e$ under $f$ is contained in exactly $q_e+1$ chambers.







\b{Th}\label{main1}
{Let $\De_1$, $\De_2$  be two Fuchsian buildings, and $G$ a group acting properly  and cocompactly 
on both $\De_1$  and  $\De_2$. Then $\De_1$  and  $\De_2$   are equivariantly isomorphic.}

\end{Th}

One cannot claim that $\De_1$  and  $\De_2$   are equivariantly isometric as in the classical case.
 Two convex compact polygons $R_1$ and $R_2$ in the hyperbolic plane are 
  said  to  be  \e{angle-preserving isomorphic}
if there is a combinatorial
isomorphism from $R_1$ to $R_2$ that preserves the angles at the vertices.
Let  $\De$ be  a Fuchsian building with chamber $R$,   and $R'$  a polygon angle-preserving isomorphic to
 $R$.  Then it is not hard to see that one obtains a Fuchsian building  $\De'$ 
by replacing 
 each chamber of $\De$ by a copy of $R'$.  The building $\De'$ is clearly isomorphic to $\De$,
 but is not isometric to $\De$  when $R$ and $R'$ are not isometric.  
A  convex compact polygon $P$ is called \e{normal} if it has an inscribed circle that touches all
 its edges. By a result in   A. Beardon's book (see \cite{Bea}, theorem 7.16.2.) each  convex compact polygon
 $P$ is angle-preserving isomorphic to a unique normal polygon $P'$.
Theorem \ref{main1} can be restated as follows.

\b{Th}\label{main'}
{Let $\De_1$, $\De_2$  be two Fuchsian buildings whose chambers are normal, and $G$ a group acting properly  and cocompactly 
on both $\De_1$  and  $\De_2$. Then $\De_1$  and  $\De_2$   are equivariantly isometric.}

\end{Th}




Our proof  of  Theorem \ref{main1}
 uses D. Sullivan's approach (\cite{S2}) to Mostow rigidity 
and M. Bourdon's   work (\cite{Bo1}, \cite{Bo2})
 on   
Fuchsian buildings.
D. Sullivan's approach  has also been 
used by   M. Bourdon (\cite{Bo1})
 to establish Mostow rigidity for right angled Fuchsian buildings.

 There are two steps in  D. Sullivan's approach.  First notice that under the conditions of Theorem \ref{main1},
     there is an  equivariant homeomorphism
 $h:\p\De_1\ra \p\De_2$  between the ideal boundaries.
In the first step one needs to show  that 
there are conformal measures on  $\p\De_1$  and   $\p \De_2$  such that 
 they have the same dimension and $h$ is nonsingular with respect to the conformal measures(that is, 
  $A\subset \p\De_1$ has measure 0 if and only if $h(A)$ has measure 0). 
  D. Sullivan's ergodic argument then shows that $h$ preserves the cross ratio almost everywhere
 with respect to the conformal measures.  In the second step one needs to  establish the implication
\lq\lq $h$ preserves the cross ratio almost everywhere"
 $\Longrightarrow$
\lq\lq $h$ extends 
 to an  isomorphism from $\De_1$ to $\De_2$". 
Step one has been established by M. Bourdon(\cite{Bo2}).
 He constructed combinatorial metrics   on the ideal boundaries   of Fuchsian buildings.  
 The conformal measures   on the ideal boundaries  are the 
    Hausdorff measures of the combinatorial metrics.   
M. Bourdon  computed the Hausdorff dimension of the combinatorial metric and showed that the 
combinatorial metric  realized the conformal dimension.  In this paper we establish the second step.

\b{Th}\label{main}
{Let $\De_1$, $\De_2$  be two Fuchsian buildings, and $h:\p\De_1\ra \p\De_2$  a homeomorphism that preserves the combinatorial cross ratio almost everywhere.    Then  $h$ extends to an isomorphism from $\De_1$  to  $\De_2$, that is, there is an isomorphism
  $f:  \De_1\ra \De_2$ such that $h$ agrees with the  boundary map of $f$.}

\end{Th}

There is a much stronger rigidity problem concerning Fuchsian buildings. 
 In fact M. Bourdon and  H. Pajot (\cite{BP2})  have  the following  quasi-isometric rigidity conjecture.

\b{conjecture}\label{quasi-ir}
{Let $\De_1$, $\De_2$  be two Fuchsian buildings whose chambers are normal. Then any quasi-isometry
$f:\De_1\ra \De_2$  lies at a finite distance from an isometry.}

\end{conjecture}



  M. Bourdon and  H. Pajot (\cite{BP2})
 have   established this  conjecture when   the chambers of $\De_1$, $\De_2$  
 are  right angled  regular polygons. Recently B. Kleiner has informed the author of his new  result on the quasiconformal group
of  a  $Q$-regular, $Q$-Loewner space.   By an argument of  M. Bourdon and  H. Pajot (\cite{BP2}),  
  this result   and Theorem \ref{main} imply that Conjecture \ref{quasi-ir} holds.

We next explain the proof of Theorem \ref{main}.

Let  
$h:\p\De_1\ra \p\De_2$   be   a homeomorphism that preserves the combinatorial cross ratio 
   (see Section \ref{dcror} for definition)
almost everywhere. 
 It  is not hard to show that for $\xi, \eta$ in the ideal boundary of a Fuchsian 
  building  $\De$, whether $\xi\eta$ is contained in the 1-skeleton can be 
  detected by the behavior of the combinatorial cross ratio (see Lemma \ref{geo1ske}). 
It follows that  if $\xi\eta$ ($\xi, \eta\in \p \De_1$) is a geodesic contained in the 
 1-skeleton of $\De_1$, then $h(\xi)h(\eta)$  is contained in the 1-skeleton of $\De_2$.   We call $h(\xi)h(\eta)$
  the image of $\xi\eta$.
 The key is to   show that  
for any   fixed vertex $v\in \De_1$, the images of  all the geodesics in the 1-skeleton of $\De_1$ through $v$ 
   intersect in a   unique  vertex $w$ in $\De_2$.  Then it is not hard to see that the map $v\ra w$ is a 1--to--1 map 
  between the vertices that extend to an isomorphism  from $\De_1$ to $\De_2$.
One first needs to show that if $c_1$, $c_2$ are two  intersecting geodesics contained in the 
  1-skeleton of an 
  apartment 
of $\De_1$ then their images intersect in a vertex.

Let $\De$ be a Fuchsian building,   and $\eta\in\p\De$ a point 
represented by a ray $\sigma:[0,\i)\ra \De$
contained in the 1-skeleton.
Let $\xi_1, \xi_2\in \De\cup \p\De-\{\eta\}$ be such that
 $\xi_1\eta\cap \sigma=\xi_2\eta\cap \sigma=\phi$. 
We parameterize  $\xi_i\eta$ ($i=1,2$) 
by $\sigma_i:  I_i\ra \De$ ($I_i\subset (-\i,\i)$) such that $\sigma_i(t)$ and $\sigma(t)$ lie on the same horosphere
 centered at $\eta$.  Since $\De$ is $\CAT(-1)$, we have $d(\sigma_i(t), \sigma(t))\ra 0$ as $t\ra \i$. 
We say $\xi_1$ and $\xi_2$   \e{lie at different sides of 
 $\eta$} if $ \sigma_1(t)$  and  $ \sigma_2(t)$      lie in  different  chambers  for all $t$.

\b{Def}\label{atdsftg}
{Let $\xi_1, \xi_2, \eta_1, \eta_2\in \p\De$  be pairwise distinct such that the geodesics 
 $\xi_1\xi_2$, $\eta_1\eta_2$ are contained in the 1-skeleton.  We say the two geodesics 
$\xi_1\xi_2$, $\eta_1\eta_2$ \e{are  at different sides} if 
$\xi_1, \xi_2$   lie at different sides of $\eta_i$ ($i=1,2$) and 
$\eta_1, \eta_2$   lie at different sides of $\xi_i$ ($i=1,2$).}

\end{Def}

A  crucial lemma (see Lemma \ref{difes})  says that 
for $\xi_1, \xi_2\in \p\De$    and  $\eta\in \p\De$ represented by 
   a geodesic ray $\sigma:[0,\i)\ra \De$
contained in the 1-skeleton
with 
  $\xi_1\eta\cap \sigma=\xi_2\eta\cap \sigma=\phi$,   
  whether  $\xi_1$, $\xi_2$ lie at different sides of  $\eta$   
  can be detected by   the combinatorial  cross ratio.
It follows that for $\xi_1, \xi_2, \eta_1, \eta_2\in \p\De_1$,  the two geodesics 
 $\xi_1\xi_2$, $\eta_1\eta_2$ are  at different sides  if and only if their images are at different sides. 
 In particular it implies that  the images of two intersecting geodesics contained in the 1-skeleton
 of 
  an  apartment 
of $\De_1$  are at different sides.  It turns out that in most cases,  two 
 geodesics at different sides 
 must intersect in a point.


Let us see what happens if two geodesics at different sides are disjoint. 
Let $\xi_1\xi_2$, $\eta_1\eta_2$ be two such geodesics.
 Continuity argument shows that for each $i=1,2$ there  are  $x_i\in \xi_1\xi_2$, $y_i\in \eta_1\eta_2$
 such that $x_i\eta_i\cap \eta_1\eta_2$  and $y_i\xi_i\cap \xi_1\xi_2$  are geodesic rays. 
 This gives rise to triangles and quadrilaterals that are contained in the 1-skeleton. 
 Hence we are led to consider  triangles and quadrilaterals that are contained in the 1-skeleton.
Each such triangle  or   quadrilateral must bound a topological 
disk which is the union of a finite number of chambers (see Section \ref{triqua1sta}).
 Then Gauss-Bonnet implies  that  there are only a few such triangles and quadrilaterals. 
For instance, there is no such triangle  or   quadrilateral 
when the chamber of $\De$ has at least 5 edges, that is, $k\ge 5$.  Hence in this case  two 
geodesics at different sides 
 must intersect in a point.  The same is true when the chamber is  not a right  triangle.
 The right  triangle case is more involved, and it is indeed possible to have 
 disjoint geodesics that are at different sides. But two 
disjoint geodesics that are at different sides in $\De_2$ cannot occur as the images of two 
 intersecting geodesics contained in a common apartment of $\De_1$. 
For the proof one needs to use some special properties 
 of vertex links, which are generalized polygons.

The proof of Theorem \ref{main}  is relatively simple  when the chambers of $\De_1$ and $\De_2$ have at least 5 edges,
please see Remark \ref{insski}.

The paper is organized as follows.  In Section \ref{pre} we recall the definition of Fuchsian buildings and some facts 
 concerning generalized polygons.  In Section \ref{crossr}
we recall combinatorial cross ratio defined on the ideal boundary of   a   Fuchsian building, and show that
 two important geometric properties can be detected by the behavior of the combinatorial cross ratio. 
In Section \ref{sulli}
 we review  D. Sullivan's approach to Mostow rigidity and   M. Bourdon's    work on  Fuchsian buildings.
In Section \ref{triqua1sta}
we   state the main facts about  triangles and quadrilaterals that are contained in the 1-skeleton.
   The proofs  of these results are contained in 
 the last section (Section \ref{triqua}) to prevent a major disruption of the exposition.
In  Section \ref{gedads}
we give sufficient conditions for two geodesics to be at different sides and study the intersection 
of two geodesics at different sides. 
In  Sections 7-9 we use the tools developed in the previous sections to show that
 the images of all the geodesics  contained in the 1-skeleton of 
  $\De_1$ and through a fixed vertex of $\De_1$ intersect in a   unique  vertex. 
 The proof is divided into three cases: when the chambers are  not  right  triangles (Section \ref{nonrightan});
the chambers are right  triangles but different from $(2,3,8)$ (Section \ref{rightan});
the chambers are $(2,3,8)$ (Section \ref{secsingular}).

\begin{figure}[h]
\noindent
\begin{minipage}[b]{.46\linewidth}
\centering\epsfig{file=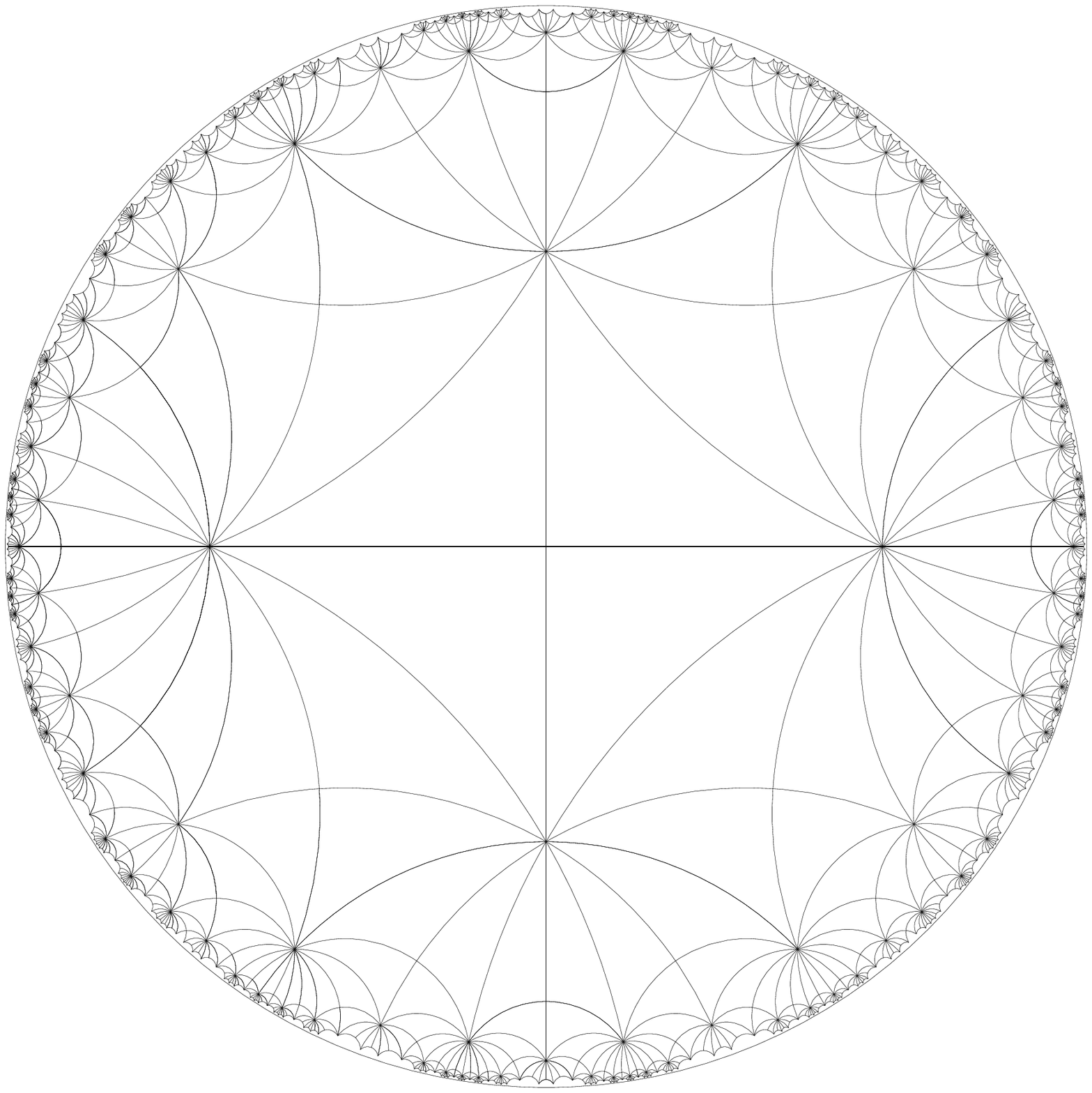,  height=7cm,  width=7cm}\caption{Tessellation of $\H^2$ by $(2,6,8)$\label{te268}}
\end{minipage}\hfill
\begin{minipage}[b]{.46\linewidth}
\centering\epsfig{file=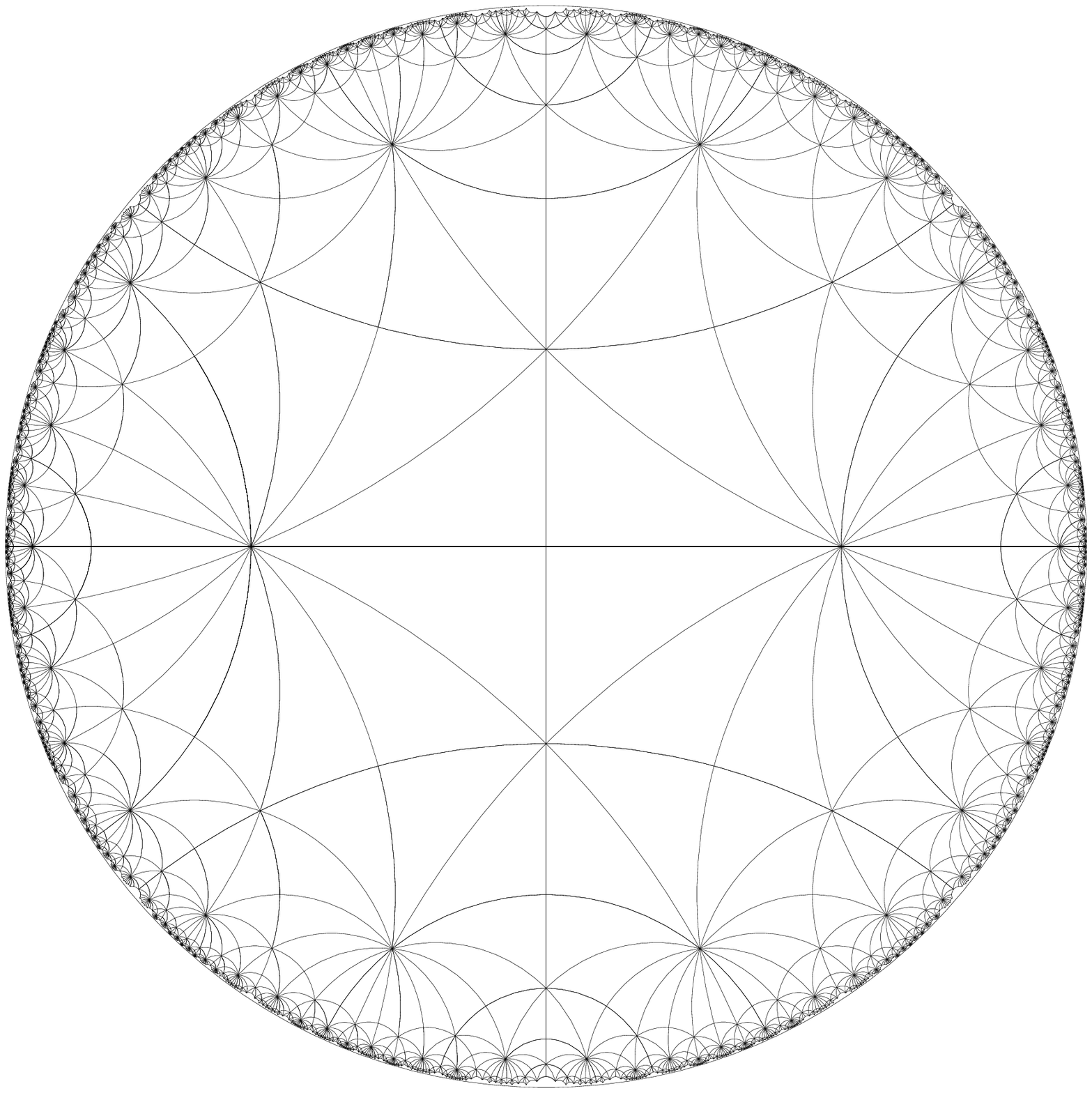,  height=7cm,  width=7cm}\caption{Tessellation of $\H^2$ by $(2,4,8)$\label{te248}}
\end{minipage}
\end{figure}

\begin{figure}[h]
\centering\epsfig{file=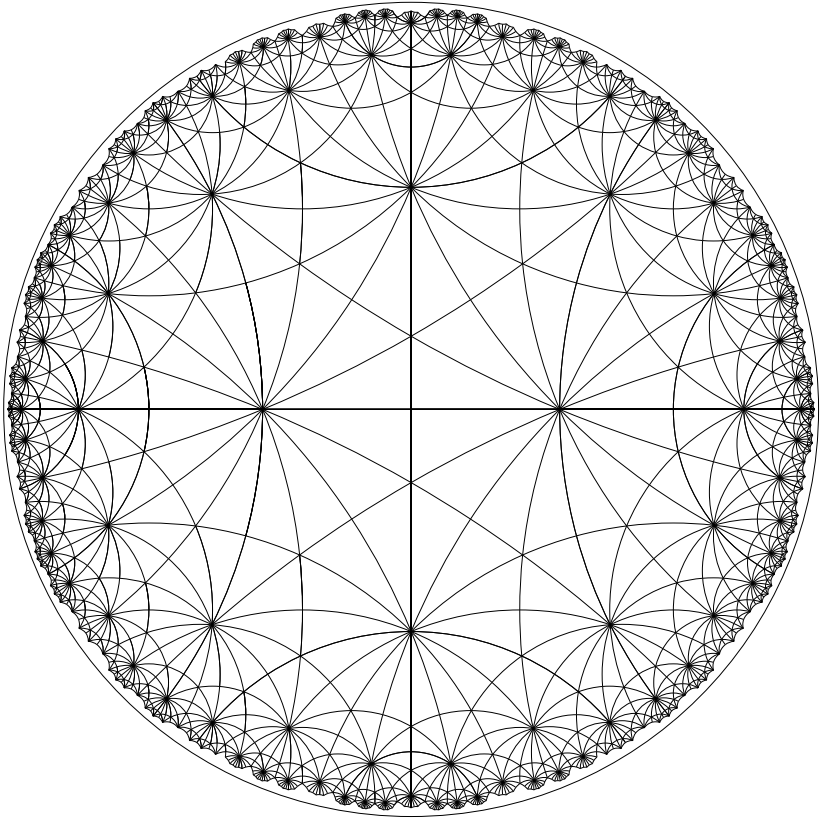,  height=13cm,  width=13cm}
\caption{Tessellation of $\H^2$ by $(2,3,8)$,  created by Martin Deraux}\label{te238}
\end{figure}

\noindent
\textbf{Acknowledgment.} \textit{I would like to thank Bruce  Kleiner for helpful discussions on the subject.
I also appreciate the help from Martin  Deraux, who created those beautiful hyperbolic tessellation pictures in 
 Figures \ref{te268}, \ref{te248}  and \ref{te238}.}


\subsection{Fuchsian buildings}\label{pre}

In this section we  define   Fuchsian buildings and review some facts about 
generalized polygons. Fuchsian buildings were introduced by M. Bourdon (\cite{Bo2}).
 Our presentation  closely follows \cite{Bo2}.

\subsubsection{Convex polygons and reflections groups in $\H^2$} \label{poly}

In this paper, $R$ denotes a compact convex polygon in $\H^2$, whose  angles are of the form 
$\pi/m$, $m\in \N$, $m\ge 2$, and  $W$ the Coxeter group generated by the 
 reflections  about the edges of $R$ (see \cite{M}, theorem IV.H.11).

We label the edges of $R$ cyclically by $\{1\}$,  $\{2\}$,  $\cdots$,   $\{k\}$ and the vertices by 
   $\{1,2\}, \cdots$,    $\{k-1, k\},
 \{k,1\}$
  such that the edges $i$ and $i+1$ intersect at the vertex   $\{i,i+1\}$.   
The angle of $R$ at the vertex $\{i,i+1\}$, $i\in \Z/k\Z$,  is denoted by  
$$\alpha_{i,i+1}=\pi/m_{i,i+1}, \text{with} \;\;  m_{i,i+1}\in \N, \;  m_{i,i+1}\ge 2.$$

It is well-known that (see \cite{M}, theorem IV.H.11)  the images of $R$ under $W$ form  a tessellation 
of $\H^2$  and the quotient $\H^2/W$  equals $R$. It follows that there is a labeling of the edges 
 and vertices of the tessellation  that is $W$-invariant and compatible with that of $R$.
We let $A_R$ be the obtained labeled 2-complex.

Two convex compact polygons $P_1$ and $P_2$ in $\H^2$
are \e{angle-preserving isomorphic} if there is a  combinatorial  isomorphism from 
  $P_1$ to $P_2$ that preserves the angles at the vertices. 
 A  convex compact polygon $P$ is called \e{normal} if it has an inscribed circle that touches all
 its edges. By a result in  A. Beardon's book (see \cite{Bea}, theorem 7.16.2.) each  convex compact polygon
 $P$ is angle-preserving isomorphic to a unique normal polygon $P'$.  Notice that 
 all   compact triangles in  $\H^2$  are normal.  We shall use $(m_1, m_2, m_3)$ to denote the 
 triangle (unique up to isometry) with 
  angles  $\pi/{m_1}$,  $\pi/{m_2}$    and    $\pi/{m_3}$.

\subsubsection{2-dimensional hyperbolic buildings} \label{2hyb}

Let $R$ be a fixed polygon as  defined in Section \ref{poly}.

\b{Def}\label{def2h}
{Let $\Delta$ be a   connected  cellular  2-complex,  whose 
edges  and  vertices  are labeled by $\{1\}$,  $\{2\}$,  $\cdots$,  $\{k\}$ and 
    $\{1,2\}, \cdots,  \{k-1, k\},
 \{k,1\}$ respectively, such that  each 2-cell (called   a  chamber) is isomorphic to $R$ as labeled 2-complexes. 
 $\Delta$ is called a \e{2-dimensional hyperbolic building} if it has a family of 
   subcomplexes (called apartments) isomorphic to $A_R$ (as labeled 2-complexes) with the following properties:\newline
(1) Given any two chambers, there is an apartment containing both;\newline
(2) For any two apartments $A_1$, $A_2$   that share a chamber, there is  an isomorphism
 of labeled 2-complexes 
 $f: A_1\ra A_2$ which pointwise fixes   $A_1\cap A_2$.\newline
$\Delta$ is called a \e{Fuchsian building}  if in addition,  there are integers $q_i\ge 2$, 
  $i=1$,   $2$,   $\cdots, k$
 such that each edge of $\Delta$  labeled by $i$  is contained in exactly $q_i+1$  chambers.}

\end{Def}

Let  $\Delta$ be a 2-dimensional hyperbolic building, $A$ an apartment   and $C\subset A$  a chamber. 
The retraction  onto $A$ centered at $C$ is a label-preserving cellular map
$r_{A,C}:\De\ra A$  defined as follows.
 For any chamber $C'$ of $\De$, choose an apartment $A'$ containing both $C$ and $C'$. 
 Since $A\cap A'$   contains  $C$, there is a   label-preserving
   isomorphism $f:A'\ra A$  that pointwise fixes 
 $A\cap A'\supset C$. Set $r_{A,C}|_{C'}=f|_{C'}$. One checks  that  $r_{A,C}$ is well-defined.

Since a  chamber is label-preserving
isomorphic to $R$, there is a metric
 on each chamber making it isometric to $R$   and  we can  glue the chambers together along the edges by isometries.
 We equip  $\De$ with  the path metric.    Then  $\Delta$ is a 
 $\CAT(-1)$    space  (see \cite{BH} for definition)   and in particular a Gromov hyperbolic space.
   The 
 boundary at infinity $\p \Delta$ is homeomorphic to the Menger curve (see \cite{Ben}).
Each  apartment is a  convex   subset  of   $\Delta$. 
Since the  limits  of apartments  are apartments,  the local finiteness of
 Fuchsian buildings implies that any geodesic in a Fuchsian building is contained in 
 an apartment.

Let  $\Delta$ be  a Fuchsian building.  By  Section \ref{lin} below, 
$m_{i, i+1}\in \{2,3,4,6,8\}$; $q_i=q_{i+1}$ if $m_{i, i+1}=3$,  and $q_i\not=q_{i+1}$ if $m_{i, i+1}=8$.

Let $\Delta$ be a Fuchsian building,  and $A\subset \De$ an apartment of $\De$.  
A \e{wall}  in $A$ is a complete geodesic contained in the 1-skeleton of $A$.
Let $W$ be a wall of $A$, $v\in W$ a vertex  and $e_1, e_2\subset W$ the two edges of $W$ containing $v$.
By assumption the vertex $v$ is labeled by $\{i,i+1\}$ for some $i$.  If $m_{i,i+1}$ is even, then
 $e_1, e_2$ are both labeled by $i$ or both labeled by $i+1$.  If $m_{i,i+1}=3$, then
 one of $e_1, e_2$  is labeled by $i$, the other is labeled by $i+1$ and $q_i=q_{i+1}$.  In any case, the number of chambers 
 containing $e_1$  equals the number of chambers containing $e_2$. It follows that there is  some $i, 1\le i\le k$
  such that each edge in $W$ is contained in exactly $q_i+1$ chambers.  For any wall $W$  and any 
  edge $e\subset W$, we set $l(e)=l(W)=\log q_i$  if $e$ is labeled by $i$.

Let $C, C'$ be   two  chambers. 
A \e{gallery}  from $C$ to $C'$  with length $n$ is a sequence 
${\mathcal G}=(C_0=C, e_1, C_1, \cdots, C_{i-1}, e_{i}, C_{i},  \cdots,  C_{n-1}, e_n, C_n=C')$,  where each $C_i$ is a chamber
 and $e_{i-1}$ is a  common  edge  of    $C_{i-1}$  and    $C_{i}$.  Notice  that  here we allow 
 $C_{i-1}$  and $ C_{i}$  to be the same.  A gallery is minimal if it has the smallest length among  all 
galleries   from $C$ to $C'$.  We define $l({\mathcal{G}})=l(e_1) + \cdots + l(e_n)$.

Let $A$ be an apartment and $C, C'\subset A$ two chambers. A wall $W\subset A$ is said to separate 
 $C$ and $C'$ if the interiors of $C$ and $C'$ lie in different components of $A-W$.
Let ${\mathcal{W}}_A(C,C')$ be the set of walls  in $A$  separating $C$ and $C'$.
 Note  that  each gallery in $A$ from $C$ to $C'$ must cross each wall in ${\mathcal{W}}_A(C,C')$   at least once. 
 On the other hand, it is easy to find a gallery from $C$ to $C'$ 
that crosses
each wall in ${\mathcal{W}}_A(C,C')$   exactly once as follows.  Pick  generic points $x\in \interior(C)$,
 $y\in \interior(C')$  such that $xy$ does not contain any vertices.  Then we obtain a 
gallery   ${\mathcal{G}}_{x,y}$     from $C$ to $C'$ 
 by recording in linear order
the    chambers and edges   intersected  by   $xy$.   Convexity implies that 
 $xy$ intersects each wall in  ${\mathcal{W}}_A(C,C')$ at   exactly   one   point, and  so the gallery 
  crosses 
each wall in ${\mathcal{W}}_A(C,C')$   exactly once.

 Finally let us remark that 
  there are uncountably  many   isomorphism classes of  Fuchsian buildings (\cite{GP}). There are also 
  many constructions of Fuchsian buildings admitting 
   proper and cocompact group actions:    complexes  of groups (\cite{Bo2}, \cite{GP}), 
  ramified coverings of Euclidean
  buildings,  blue-prints of Ronan-Tits  (\cite{RT}), etc.

\subsubsection{Generalized polygons and vertex links} \label{lin}

The reader is referred to 
  \cite{R} (chapter 3) and \cite{T} (section 3) for more details on the material in this subsection.



Generalized polygons are  also called   rank two  spherical buildings.  
  Let $L$ be a 
 connected  graph whose vertices are colored black and white such that 
 the two vertices of each edge always have different colors.
     $L$ 
 is called  a   \e{generalized $m$-gon} ($m\in \N$, $m\ge 2$)
  if    it has   the following properties:\newline
(1) Given any two   edges, there is a circuit with  combinatorial length $2m$ containing both;\newline
(2) For two circuits  $A_1$, $A_2$ of combinatorial  length $2m$  that share an edge, there is an isomorphism
 $f:A_1\ra A_2$ that pointwise fixes $A_1\cap A_2$.\newline
When $L$ is a generalized $m$-gon, an edge is called a chamber and a circuit with combinatorial
length $2m$ is called 
  an apartment.  Two vertices of the same color are often said to have the same type. 
 A theorem of Feit-Higman (\cite{R}, p.30) says that if $L$ is a finite generalized  m-polygon, then $m$ belongs
 to $\{2,3, 4,6,8\}$.

A   generalized $m$-gon  is \e{thick} if each vertex is contained  in at least 3  chambers.  
 Given a   thick   generalized $m$-gon  $L$, there are integers $s,t\ge 2$  such that 
 each black vertex is contained in exactly $s+1$   chambers and each white vertex
is contained in exactly $t+1$   chambers (see p.29 of \cite{R}).    In this case, we say 
  $L$   has parameters $(s,t)$.
 It is known that 
$s=t$  when 
$m$ 
is odd (p.29 of \cite{R}),  and $s\not=t$  when $m=8$ (see \cite{FH}).

Let $L$ be a generalized  $m$-gon.  We put a metric on each  edge such that it is isomeric to
the closed interval with length $\pi/m$.    
 Note that  all the apartments have length $\pi$.  
We equip $L$ with the path metric.
 Then $L$ becomes   a $\CAT(1)$ space (see \cite{BH} for definition). 
 All the apartments are convex in $L$, that is, if $A$ is an apartment and 
 $x,y\in A$ with $d(x,y)<\pi$,   then the geodesic segment $xy$ also lies in $A$.  It follows that if two apartments 
  $A_1$, $A_2$  share a
    chamber, then either $A_1=A_2$ or $A_1\cap A_2$ is a segment.

An injective edge path with combinatorial
length $m$ in a generalized $m$-gon 
is called a half apartment.

\b{Le}\label{halfa}
{Let $L$ be a  thick generalized  $m$-gon, and $A$, $A'$ two apartments in $L$.  Then there  is  
 a finite sequence of apartments $A_0=A$, $A_1, \cdots, A_n=A'$ such that 
$A_i\cap A_{i+1}$ \e{(}$0\le i \le n-1$\e{)} is a
  half apartment.}
\end{Le}

\b{proof}
Let $e\subset A$, $e'\subset A'$ be two chambers. By considering an apartment containing
  $e$ and $e'$, we may assume $A\cap A'$  is  a segment with combinatorial
length $l$, $1\le l\le m$.
 If $l=m$, we are done.  Assume $l<m$.   We claim  in this case there is an apartment $A''$ such that 
 $A\cap A''$ and $A'\cap A''$ have combinatorial
length $>l$.  The lemma clearly follows from the claim.
Let $v_1$, $v_2$ be the endpoints of the segment  $A\cap A'$, $e_1\subset A$ the edge of $A$ 
  incident to $v_1$ but not contained in  $A\cap A'$ and 
  $e_2\subset A'$ the edge of $A'$ 
  incident to $v_2$ but not contained in  $A\cap A'$.   Let $A''$ be an apartment containing
 $e_1$ and $e_2$. 
  The convexity of apartments   implies $A\cap A'\subset A''$, and the claim follows.

\end{proof}

Two vertices of a generalized $m$-gon  are \e{opposite} if the combinatorial
distance between them is $m$.

\b{Prop}\label{oppoti} \e{(\cite{T}, p.57)}
{Let $L$ be a  thick generalized  m-gon  with  $m=3,4$ or $5$.  Then for any apartment $A$, 
 there is a  vertex opposite to all the vertices in $A$ of the same type.}

\end{Prop}

\b{Prop}\label{oppoti2} \e{(\cite{T}, p.55)}
{Let $L$ be a  thick generalized  polygon.  Then for any two vertices 
 $v_1, v_2$ of the same type, 
 there is a  vertex  $v\in L$ opposite to  both $v_1$ and $v_2$. }

\end{Prop}

Let $\Delta$  be a Fuchsian building and $v\in \Delta$ a vertex labeled by $\{i,i+1\}$. 
  The link $\Link(\Delta, v)$ is a graph defined as follows.  The vertex set  of $\Link(\Delta, v)$  is in 
1-to-1 correspondence with  the set of edges of $\Delta$ incident to $v$.  Two vertices of 
$\Link(\Delta, v)$  are connected by an edge if the edges of  $\De$ 
corresponding to the two vertices 
 are contained in a common
   chamber  of  $\De$.    A vertex of $\Link(\Delta, v)$ is black if the corresponding edge in $\De$
 is labeled by $i$, and is white if the corresponding edge in $\De$
 is labeled by $i+1$.  Then 
  it  follows from the definitions that $\Link(\Delta, v)$   is a   finite  thick 
   generalized $m_{i, i+1}$-gon with parameters $(q_i, q_{i+1})$. Consequently, 
   $m_{i, i+1}\in \{2,3,4,6,8\}$; $q_i=q_{i+1}$ if $m_{i, i+1}=3$,  and $q_i\not=q_{i+1}$ if $m_{i, i+1}=8$.


\subsection{Combinatorial cross ratio}\label{crossr}

In this section we  study     combinatorial cross ratio on the ideal boundary of  a 
Fuchsian building, and discuss two geometric properties that can be detected by 
 combinatorial cross ratio. Combinatorial cross ratio  was introduced and studied by M. Bourdon
(\cite{Bo1}, \cite{Bo2}).  Our definitions are slightly different from M. Bourdon's in that we 
only use vertices in the dual graph.

\subsubsection{Dual graph} \label{dualg}

Let $\Delta$ be a Fuchsian building.  The dual graph  ${\mathcal{G}}_{\De}$ of $\Delta$ is 
  a  metric graph
defined as follows.  The vertex set of ${\mathcal{G}}_{\De}$   is in 1-to-1 correspondence with the set
  of chambers of $\Delta$.   Two vertices of ${\mathcal{G}}_{\De}$    are connected by an
  edge   with length $\log q_i$  if the corresponding chambers  in $\De$ share an edge  labeled by $i$.
 We equip ${\mathcal{G}}_{\De}$  with the path metric, and denote the distance  
between $x, y \in {\mathcal{G}}_{\De}$  by  $|x-y|$. 
 ${\mathcal{G}}_{\De}$  is clearly quasi-isometric to $\De$, hence is Gromov hyperbolic and 
  $\p {\mathcal{G}}_{\De}=\p \De$.

Abusing notation,  for any chamber $C$, we shall also use $C$ to denote the   corresponding  vertex
 in ${\mathcal{G}}_{\De}$.

\b{Le}\label{dising}
{Let $C$, $C'$ be two  chambers and $A$ an apartment containing both. 
   Then $|C-C'|=\Sigma_{W\in {\mathcal{W}}_A(C,C')}l(W)$. 
Furthermore, if $\omega'$ is an  edge path in  ${\mathcal{G}}_{\De}$  from $C$ to $C'$
  with length $>|C-C'|$, then $\length(\omega')\ge |C-C'|+\log q_i$
   for some  $i$.}

\end{Le}
\b{proof}
Let $\omega$ be a shortest path in ${\mathcal{G}}_{\De}$  from $C$ to $C'$.  The path 
$\omega$  gives rise to a gallery ${\mathcal{G}}_\omega$ in $\De$ 
    from $C$ to $C'$ with $l({\mathcal{G}}_\omega)=|C-C'|$.   Let ${\mathcal{G}}'_\omega=r_{A,C}({\mathcal{G}}_\omega)$. 
  Since  $r_{A,C}$  is label-preserving,  
     ${\mathcal{G}}'_\omega$   is    a gallery   in $A$   from  $C$ to $C'$ with 
$l({\mathcal{G}}'_\omega)=l({\mathcal{G}}_\omega)$.  On the other hand, for generic points $x\in C$, $y\in C'$, 
  $xy$ does not contain any vertices and ${\mathcal{G}}_{x,y}$ is a minimal gallery   in $A$ from 
  $C$ to  $C'$.  
It follows that $l({\mathcal{G}}_{x,y})\le l({\mathcal{G}}'_\omega)$. 
 It is clear that ${\mathcal{G}}_{x,y}$  determines  an edge path $\omega_{x,y}$ in 
 ${\mathcal{G}}_{\De}$  from $C$ to $C'$  such that  $\length({\omega_{x,y}})=l({\mathcal{G}}_{x,y})$.
Now the minimality of 
$\omega$ implies that $|C-C'|=l({\mathcal{G}}_{x,y})=\Sigma_{W\in {\mathcal{W}}_A(C,C')}l(W)$.

Let  $\omega'$   be   an  edge path in  ${\mathcal{G}}_{\De}$  from $C$ to $C'$
  with length $>|C-C'|$,
  and denote  ${\mathcal{G}}_{A, \omega'}=r_{A,C}({\mathcal{G}}_{\omega'})$.  Then ${\mathcal{G}}_{A, \omega'}$ is a gallery in $A$
  from $C$ to $C'$ with 
$$l({\mathcal{G}}_{A, \omega'})=l({\mathcal{G}}_{\omega'})=\length(\omega')>|C-C'|=\Sigma_{W\in {\mathcal{W}}_A(C,C')}l(W).$$ 
Since ${\mathcal{G}}_{A, \omega'}$  crosses each wall in  ${\mathcal{W}}_A(C,C')$  at least once,
  there is some wall $W$ of  $A$ such that  $l({\mathcal{G}}_{A, \omega'})$ contains the term 
 $l(W)$ in addition to all the terms in $\Sigma_{W\in {\mathcal{W}}_A(C,C')}l(W)$. 
 Hence we have $l({\mathcal{G}}_{A, \omega'})\ge |C-C'|+l(W)$.


\end{proof}

\b{Le}\label{discgeo}
{Let $C$, $C'$ be two chambers in $\De$, and $x\in \interior(C)$, $y\in \interior(C')$.
  Let $C_1=C, \cdots, C_n=C'$ be the sequence of chambers in linear order 
 with $xy\cap \interior(C_i)\not=\phi$.  Then the sequence $C_1=C, \cdots, C_n=C'$ 
  defines a discrete geodesic   in  ${\mathcal{G}}_{\De}$ 
in the following sense: 
  $|C_{i_1}-C_{i_3}|=|C_{i_1}-C_{i_2}|+|C_{i_2}-C_{i_3}|$
for any $1\le i_1\le i_2\le i_3\le n$.}

\end{Le}

\b{proof}
Let $A$ be an apartment containing $C$ and $C'$.  Then the convexity of apartments implies that
 $C_i\subset A$ for all $i$. 
By Lemma \ref{dising}, we have 
$|C_{i_1}-C_{i_3}|=\Sigma_{W\in {{\mathcal{W}}_{A(C_{i_1},C_{i_3})}}}l(W)$, 
$|C_{i_1}-C_{i_2}|=\Sigma_{W\in {{\mathcal{W}}_{A(C_{i_1},C_{i_2})}}}l(W)$  and 
$|C_{i_2}-C_{i_3}|=\Sigma_{W\in {{\mathcal{W}}_{A(C_{i_2},C_{i_3})}}}l(W)$. 
Since the geodesic $xy$ intersects the interiors of $C_{i_1}$, $C_{i_2}$   and  $C_{i_3}$,
we see  ${\mathcal{W}}_{A(C_{i_1},C_{i_3})}$ is the disjoint union of  ${\mathcal{W}}_{A(C_{i_1},C_{i_2})}$
  and  ${\mathcal{W}}_{A(C_{i_2},C_{i_3})}$. The lemma follows.

\end{proof}

\subsubsection{Combinatorial Gromov product} \label{gromovp}

Let $C$ be a vertex of   ${\mathcal{G}}_{\De}$. Recall the Gromov product of $x\in {\mathcal{G}}_{\De}$
  and $y\in {\mathcal{G}}_{\De}$    based at $C$ is defined by:   $\{x|y\}_C=\frac{1}{2} \{|x-C|+|y-C|-|x-y|\}$. 
 For $\xi,\eta\in \p\De=\p {\mathcal{G}}_{\De}$, the combinatorial  
Gromov product of $\xi$ and  $\eta$   based at 
 $C$ is:  
$$\{\xi|\eta\}_C=\sup\; \liminf_{i,j\ra \i} \{x_i|y_j\}_C,$$
 where the supreme is taken over all sequences  of  vertices in  ${\mathcal{G}}_{\De}$   with 
$\{x_i\}\ra \xi$, $\{y_j\}\ra \eta$.

A sequence  $\{C_i\}_{i=1}^{\i}\subset {\mathcal{G}}_{\De}$  of  vertices  is called a  \e{geodesic sequence}
   starting from $C_1$ 
  if there is a geodesic ray  $\sigma$ in $\De$ starting from the interior of $C_1$ such that
   $\{C_i\}$ is the sequence of chambers (in linear order) whose interiors have nonempty intersection with
   $\sigma$.  For $\xi,\eta\in \p\De$, 
 the modified Gromov product of $\xi$ and  $\eta$   based at 
 $C$ is:
$$\{\xi|\eta\}'_C=\sup\; \liminf_{i,j\ra \i} \{C_i|D_j\}_C,$$
 where the supreme is taken over all  geodesic sequences $\{C_i\}\ra \xi$, $\{D_j\}\ra \eta$
  starting from  $C$.
By definition,   $\{\xi|\eta\}'_C\le \{\xi|\eta\}_C$ always holds.

A point $\xi\in\p\De$ is called a \e{singular point} if it is represented by a ray contained in the 1-skeleton.
 Otherwise $\xi\in\p\De$ is called a \e{regular point}. Let $B\subset \p\De$  be the set of regular points in 
  $\p\De$.   For $x\in \De$ and $\xi\in \p\De$,   let  $c_{x,\xi}:  [0,\i)\ra \De$ denote 
 the geodesic ray from $x$ to $\xi$. 
  Since  the angles of a chamber  belong to  $\{\pi/2,\pi/3,\pi/4,\pi/6,\pi/8\}$,
   it is not hard to see that any ray representing a regular point cannot lie in a small neighborhood
 of the 1-skeleton:

\b{Le}\label{1skel}
{There is some $\epsilon_0>0$ depending only  on    $\De$ 
with the following property.
  For any  $\xi\in B$,   $x\in \De$  and any $a>0$, 
  $c_{x,\xi}([a, +\i))\cap (\De-N_{\epsilon_0}(\De^{(1)}))\not=\phi$,
   where  $N_{\epsilon_0}(\De^{(1)})$  is the $\epsilon_0$-neighborhood   of the 1-skeleton 
 $\De^{(1)}$   of $\De$.}
\end{Le}




\qed

\b{Le}\label{modgp}
{Let $\xi\not=\eta\in \p\De$   be regular points.
  If   $\{C_i\}\ra \xi$, $\{D_j\}\ra \eta$  are 
geodesic sequences  starting  from  $C$,  then there is some $i_0\ge 1$ such that 
$\{\xi|\eta\}'_C=\{C_i|D_j\}_C$  for all
   $i,j\ge i_0$.}
\end{Le}

\b{proof}
 By  Lemma  \ref{discgeo}
     the sequences $\{C_i\}$, $\{D_j\}$  are 
  discrete geodesics in ${\mathcal{G}}_{\De}$. 
Triangle inequality implies that $\{C_{i_1}|D_{j_1}\}_C\le \{C_{i_2}|D_{j_2}\}_C$
  for $i_1\le i_2$, $j_1\le j_2$. 
  Suppose $\{C_{i_1}|D_{j_1}\}_C< \{C_{i_2}|D_{j_2}\}_C$
  for  some   $i_1\le i_2$, $j_1\le j_2$.
Let $\omega_1$ be a shortest path in ${\mathcal{G}}_{\De}$
from $C_{i_2}$ to $C_{i_1}$, $\omega_2$ be a shortest path 
from $C_{i_1}$ to $D_{j_1}$  and $\omega_3$ be a shortest path 
from $D_{j_1}$ to $D_{j_2}$.   Then $\length(\omega_3*\omega_2*\omega_1)> |C_{i_2}-D_{j_2}|$.
  Now  Lemma \ref{dising} implies that
 $\length(\omega_3*\omega_2*\omega_1)\ge  |C_{i_2}-D_{j_2}|+\log q_i$  for some $i$.
 Consequently   $\{C_{i_2}|D_{j_2}\}_C\ge  \{C_{i_1}|D_{j_1}\}_C+\log q_i$.  
  It follows that  there is some integer $a\ge 1$ such that 
 $\liminf_{i,j\ra \i} \{C_i|D_j\}_C=\{C_i|D_j\}_C$ for  all  $i, j\ge a$.

Let    $\{C'_i\}\ra \xi$, $\{D'_j\}\ra \eta$  be  two  arbitrary 
geodesic sequences  starting  from  $C$.  The above paragraph   shows that there  is   an  integer $b\ge 1$ such that 
  $\liminf_{i,j\ra \i} \{C'_i|D'_j\}_C=\{C'_i|D'_j\}_C$ for  all  $i, j\ge b$. 
By the definition of geodesic sequences, there are points $x, x', y, y'\in \interior(C)$ such that
 $\{C_i\}$, $\{C'_i\}$,    $\{D_j\}$,    and    $\{D'_j\}$   are  
 the sequences of chambers (in linear order) whose interiors have nonempty intersection with
   $c_{x, \xi}$,   $c_{x',  \xi}$,  $c_{y, \eta}$  and   $c_{y',  \eta}$  respectively. 
We reparameterize the geodesics $c_{x,  \xi}$,   $c_{x', \xi}$,  $c_{y, \eta}$  and   $c_{y',  \eta}$ 
 such that $c_{x, \xi}(t)$,   $c_{x',  \xi}(t)$  lie on the same horosphere centered at $\xi$, and 
 $c_{y,  \eta}(t)$,    $c_{y',  \eta}(t)$  lie on the same horosphere centered at $\eta$. 
 Since $\De$ is a $\CAT(-1)$ space,  $d( c_{x,  \xi}(t),  c_{x',  \xi}(t))\ra 0$  and  
$d(c_{y,  \eta}(t),  c_{y',  \eta}(t))\ra 0$ as $t\ra +\i$.    On the other hand, 
Lemma \ref{1skel}
 implies that  there are arbitrarily large $t$, $t'$  with  $c_{x, \xi}(t)\notin N_{\epsilon_0}(\De^{(1)})$  and 
  $c_{y,  \eta}(t')\notin N_{\epsilon_0}(\De^{(1)})$. 
  It follows that there are $i\ge a$, $j\ge a$ and $i'\ge b$, $j'\ge b$ with 
 $C_i=C'_{i'}$ and $D_j=D'_{j'}$.   Consequently,  
 $\liminf_{i,j\ra \i} \{C'_i|D'_j\}_C=\{C'_{i'}|D'_{j'}\}_C=\{C_i|D_j\}_C=\liminf_{i,j\ra \i} \{C_i|D_j\}_C$
  and the lemma   follows.

\end{proof}

\b{Le}\label{expregr}
{Let $\xi,\eta\in \p\De$   be regular points.   Then   $\{\xi|\eta\}_C=\{\xi|\eta\}'_C$.}
\end{Le}

\b{proof}
Let $\{x_i\}\ra \xi$, $\{y_j\}\ra \eta$  be  two  arbitrary 
 sequences of vertices.   We claim  that  $\{x_i|y_j\}_C=\{\xi|\eta\}'_C$
for sufficiently large 
$i$ and $j$.  The lemma follows from the claim.

Fix some $x\in \interior(C)$ and let $\{C_k\}$ and $\{D_l\}$  be the geodesic sequences 
 corresponding to the rays $x\xi$ and $x\eta$ respectively. By Lemma \ref{modgp},
 there is some $i_0\ge 1$ such that 
$\{\xi|\eta\}'_C=\{C_k|D_l\}_C$  for all
   $k, l\ge i_0$.  Since $\xi$ is regular and the asymptotic distance between
 $x\xi$ and $\eta\xi$ is 0, Lemma \ref{1skel} implies  that there is some $k\ge i_0$ such that
   both $x\xi\cap (C_k-N_{\frac{\epsilon_0}{2}}(\De^{(1)}))$ and 
$\xi\eta\cap (C_k-N_{\frac{\epsilon_0}{2}}(\De^{(1)}))$   are nonempty, where $\epsilon_0$ is as in 
Lemma \ref{1skel}.  Similarly there is some $l\ge i_0$ such that
   both $x\eta\cap (D_l-N_{\frac{\epsilon_0}{2}}(\De^{(1)}))$ and 
$\xi\eta\cap (D_l-N_{\frac{\epsilon_0}{2}}(\De^{(1)}))$   are nonempty. 
 It follows that for  sufficiently large $i, j$ (so that $x_i$ is sufficiently  close to $\xi$ and 
 $y_j$ is  sufficiently  close to $\eta$),  and any $x'\in x_i$, $y'\in y_j$,  
we have   $xx'\cap \interior(C_k)\not=\phi$, 
 $x'y'\cap \interior(C_k)\not=\phi$,
 $xy'\cap \interior(D_l)\not=\phi$,
and  $x'y'\cap \interior(D_l)\not=\phi$.   Now Lemma \ref{discgeo}  implies  that 
$\{x_i|y_j\}_C=\{C_k|D_l\}_C=\{\xi|\eta\}'_C$.

\end{proof}

A similar argument shows that the combinatorial 
Gromov product is locally constant on  the set of regular points:
   
\b{Le}\label{localc}
{Let $\xi,\eta\in \p\De$   be regular points.  Then 
there are neighborhoods $U\owns \xi$, $V\owns \eta$ such that
  $\{\xi'|\eta'\}_C=\{\xi|\eta\}_C$  for all $\xi'\in B\cap U$, $\eta'\in B\cap V$.
}

\end{Le}

\qed

The following  result  follows easily from 
Lemmas \ref{expregr},  \ref{modgp} and \ref{discgeo}.

\b{Le}\label{simgp}
{Let  $C$ be  a chamber  and 
$\xi_1, \xi_2\in B$.    
  If $\xi_1\xi_2\cap \interior(C)\not=\phi$, then $\{\xi_1|\xi_2\}_C=\{\xi_1|\xi_2\}'_C=0$.}

\end{Le}

\subsubsection{Combinatorial cross ratio} \label{dcror}

For  a regular point  $\xi\in B$,  we can define the  Busemann function 
  $B_\xi$  as follows. 
 Let $\xi\in B$,   and  
     $C, D$ be chambers.    Let 
   $\{C_i\}\ra  \xi$ be   a geodesic sequence  starting  from  some chamber,  and set 
 $$ B_\xi(C,D)=\lim_{i\ra +\i}(|D-C_i|-|C-C_i|).$$ 
 Note  that    the triangle inequality implies that 
 the limit exists. The arguments in the proofs   of Lemma \ref{expregr}
  and Lemma \ref{modgp}
 show that   $B_\xi(C,D)$   is  independent of the choice of   $\{C_i\}$ 
  and there is some $i_0\ge 1$ such
  that $ B_\xi(C,D)=|D-C_i|-|C-C_i|$  for  all $i\ge i_0$. 
 It follows that $ B_\xi(C,E)= B_\xi(C,D)+ B_\xi(D,E)$  for any three
  chambers  $C, D, E$.  It is also easy to see that 
$B_\xi(C,D)$  is locally constant as a function of regular points $\xi$: 
$B_{\xi'}(C,D)=B_\xi(C,D)$  for all regular points $\xi'$
  sufficiently close to  $\xi$.

Next we look at how the combinatorial Gromov product changes when the base point changes.

\b{Le}\label{changba}
{Let $E, E'$ be two chambers and $\xi, \eta\in B$.  Then
$\{\xi|\eta\}_{E'}=\{\xi|\eta\}_{E}+\frac{1}{2}B_\xi(E,E')+\frac{1}{2}B_\eta(E,E')$.}

\end{Le}

\b{proof}
Let $\{C_i\}$,   $\{D_j\}$  be two geodesic sequences from $E$ to  $\xi$, $\eta$ respectively,
  and $\{C'_i\}$,   $\{D'_j\}$  be two geodesic sequences from $E'$ to  $\xi$, $\eta$ respectively.
   Lemma \ref{modgp}
 and the first paragraph of this  subsection show that there is some $i_0\ge 1$ such that
 $\{\xi|\eta\}_{E}=\{C_i|D_j\}_E$,  $\{\xi|\eta\}_{E'}=\{C'_k|D'_l\}_{E'}$,
  $B_\xi(E,E')=|E'-C_i|-|E-C_i|$ and  $B_\eta(E,E')=|E'-D_j|-|E-D_j|$
  for all  $i,j,k,l\ge i_0$.  
  On the other hand, since $\xi$ and $\eta$ are regular points, Lemma \ref{1skel}
  implies that there are $i_1, j_1, k_1,l_1\ge i_0$ with
  $C_{i_1}=C'_{k_1}$ and $D_{j_1}=D'_{l_1}$.  Now the lemma follows.

\end{proof}

 Lemma \ref{changba}  ensures   that  the following definition is well-defined.

\b{Def}\label{defcra}
{Let $\xi_1,\xi_2, \eta_1, \eta_2\in B$ be regular points  that are pairwise distinct.
  The combinatorial cross ratio $\{\xi_1\xi_2\eta_1\eta_2\}$ is defined by:
$$ \{\xi_1\xi_2\eta_1\eta_2\}=-\{\xi_1|\eta_1\}_C-\{\xi_2|\eta_2\}_C+\{\xi_1|\eta_2\}_C+\{\xi_2|\eta_1\}_C,$$
  where $C$ is any chamber of $\De$.}

\end{Def}

Denote $B^4=\{(\xi_1, \xi_2, \xi_3, \xi_4):  \xi_i\in B,\; \xi_i\not=\xi_j \;\; \text{for}\;\; i\not=j\}$.  
Lemma \ref{localc} implies  that 
    combinatorial cross ratio is locally constant: 
$\{\xi'_1\xi'_2\eta'_1\eta'_2\}=\{\xi_1\xi_2\eta_1\eta_2\}$ for all  
$(\xi'_1,  \xi'_2,   \eta'_1,  \eta'_2)\in B^4$  sufficiently close to 
$(\xi_1,  \xi_2,   \eta_1,  \eta_2)\in B^4$.

We record the following simple lemma  for later use.

\b{Le}\label{lbusem}
{Let $e\subset \De$ be an edge labeled by $i$,   $x\in \interior(e)$
  and $C_1,  C_2$ two chambers containing $e$.   
Let $\xi\in B$ be a regular point  and $C$ the chamber that contains the initial segment of 
 $x\xi$.  Then\newline
\e{(1)}  $B_\xi(C_1, C_2)=0$ if $C\not=C_1,  C_2$;\newline
\e{(2)}  $B_\xi(C_1, C_2)=\log q_i$ if $C=C_1$;\newline
\e{(3)}  $B_\xi(C_1, C_2)=-\log q_i$  if  $C=C_2$.}

\end{Le}

\b{proof}
Let $C'\not=C$ be a chamber containing $e$.  It is clear that  $c_{x, \xi}: [0, \i)\ra \De$ 
   can be extended to a geodesic $c:(-\epsilon, \i)\ra \De$ for some $\epsilon>0$ such that 
$c(-\epsilon, 0)\subset \interior(C')$.   Now the lemma follows from Lemma \ref{discgeo}.

\end{proof}

\subsubsection{Geometric properties detected by combinatorial  cross ratio} \label{geopdcr}

In this section we discuss some geometric properties of a Fuchsian building that can be 
 detected by  the combinatorial cross ratio. The results here are crucial  to the proof of 
   Theorem \ref{main}.

 Let  $\xi,\eta\in \p\De$.  The following result says  that one can decide 
 whether $\xi\eta$ lies in the 1-skeleton  by  looking at the behavior of the combinatorial 
cross ratio.  

\b{Le}\label{geo1ske}
{Let $\xi,\eta\in \p\De$.  \newline
\e{(1)} Suppose  $\xi\eta\not\subset  \De^{(1)}$.  Then    there are  
neighborhoods $U_0\owns \xi$,  $V_0\owns \eta$  in  $\p\De$  
   such that
 $\{\xi_1\xi_2\eta_1\eta_2\}=0$  for
 all  $(\xi_1, \xi_2, \eta_1, \eta_2)\in  B^4\cap  (U_0\times U_0\times V_0\times V_0)$;\newline
\e{(2)} Suppose $\xi\eta\subset  \De^{(1)}$.  
  Then 
given any 
neighborhoods $U_0\owns \xi$,  $V_0\owns \eta$  in  $\p\De$,     
  there   exist some integer $i$, $1\le i\le k$,   open subsets $U, V\subset \p\De$  
and open subsets $W_1, W_2\subset U\times U\times V\times V$ with the following properties:
$\xi\in U\subset U_0$,  $\eta\in V\subset V_0$;  
$\{\xi_1\xi_2\eta_1\eta_2\}$   is    an integral  multiple of $\frac{1}{2}{\log {q_i}}$  
 for all   $(\xi_1, \xi_2, \eta_1, \eta_2)\in  B^4\cap  (U\times U\times V\times V)$; 
$\{\xi_1\xi_2\eta_1\eta_2\}=-\frac{1}{2}{\log  {q_i}}$  for all 
$(\xi_1, \xi_2, \eta_1, \eta_2)\in B^4\cap W_1$,
  and  $\{{\xi_1}{\xi_2}{\eta_1}{\eta_2}\}=0$  for all
$({\xi_1}, {\xi_2}, {\eta_1}, {\eta_2})\in   B^4\cap W_2$.}

\end{Le}

\b{proof} (1).  
 Suppose $\xi\eta$ is not contained in  $\De^{(1)}$. Then $\xi\eta$ intersects the interior of
 some chamber $C$.  There are neighborhoods $U_0\owns \xi$ and $V_0\owns \eta$ such that
 $\xi'\eta'$ meets the interior of $C$ for all  $\xi'\in U_0$, $\eta'\in V_0$. 
    Then   Lemma \ref{simgp}
 implies that $\{\xi'|\eta'\}_C=0$ for all $\xi'\in B\cap U_0$,  
 $\eta'\in B\cap V_0$.  Now (1) follows. 


(2).  
Now suppose $\xi\eta$ is contained in $\De^{(1)}$.  Let $e\subset \xi\eta$ be an edge, which is labeled by some $i$.
Given any 
neighborhoods $U_0\owns \xi$,  $V_0\owns \eta$  in  $\p\De$,  we choose 
    open subsets $U, V\subset \p\De$  with  
   $\xi\in U\subset U_0$,  $\eta\in V\subset V_0$  such that for any
 $\xi'\in U$ and  $\eta'\in V$, the geodesic $\xi'\eta'$ 
 either  contains  $e$ or 
intersects the interior of 
  some chamber containing $e$. Let $C$ be a fixed chamber containing $e$,  and $\xi'\in B\cap U$, 
 $\eta'\in B\cap V$. Then there is some  chamber $C'$ containing $e$ with $\xi'\eta'\cap \interior(C')\not=\phi$.
  Lemma \ref{simgp}
implies that $\{\xi'|\eta'\}_{C'}=0$.  On the other hand, Lemma \ref{lbusem}
 shows that $B_{\xi'}(C', C)$ and  $B_{\eta'}(C', C)$  are  integral multiples of 
  $\log q_i$.  It follows from Lemma \ref{changba}
   that $\{\xi'|\eta'\}_{C}$  is an integral multiple of 
  $\frac{1}{2}{\log q_i}$.  Consequently  
$\{\xi_1\xi_2\eta_1\eta_2\}$   is    an integral  multiple of $\frac{1}{2}{\log {q_i}}$  
 for all   $(\xi_1, \xi_2, \eta_1, \eta_2)\in  B^4\cap  (U\times U\times V\times V)$.

We claim there   exist 
$(\xi_1, \xi_2, \eta_1, \eta_2), 
(\bar{\xi_1}, \bar{\xi_2}, \bar{\eta_1}, \bar{\eta_2})\in   B^4\cap (U\times U\times V\times V)$ 
  such that $\{\xi_1\xi_2\eta_1\eta_2\}=-\frac{1}{2}{\log  {q_i}}$  and
  $\{\bar{\xi_1}\bar{\xi_2}\bar{\eta_1}\bar{\eta_2}\}=0$.
  Since the combinatorial cross ratio is locally constant, there are 
open subsets $W_1, W_2\subset U\times U\times V\times V$
   with $(\xi_1, \xi_2,  \eta_1,   \eta_2)\in W_1$,
 $(\bar{\xi_1},  \bar{\xi_2},  \bar{\eta_1},  \bar{\eta_2})\in W_2$  such that 
 the combinatorial cross ratio is  $-\frac{1}{2}{\log  {q_i}}$  on $B^4\cap W_1$  and 
is  0 on $B^4\cap W_2$.    We next prove the claim.

Let $A$ be an apartment containing $\xi\eta$.
 Denote by $C_1$, $C_2$ the two 
 chambers in $A$ containing $e$, and $C_3\not=C_1, C_2$    some other  chamber containing $e$.  
  Let $m$ be the midpoint of $e$.  
We choose  
  $\xi'_1\in U\cap \p A$,     $\eta_1',   \eta_2'\in  V\cap \p A$    and $\xi'_2\in U$  
   with the following properties:\newline
(a) $C_1$ contains the initial segments of $m\xi'_1$ and $m\eta'_1$, 
 $C_2$ contains the initial segment of $m\eta'_2$,
   and   $C_3$  contains the initial segment   of $m\xi'_2$;\newline
(b) $\angle_{m}(\xi'_i, \xi)=\angle_m(\eta'_j, \eta)$  for all $i,j=1,2$;    \newline
(c)  $\xi'_1\eta'_1\cap \interior(C_1)\not=\phi$.\newline
  Notice (b)  implies   that $m\in \xi'_1\eta'_2,  \xi'_2\eta'_1,    \xi'_2\eta'_2$.  
We can choose regular points $\xi_1, \xi_2\in U$, $\eta_1, \eta_2\in V$ 
  sufficiently close to $\xi'_1, \xi'_2,  \eta'_1, \eta'_2$   respectively  
  such    that  the following   conditions  are satisfied:\newline
(I) $\xi_1\eta_j\cap \interior(C_1)\not=\phi$ ($j=1,2$),  
  $\xi_2\eta_1\cap \interior(C_1)\not=\phi$;   \newline
(II) $\xi_2\eta_2\cap \interior(C_k)\not=\phi$ for $k=2,3$.\newline
   By   Lemma \ref{simgp}   
  $\{\xi_1|\eta_j\}_{C_1}=\{\xi_2|\eta_1\}_{C_1}=\{\xi_2|\eta_2\}_{C_2}=0$.
 Notice  (II) implies that  $\xi_2\eta_2\cap \interior(e)\not=\phi$  and  
 $\xi_2\eta_2\cap \interior(C_1)=\phi$  since $C_2\cup C_3$ 
 is convex in $\De$ and $e$ separates $C_2$  from  $C_3$.   
   It follows from  
  Lemma \ref{lbusem}   that 
$B_{\xi_2}(C_2, C_1)=0$ and $B_{\eta_2}(C_2, C_1)=\log q_i$. 
 Now Lemma \ref{changba}
 implies that $\{\xi_2|\eta_2\}_{C_1}=\frac{1}{2}{\log q_i}$.
   Consequently            $\{\xi_1\xi_2\eta_1\eta_2\}=-\frac{1}{2}{\log q_i}$.
Finally if one chooses $\bar{\xi_1}, \bar{\xi_2}\in U$ sufficiently close to 
 $\xi'_1$,  and also chooses  $\bar{\eta_1}, \bar{\eta_2}\in V$ sufficiently close to 
 $\eta'_1$,  then we  have  $\{\bar{\xi_1}\bar{\xi_2}\bar{\eta_1}\bar{\eta_2}\}=0$. 
 The proof  is now complete.

\end{proof}

Let $\eta\in \p\De$ be a singular point. Then $\eta$ is represented 
by a ray $\sigma:[0,\i)\ra \De$ contained in $\De^{(1)}$.  Let $\xi_1, \xi_2\in \De\cup \p\De-\{\xi\}$ be such that
 $\xi_1\eta\cap \sigma=\xi_2\eta\cap \sigma=\phi$. 
  We parameterize  $\xi_i\eta$ ($i=1,2$) 
by $\sigma_i:  I_i \ra \De$  such that $\sigma_i(t)$ and $\sigma(t)$ lie on the same horosphere
 centered at $\eta$,   where $I_i\subset \R$ is  an interval.  
 Since $\De$ is $\CAT(-1)$, we have $d(\sigma_i(t), \sigma(t))\ra 0$ as $t\ra \i$. 
 Fix a very large $t_0$ such that $\sigma(t_0)$ is the midpoint of  some  edge  $e\subset \sigma$. Then
  $ \sigma_i(t_0)$ must lie in the interior of some chamber  $C_i$ containing $e$.  We say 
\e{$\xi_1$ and $\xi_2$ lie at different sides of 
 $\eta$} if  $C_1\not=C_2$,  and say 
 \e{$\xi_1$ and $\xi_2$ lie  on the same side  of    $\eta$}
 if   $C_1=C_2$.  It is not hard to see that this definition is well-defined.


The following result says  that, for a singular point $\eta$ represented by a ray 
$\sigma:[0,\i)\ra \De$ contained in $\De^{(1)}$,  and 
  $\xi_1, \xi_2\in \p\De$ with 
  $\xi_1\eta\cap \sigma=\xi_2\eta\cap \sigma=\phi$, whether  $\xi_1$, $\xi_2$ lie at different sides 
 of $\eta$  can be detected by the combinatorial cross ratio.

\b{Le}\label{difes}
{Let $\eta\in \p\De$ be a singular point   represented 
by a ray $\sigma:[0,\i)\ra \De$ contained in $\De^{(1)}$,  and 
$\xi_1, \xi_2\in \p\De$ with 
  $\xi_1\eta\cap \sigma=\xi_2\eta\cap \sigma=\phi$.  \newline
\e{(1)}  If  $\xi_1$ and $\xi_2$ lie  on the same side  of    $\eta$,  then there are pairwise disjoint open
 subsets  $U_1, U_2, V_1, V_2\subset \p\De$  that  
satisfy the following conditions:  $\xi_1\in U_1$, $\xi_2\in U_2$, $\eta\in V_2$;  
the  combinatorial 
cross ratio  is constant on $B^4\cap (U_1\times U_2\times V_1\times V_2)$.\newline
\e{(2)} If  $\xi_1$ and $\xi_2$ lie  on   different  sides  of    $\eta$,
then for any pairwise disjoint open
 subsets  $U_1, U_2, V_1, V_2\subset \p\De$  satisfying  $\xi_1\in U_1$, $\xi_2\in U_2$, $\eta\in V_2$,
there are open subsets $W_1, W_2\subset U_1\times U_2\times V_1\times V_2$  
  and constants $c_1\not=c_2$ 
 such that:   $\{\xi_1'\xi_2'\eta'_1\eta'_2\}=c_1$  for 
all $(\xi_1',  \xi_2',  \eta'_1,  \eta'_2)\in B^4\cap W_1$  and 
 $\{\xi_1'\xi_2'\eta'_1\eta'_2\}=c_2$  for 
all $(\xi_1',  \xi_2',  \eta'_1,  \eta'_2)\in B^4\cap W_2$.}

\end{Le}

\b{proof}
(1). Suppose $\xi_1$, $\xi_2$ lie at the same side of $\eta$. Then  there is a chamber $C$ that contains 
 an edge $e\subset \sigma$ such that both $\xi_1\eta$ and    $\xi_2\eta$  intersect the interior of $C$. 
Then there are  pairwise disjoint 
neighborhoods $U_i\owns \xi_i$($i=1,2$), $V\owns \eta$ such that  for any $\xi_1'\in  B\cap U_1$, 
$\xi_2'\in B\cap U_2$ and $\eta'\in B\cap V$, both $\xi_1'\eta'$ and    $\xi_2'\eta'$  intersect the interior of $C$.
 Lemma \ref{simgp}  implies 
  that $\{\xi'_1|\eta'\}_C=\{\xi'_2|\eta'\}_C=0$.
Now choose   disjoint open subsets $V_1, V_2\subset V$ such that $\eta\in V_2$.   
   Then it follows from the definition   that
$\{\xi_1'\xi_2'\eta'_1\eta'_2\}=0$  for all  
$(\xi_1',  \xi_2',  \eta'_1,  \eta'_2)\in  B^4\cap  (U_1\times U_2\times V_1\times V_2)$.

(2). Now suppose $\xi_1$, $\xi_2$ lie at different  sides of $\eta$.  
We claim for any pairwise disjoint open
 subsets  $U_1, U_2, V_1, V_2\subset \p\De$  satisfying  $\xi_1\in U_1$, $\xi_2\in U_2$, $\eta\in V_2$,
there exist  $\xi'_1\in B\cap U_1, \xi'_2\in B\cap  U_2, \eta'_1\in B\cap V_1$ 
     and $\eta'_2, \eta''_2\in B\cap V_2$ 
such that $\{\xi_1'\xi_2'\eta'_1\eta'_2\}\not=\{\xi_1'\xi_2'\eta'_1\eta''_2\}$.
Let $c_1=\{\xi_1'\xi_2'\eta'_1\eta'_2\}$,  $c_2=\{\xi_1'\xi_2'\eta'_1\eta''_2\}$.
  Since the combinatorial 
cross ratio  is  locally constant, there are open subsets $W_1, W_2\subset U_1\times U_2\times V_1\times V_2$
  satisfying  $(\xi_1',  \xi_2',  \eta'_1,  \eta'_2)\in W_1$,
$(\xi_1',   \xi_2',   \eta'_1,   \eta''_2)\in W_2$
  such that the combinatorial 
cross ratio  takes  the constant  value   $c_i$ ($i=1,2$) on $B^4\cap W_i$.   
 Next we prove the claim.

Since $\xi_1$, $\xi_2$ lie at different  sides of $\eta$, 
  there are chambers $C_1\not=C_2$  both  containing  an edge $e\subset \sigma$
  such that $\xi_i\eta\cap \interior(C_i)\not=\phi$  ($i=1,2$).
  By shrinking the neighborhoods  $V_2$, $U_1$, $U_2$ we may assume
   $\xi'_i\eta'\cap \interior(C_i)\not=\phi$ ($i=1,2$) for 
 all  $\eta'\in B\cap  V_2$, $\xi'_i\in B\cap  U_i$.  
Lemma \ref{simgp}  implies 
    that $\{\xi_i'|\eta'\}_{C_i}=0$. 
 We fix  some   $\xi'_i\in B\cap U_i$($i=1,2$) and $\eta'_1\in V_1$.
  By the definition of combinatorial
cross ratio  we only need to find  $\eta'_2, \eta''_2\in  B\cap V_2$ such that 
 $\{\xi'_2|\eta'_2\}_{C_1}\not= \{\xi'_2|\eta''_2\}_{C_1}$. 
 By Lemma \ref{changba}, 
$\{\xi'_2|\eta'\}_{C_1}=\{\xi'_2|\eta'\}_{C_2}+\frac{1}{2}B_{\xi'_2}(C_2,C_1)+\frac{1}{2}B_{\eta'}(C_2,C_1)$.
  $B_{\xi'_2}(C_2,C_1)$  is independent of $\eta'$ and we have observed that 
  $\{\xi'_2|\eta'\}_{C_2}=0$.  
Let $m$  be the  midpoint of $e$. 
For $\eta'\in B\cap V_2$, the initial segment of $m\eta'$ could lie in $C_1$, $C_2$ or some other
 chamber $C_3$ containing $e$.  Lemma \ref{lbusem}
 shows that the values for $B_{\eta'}(C_2,C_1)$ are different in these three cases. 
 Consequently  the values  of  $\{\xi_1'\xi_2'\eta'_1\eta'\}$ are also different,  and the proof is complete.

\end{proof}

\subsection{Sullivan's approach   and Bourdon's  work on  Fuchsian buildings}\label{sulli}
    
In this section we review Bourdon's  work on  Fuchsian buildings
  and  Sullivan's approach  to  Mostow rigidity.

\subsubsection{Combinatorial metrics on the boundary} \label{confm}

Combinatorial metrics   were introduced and studied by M. Bourdon (\cite{Bo1}, \cite{Bo2}). 
 Here we recall the definition and some facts about   combinatorial metrics.
Let $A_R$  be the labeled 2-complex of Section \ref{poly} and ${\mathcal{G}}_{A_R}$ its dual graph. 
  An edge of ${\mathcal{G}}_{A_R}$ has length $\log q_i$ if the corresponding chambers share
 an edge  labeled by $i$. Let $|\cdot|$ be the induced path metric.  
  For any integer $n\ge 1$,  let  $a(n)$ be the number of vertices that are at distance at 
 most $n$ from   a fixed vertex  $C$  of  ${\mathcal{G}}_{A_R}$.
 The growth rate of 
${\mathcal{G}}_{A_R}$  is defined by: 
$$\tau={\lim \sup}_{n\ra \i}\frac{1}{n}\log a(n).$$

Let $C$ be a chamber of $\De$. For $\xi\in \p\De$ and $r>0$,  we introduce the following subset of 
 $\p\De$:
$$B_C(\xi,r)=\{\eta\in \p\De:  e^{-\tau\{\xi|\eta\}_C}\le r\}.$$
 For a continuous path $\gamma$ in $\p\De$, we define  
$$l(\gamma)=\lim_{r\ra 0} \inf \{\Sigma_i r_i\},$$
 where the infimum is taken over all finite coverings $\{B_C(\xi_i, r_i)\}$ of  $\gamma$, with $\xi_i\in \gamma$
     and $r_i\le r$.    Finally, for  $\xi,\eta\in \p\De$   we  set 
$$\delta_C(\xi,\eta)=\inf {l(\gamma)},$$
 where $\gamma$ varies over all continuous path in $\p\De$  from $\xi$ to $\eta$.

It is shown by M. Bourdon (\cite{Bo1}, 3.1.4)  that $\delta_C$ defines a metric on $\p\De$ and has the following property: 
there is a constant $\lambda\ge 1$ such that 
$$\frac{1}{\lambda}e^{-\tau\{\xi|\eta\}_C}\le \delta_C(\xi,\eta)\le \lambda e^{-\tau\{\xi|\eta\}_C}$$  for
  all $\xi,\eta\in \p\De$.

We shall show that for any two chambers $C$, $D$, the combinatorial metrics $\delta_C$ and 
$\delta_D$ are conformally equivalent.  We first recall the definition of conformal 
 maps.

 Let $f:X\ra Y$ be a homeomorphism between metric spaces. 
The map $f$ is \e{quasi-symmetric}  if there exists a homeomorphism  
$\phi:[0,\i)\ra [0,\i)$ so that 
$$d_X(x,a)\le t d_X(x,b)\Rightarrow d_Y(f(x), f(a))\le \phi(t) d_Y(f(x), f(b))$$
for all $x, a,b\in X$ and all $t\in [0,\i)$.

Let $f:X\ra Y$ be a homeomorphism between metric spaces.
For any $x\in X$ and $r>0$, let
$$L_f(x,r)=\sup \{d_Y(f(x),f(x')): d_X(x, x')\le r\},$$
$$l_f(x,r)=\inf \{d_Y(f(x),f(x')): d_X(x, x')\ge r\},$$
$$L_f(x)=\lim \sup_{r\ra 0} \frac{L_f(x,r)}{r},$$
$$l_f(x)=\lim \inf_{r\ra 0} \frac{l_f(x,r)}{r}.$$
Assume $X$ and $Y$ have finite Hausdorff dimensions. Denote by $H_X$ and $H_Y$ their Hausdorff dimensions  and by $\mathcal{H}_X$
  and $\mathcal{H}_Y$  their Hausdorff  measures
(see \cite{F} for definitions).  We say that $f$ is  \e{conformal} if $f$ is quasi-symmetric and satisfies\newline
(i)  $L_f(x)=l_f(x)\in (0,\i)$ for $\mathcal{H}_X$-almost every $x\in X$;\newline
(ii)  $L_{f^{-1}}(y)=l_{f^{-1}}(y)\in (0,\i)$ for $\mathcal{H}_Y$-almost every $y\in Y$.

Recall (see the proof of Lemma 2.2.7 in \cite{Bo2}) that for any chamber $C$, the set 
  $B$ has  full $\mathcal{H}_{\delta_C}$-measure in 
$\p\De$.

\b{Le}\label{conbase}
{Let $C$, $D$ be two chambers of $\De$.   Then
 $\delta_C$ and  $\delta_D$  are conformally equivalent, that is,   the identity map
 $\id:  (\p\De, \delta_C)\ra (\p\De, \delta_D)$  is a conformal map.  Furthermore,  for any regular point $\xi$,
 there is a neighborhood $V$ of $\xi$ in $\p\De$ such that 
$\delta_D(\xi,\eta)=e^{-\tau B_\xi(C,D)}\delta_C(\xi,\eta)$   for all $\eta\in V$.}

\end{Le}

\b{proof}
Triangle inequality implies $ -|C-D|\le   \{\xi|\eta\}_C-\{\xi|\eta\}_D\le |C-D|$.
Since for  any  chamber $C'$, the combinatorial distance  $\delta_{C'}(\xi,\eta)$ is comparable with
$e^{-\tau\{\xi|\eta\}_{C'}}$,  we see  $\id:  (\p\De, \delta_C)\ra (\p\De, \delta_D)$ 
 is bi-Lipschitz, in particular, it is quasi-symmetric.

We first assume $C$ and $D$ share an edge $e$. 
Then $e$ is  labeled by  some $i$.
 Denote by  $m$  the midpoint of $e$.
Let  $\xi\in \p\De$ be a regular point.  
The initial segment of
   $m\xi$ lies in $C$, $D$ or some chamber $E\not=C, D$
 containing $e$. 
 Suppose  the initial segment of
   $m\xi$ lies in $C$.  Then  there is a neighborhood $U\subset \De\cup \p\De$ of $\xi$ such that 
   the initial segment of
   $m\eta$ lies in $C$  for all $\eta \in U$.  All these geodesic rays (or segments) can be extended into 
 $D$. It follows that $|D-C'|=\log q_i+|C-C'|$ for all chambers $C'\subset U$.  Now it is clear that 
 $\{\eta_1|\eta_2\}_D=\log q_i+\{\eta_1|\eta_2\}_C$ for all $\eta_1, \eta_2\in U\cap \p\De$.
 It implies that for $\xi'\in U\cap \p\De$ and sufficiently small $r>0$ we have
 $B_D(\xi',r)=B_C(\xi', re^{\tau\log q_i})$.  The definition of $\delta_C$ now implies  that 
 $\delta_D(\xi, \eta)=e^{-\tau\log q_i}\delta_C(\xi, \eta)$ for all $\eta\in U\cap \p\De$. Note
$\log q_i=B_\xi(C,D)$. Hence 
$\delta_D(\xi, \eta)=e^{-\tau B_\xi(C,D)}\delta_C(\xi, \eta)$ for all $\eta\in U\cap \p\De$, 
  which shows $L_f(\xi)=l_f(\xi)=e^{-\tau B_\xi(C,D)}$. 
Here $f=\id:(\p\De, \delta_C)\ra (\p\De, \delta_D)$. 
   Similarly one arrives at  the same conclusion  when the initial segment of
 $m\xi$ lies in $D$ or some $E\not=C, D$.  Since the regular set $B$ has full $\mathcal{H}_{\delta_C}$-measure in 
$\p\De$,  $\id:  (\p\De, \delta_C)\ra (\p\De, \delta_D)$  is a conformal map.

 Now  let $C$, $D$ be  two  arbitrary  chambers.  Then 
$\delta_C$ and  $\delta_D$  are  still conformally equivalent
  since  there exists  a  gallery from $C$ to $D$.  
Since $B_\xi(C_1, C_3)=B_\xi(C_1, C_2)+B_\xi(C_2, C_3)$  holds  for any $\xi\in B$ and 
   any three chambers $C_1$, 
$C_2$, $C_3$,  the formula 
  $\delta_D(\xi,\eta)=e^{-\tau B_\xi(C, D)}\delta_C(\xi,\eta)$  holds for all $\eta$ sufficiently close to $\xi$.

\end{proof}

 If $g:\De\ra \De$ is an isometry, then it is clear that
  the induced map on the boundary $g: (\p\De, \delta_C) \ra  (\p\De, \delta_{g(C)})$ is an isometry for any chamber $C$.
  By Lemma \ref{conbase}, $g: (\p\De, \delta_C) \ra  (\p\De, \delta_{C})$ is   a    conformal map and 
  $\delta_C(g(\xi), g(\eta))=e^{-\tau B_\xi(C,  g^{-1}(C))}\delta_C(\xi,\eta)$ for any fixed regular point $\xi$ and all 
 $\eta$  in a small neighborhood of $\xi$.

Since $\delta_C$ and $\delta_D$  are bi-Lipschitz,  they have the same Hausdorff dimension. Let $H$ be their common
Hausdorff dimension.
Let ${\mathcal{H}}_C$ be the Hausdorff measure  of $\delta_C$.  Lemma \ref{conbase}  implies  that  
${\mathcal{H}}_C$ and ${\mathcal{H}}_D$  are in the same measure class and 
$ {\mathcal{H}}_D(\xi)=e^{-H\tau B_\xi(C, D)}{\mathcal{H}}_C(\xi)$
 for all regular points $\xi$.  If $g:\De\ra \De$ is an isometry, then we have 
 $ g^*{\mathcal{H}}_C(\xi)=e^{-H\tau B_\xi(C,  g^{-1}(C))}{\mathcal{H}}_C(\xi)$.
Hence the measures $\{{\mathcal{H}}_C\}$   are the so-called conformal measures with respect to 
$Isom(\De)$.

\subsubsection{Sullivan's approach } \label{boucon}

Let $\De_1$, $\De_2$ be two Fuchsian buildings, and $G$  a group acting properly and
 cocompactly on both  $\De_1$  and  $\De_2$.   Then there is an equivariant homeomorphism
 $h:  \p\De_1\ra \p\De_2$.     There are two steps in Sullivan's approach.
  The first step is to show  $h$ preserves the combinatorial cross ratio 
almost everywhere (with respect to   certain  measures).
  The second step is to show the implication: 
\lq\lq $h$ preserves the combinatorial cross ratio almost everywhere" $\Rightarrow$ \lq\lq $h$ extends to an isomorphism
 from $\De_1$  to  $\De_2$".

The measures appearing in Step 1 are the so-called conformal measures. 
 In our case they are 
 Hausdorff measures $\{{\mathcal{H}}_C\}$  
 of  the   combinatorial metrics,
  as introduced in   Section \ref{confm}. Sullivan's   ergodic  arguments (see \cite{S2}) 
   show that if the 
  Hausdorff dimensions of the  combinatorial metrics on $\p\De_1$ and $\p\De_2$ are the same   and
 $h$ is nonsingular with respect to the Hausdorff measures, then 
$h$ preserves the combinatorial cross ratio almost everywhere with respect to the Hausdorff measures.
Here \lq\lq $h$ is nonsingular" means a subset $X\subset \p\De_1$  has measure 0 if and only if 
$h(X)$ has measure 0.  Hence to complete Step 1, one has to show \newline
(1) the 
  Hausdorff dimensions of the  combinatorial metrics on $\p\De_1$ and $\p\De_2$ are the same;\newline
(2) $h$ is nonsingular with respect to the Hausdorff measures.  \newline
These two facts have been verified by M. Bourdon  (see \cite{Bo1}, \cite{Bo2}),  and Step 1 is 
 established.
Step  2 follows from  
 Theorem \ref{main}.  

\b{remark}\label{modicm}
{As mentioned at the beginning of Section \ref{crossr}  our definitions of combinatorial cross ratio
 and combinatorial metric are slightly different from  those of M. Bourdon since we only used vertices in the dual graph
 while  M. Bourdon used the whole dual graph. Since the vertices form a net in the dual graph, our combinatorial cross ratio
 differs from M. Bourdon's by an  additive constant that depend only on the building. It follows that the two combinatorial metrics are bi-Lipschitz  equivalent  with the Lipschitz constant depending only on the building, and the corresponding Hausdorff measures are in the same measure class.}  

\end{remark}

\subsection{Triangles and quadrilaterals in the 1-skeleton}\label{triqua1sta}

In this section we state the main results concerning triangles and quadrilaterals contained in the 1-skeleton
  of $\De$.  Their proofs are  somehow  tedious and  are  contained  in Section \ref{triqua}.

Let $\De$ be  a  Fuchsian building.  A subset $T\subset \De$ is called a  \e{triangle} if there are points
 $x,y,z\in \De$ such that $T=xy\cup yz\cup zx$.  In this case, we also use the notation $T=(x,y,z)$.
  The three points $x, y,z$ shall be called the corners of $T$, and $xy$, $yz$, $zx$  called 
  the sides of $T$.
 Similarly, we say  a subset $Q\subset \De$ is a   \e{quadrilateral}
  if there are $x_1, x_2, x_3, x_4\in  \De$ such that $Q=x_1x_2\cup x_2x_3\cup x_3x_4\cup x_4x_1$.  
 We use the notation $Q=(x_1, x_2, x_3, x_4)$,   and call $x_i$ a corner of $Q$  and $x_ix_{i+1}$ a side 
of  $Q$. 
We are mainly interested in triangles and quadrilaterals that are contained in the 1-skeleton of $\De$.

Let $D\subset \De$  be a finite subcomplex that is homeomorphic to a  compact surface  with boundary. 
A vertex $v\in D$ is a  \e{special point} if $v\in \interior(D)$  and  $\length(\Link(D,v))>2\pi$, 
 or $v\in \p D$   and  $ \length(\Link(D,v))>\pi$. 
The following result   implies that each triangle in $\De^{(1)}$  bounds a convex disk.


\b{Prop}\label{tri1}
{Let $\De$  be a Fuchsian building and $T\subset \De^{(1)}$ a triangle that is 
homeomorphic to a   circle.  Then there is  a  finite subcomplex $S(T)$ 
with the following properties:\newline
\e{(1)}  $S(T)$  is  homeomorphic  to  a  closed disk  with boundary $T$; \newline
\e{(2)}   $S(T)$ has no special points.\newline
In particular,  $S(T)$  is convex in $\De$.  }

\end{Prop}

We  list all the triangles in $\De^{(1)}$.
 Recall for integers $m_1, m_2, m_3\ge 2$, $(m_1, m_2, m_3)$ denotes  the triangle (unique up to isometry) 
  in $\H^2$  with angles $\pi/m_1$, $\pi/m_2$  and   $\pi/m_3$.

\b{Prop}\label{tri2}
{Let $\De$ be a Fuchsian building  with chamber $R$,   and $T\subset \De^{(1)}$ a triangle   homeomorphic to
 a circle.   \newline
\e{(1)} If $R$ is not  a triangle, then there is no triangle in $\De^{(1)}$   
homeomorphic to
 a circle;\newline
\e{(2)} If $R$ is a triangle  with all angles $<\pi/2$,  then  $S(T)$ is a chamber;\newline
\e{(3)} If $R=(2, m_1, m_2)$  with $(m_1, m_2)=(6,6), (6,8)$ or   $(8,8)$, then $S(T)$ 
is  isomorphic to one of the  labeled  complexes  in Figure  \ref{268};


\begin{figure}[h]
\centering\epsfig{file=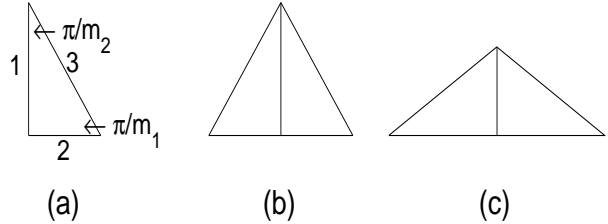,  height=3cm,  width=8cm}\caption{$S(T)$ when $R=(2, m_1, m_2)$  with $(m_1, m_2)=(6,6), (6,8)$ or $(8,8)$}\label{268}
\end{figure}

\e{(4)} If $R=(2, 4,6)$ or $(2,4,8)$, then $S(T)$ 
is  isomorphic to one of the  labeled  complexes  in Figure  \ref{248};

\begin{figure}[h]
\centering\epsfig{file=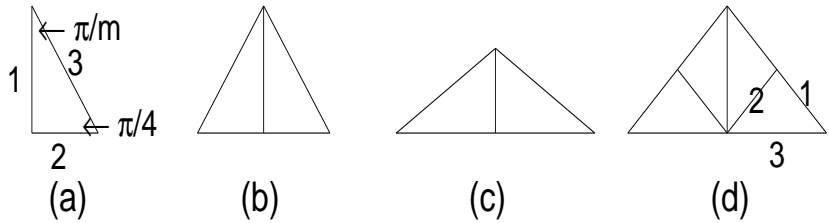,  height=3cm,  width=11cm}\caption{$S(T)$ when $R=(2, 4, m)$  with $m=6$ or $8$}\label{248}
\end{figure}

\e{(5)}  If $R=(2,3,  8)$,  then $S(T)$ 
is  isomorphic to one of the  labeled  complexes  in Figure  \ref{238}.

\begin{figure}[h]
\centering\epsfig{file=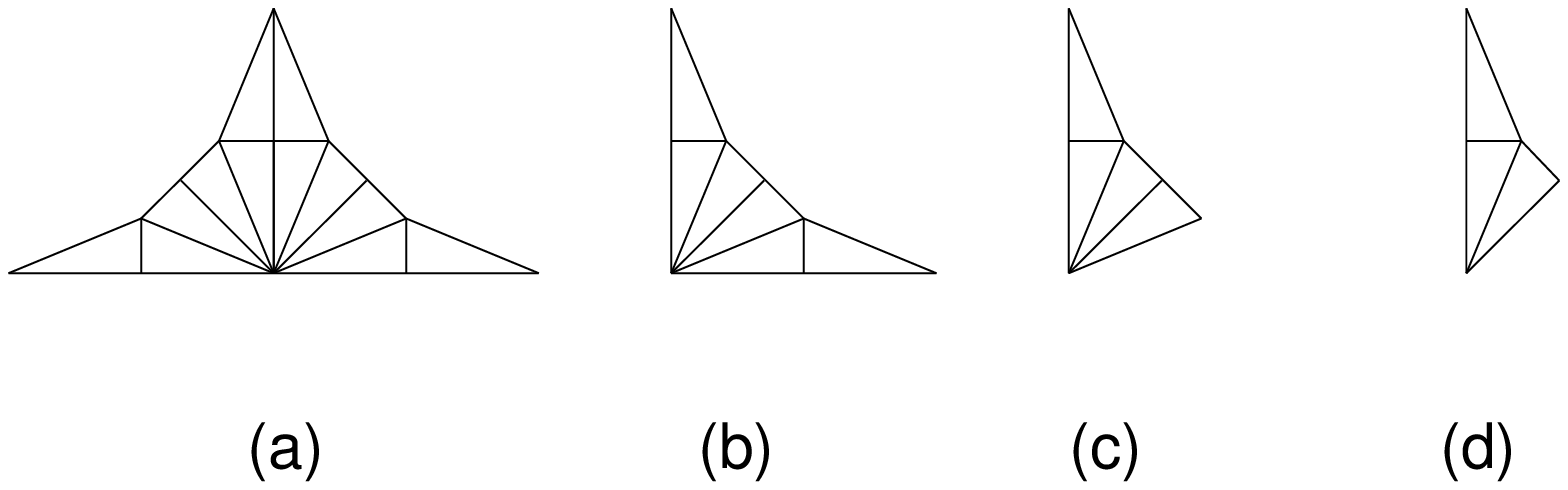,  height=3cm,  width=11cm}\label{238A}
\end{figure}

\begin{figure}[h]
\centering\epsfig{file=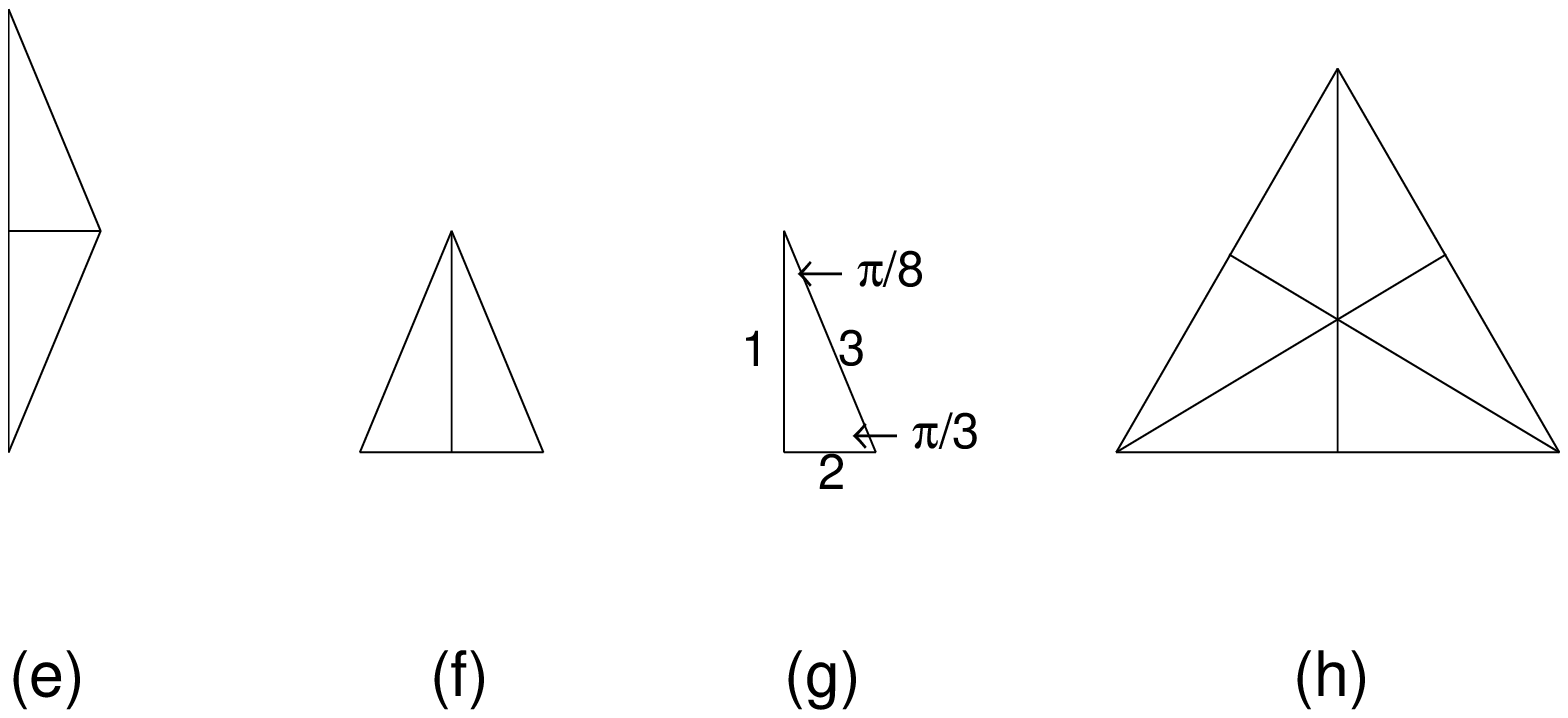,  height=4cm,  width=11cm}\caption{$S(T)$ when $R=(2, 3, 8)$  }\label{238}
\end{figure}

}

\end{Prop}

We let $A_0$  denote  the area of a chamber.
For any quadrilateral $Q=(x_1, x_2, x_3, x_4)$,    let $\Sigma(Q)=\Sigma_{i=1}^4 \angle_{x_i}(x_{i-1}, x_{i+1})$
 be the sum  of  angles at the 4 corners of $Q$.

\b{Prop}\label{quar11}
{Let $\De$ be a Fuchsian building,   and  $Q\subset \De^{(1)}$ a quadrilateral homeomorphic to a  circle. 
   Then there is a finite subcomplex $S(Q)$  
with the following properties:\newline
\e{(1)} $S(Q)$   is  homeomorphic to a closed disk with boundary $Q$;\newline
\e{(2)}   $2\pi\ge  \Sigma(Q)+  n(Q) A_0$,  where  $n(Q)$ is the number of chambers in  $S(Q)$.}

\end{Prop}

Recall the chamber $R$ is  a   compact convex $k$-gon in $\H^2$.  
We let $\alpha_0$ be the smallest angle of $R$.
Note $A_0\ge \pi/2$  if $k\ge 5$.
  There are only   6 right  triangles that can occur as the chambers of  Fuchsian buildings:
  $(2, 8, 8)$,  $(2, 6, 6)$,  $(2, 6, 8)$,  $(2, 4, 6)$,  $(2, 4, 8)$,  $(2, 3, 8)$.
Their areas are respectively:  $\pi/4$,   $\pi/6$, $5\pi/24$, $\pi/12$, $\pi/8$, $\pi/24$.


\b{Cor}\label{shoq}
{Let $\De$ be  a Fuchsian building whose chamber $R$ is  a $k$-gon with $k\ge 5$,
  and  $Q\subset \De^{(1)}$ a quadrilateral.  Then  $Q$ is  not   homeomorphic to a circle.}

\end{Cor}

\b{proof}
Suppose $Q$ is   homeomorphic to a circle.  
 Since $R$ is a $k$-gon with $k\ge 5$ and all its angles are $\le \pi/2$,
 we have $A_0\ge \pi-\alpha_0$.  Proposition \ref{quar11}  implies $2\pi\ge 4\alpha_0+n(Q)(\pi-\alpha_0)$.  
  It is easy to verify that $n(Q)\le 1$ for each of the cases $\alpha_0=\pi/2, \pi/3, \pi/4, \pi/6,\pi/8$, which means  
  $Q$ is the boundary of  a chamber.  This is a contradiction since $Q$ is a quadrilateral and $R$ is a $k$-gon with $k\ge 5$.

\end{proof}

Recall for any vertex $v\in \De$, the link $\Link(\De, v)$ is a generalized $m$-gon for some $m\ge 2$.
  We set $m(v)=m$ if $\Link(\De, v)$ is  a generalized $m$-gon
  and say $v$  is indexed by $m(v)$.
 The vertices in $\Link(\De, v)$ are divided into two types, and two vertices are of the same type if and only
 if the distance between them is an even multiple of $\pi/m$.  Let $vx, vy\subset \De^{(1)}$
  be two (nondegenerate) geodesics starting from $v$. 
  The angle  $\angle_v(x,y)$  is called \e{even} if   it 
is an even multiple of $\pi/m$,  and is  called \e{odd} if  it  is an odd multiple of $\pi/m$.
  Equivalently, $\angle_v(x,y)$  is  even 
 if and only if  the initial directions of $vx$ and $vy$ at $v$  have the same type in 
$\Link(\De, v)$.

\b{Cor}\label{q3in2ev}
{Let   $\De$ be a Fuchsian building  with chamber  $R$,  and 
 $Q=(x, y, z, w)\subset \De^{(1)}$  a quadrilateral   homeomorphic to   a  circle.
  Suppose   $Q$ has even angles at $x$ and $y$,  
   and $xy$ contains at least two vertices in the interior. 
 Then $R=(2,3,8)$.}

\end{Cor}

\b{proof}  
Let $v_1\not= v_2\in \interior(xy)$ be two vertices  with  $v_1\in \interior(xv_2)$. 
We count the chambers in $S(Q)$ that intersect $xy$. There  are at least $m(v_1)$ chambers (in $S(Q)$)
 incident to $v_1$, at least $m(v_2)$ incident to $v_2$; but there might be one chamber incident to both $v_1$ and $v_2$.
 Since $Q$ has an  even angle at $x$, there are at least 2 chambers incident to $x$, but one of them could also be incident to
  $v_1$. Similarly for $y$. Summarize, there are at least $m(v_1)+m(v_2)+1$ chambers in $S(Q)$, that is, 
$n(Q)\ge m(v_1)+m(v_2)+1$.  In particular,  $n(Q)\ge 5$ always holds. 
Also note $\Sigma(Q)\ge \frac{\pi}{m(z)}+\frac{\pi}{m(w)} + \frac{2\pi}{m(x)} + \frac{2\pi}{m(y)}$. 
    In particular, $\Sigma(Q)\ge 3\pi/4$ always holds.

We suppose $R\not=(2,3,8)$  and will derive a contradiction from this.  We consider several cases.

\noindent
Case 1: $k\ge 5$.  In this case  $A_0\ge \pi/2$.  Proposition \ref{quar11}  implies $n(Q)\le 2$, contradicting to the fact 
$n(Q)\ge 5$. 

 \noindent
Case 2: $k=4$.  If $\alpha_0=\pi/8$,
 then $A_0\ge 3\pi/8$. Proposition \ref{quar11} implies
 $n(Q)\le 3$,  contradiction.  If $\alpha_0=\pi/6$, then $A_0\ge \pi/3$ and $\Sigma(Q)\ge \pi$. 
Proposition \ref{quar11}  again  implies
 $n(Q)\le 3$. We similarly obtain contradictions when $\alpha_0=\pi/4$  or  $\alpha_0=\pi/3$.

Below we will apply Proposition \ref{quar11}  without mentioning it.

\noindent
Case 3:   $R$ is a triangle  with no right angle. In this case $m(v_1), m(v_2)\ge 3$ and we have 
 $n(Q)\ge 7$. 
If $\alpha_0=\pi/8$, then $A_0\ge 5\pi/{24}$  and $n(Q)\le 6$,  contradiction. 
   Similarly  one obtains contradictions when $\alpha_0=\pi/6$  or  $\alpha_0=\pi/4$.

\noindent
Case 4:  $R\not=(2,3,8)$ is  a right  triangle. 
We only give the details  for  $R=(2,4, 6)$, the proofs for the other 4 cases are similar.
So assume $R=(2,4,6)$. In this case, $A_0=\pi/12$. 
 Consider any complete geodesic  $c$ containing $xy$.  The vertices on $c$ are alternatingly  indexed by
 $m_1$  and  $m_2$, where $(m_1,m_2)=(2,6), (4,6)$ or $(2,4)$.  Notice none of 
  $x$, $y$ is indexed by $2$, otherwise the quadrilateral is actually a triangle   and one of its 
 sides
 contains   at least 4 edges, contradicting to Proposition \ref{tri2}. 
First  assume  $(m_1,m_2)=(2,4)$. Then $m(x)=m(y)=4$  and   
    $n(Q)\le 8$. There is at least  one vertex $v\in \interior(xy)$ indexed by $4$. 
   Also recall the angles of $Q$ at $x$ and $y$ are even. It follows that there   are  at least 4 chambers in
 $S(Q)$ incident to $v$ and at least two chambers incident to each of $x$, $y$,
  for   a  total of 8. Hence $S(Q)$ is exactly the union of these 8 chambers.  But   this  union 
    is  a hexagon, not a quadrilateral.

Next assume $(m_1,m_2)=(4,6)$. 
Suppose $\{m(x),m(y)\}=\{4,6\}$. Then $n(Q)\le 10$. There are  two  vertices
  $v_1, v_2\in \interior(xy)$ with $m(v_1)=6$, $m(v_2)=4$. It  follows that 
$n(Q)\ge 6+4+1=11$, contradiction.
One similarly obtains contradictions  when   $m(x)=m(y)=4$  or   $6$.

Finally we assume $(m_1,m_2)=(2,6)$.
 Then $m(x)=m(y)=6$  and $n(Q)\le 12$. 
There is at least one    vertex $v\in \interior(xy)$  with $m(v)=6$. We count chambers in $S(Q)$ incident to vertices
 indexed by $6$: 
 (at least)  6 chambers incident to $v$,   2  incident to each of $x$, $y$.
 So  $n(Q)\ge 10$.   Proposition \ref{quar11}  implies  the angles at $x$ and $y$ are
  $\pi/3$ and 
none of the other two angles of $Q$ is $\pi/2$. It follows that the vertices on $xw$ are alternatingly indexed by
 $2$ and $6$, and $m(w)\not=2$.  Hence $m(w)=6$. Similarly $m(z)=6$.  Notice $\interior(xw)$ contains 
 exactly one vertex  and so does $\interior(yz)$.  
There is (at least) 
   one chamber in $S(Q)$  incident to each of $z$, $w$.  We have exhibited 12 chambers in $S(Q)$. It follows that $S(Q)$ 
is the union of these 12 chambers. However, one can check that this  union  is not a quadrilateral.

\end{proof}

\b{Cor}\label{threeenenq}
{Let   $\De$ be a Fuchsian building  with chamber  $R$,  and 
 $Q=(x, y, z, w)\subset \De^{(1)}$  a quadrilateral   homeomorphic to   a  circle.
If  $Q$ has three even angles, then $R$ is a right  triangle.}

\end{Cor}

\b{proof} Suppose $R$ is not  a right  triangle.  We may assume the angles at $x$, $y$ and $z$ are even.
We count the chambers in $S(Q)$  that intersect $xy$:  at least 2 chambers incident to $x$ and $y$ each, but there might be 
a chamber incident to both $x$ and $y$.  Hence $n(Q)\ge 3$. On the other hand, 
  the assumption implies $\Sigma(Q)\ge 7\alpha_0$, where $\alpha_0$ is the smallest angle of $R$. 
Recall $R$ is  a  $k$-gon. 

\noindent 
Case 1: $k\ge 5$.  In this case, $A_0\ge \pi/2$ and $\Sigma(Q)\ge 7\pi/8$.  
 Proposition \ref{quar11}  (2)  implies $n(Q)\le 2$, contradicting to the fact 
$n(Q)\ge 3$.

\noindent 
Case 2: $k=4$.  As above, Proposition \ref{quar11}  (2)  implies $n(Q)\le 2$ if $\alpha_0=\pi/6$,  $\pi/4$  or  $\pi/3$.
 Assume $\alpha_0=\pi/8$.  Then $n(Q)\le 3$. It follows that $S(Q)$ is the union of the three chambers indicated in the first 
 paragraph.  However, this union is not a quadrilateral.

\noindent 
Case 3: $R$ is  a triangle with no right angle.   
Then Proposition \ref{quar11}  (2)  implies $n(Q)\le 5$ in each of the cases: 
 $\alpha_0=\pi/8,  \pi/6$ or $\pi/4$. Let $\alpha_0\le \beta_0\le \gamma_0\le \pi/3$ be the three angles of $R$. 
   We first assume  one of the  segments   $xy$,  $yz$   contains a vertex  $v$ in the interior, say $v\in \interior(xy)$.  
     Then $m(v)\ge 3$  and 
  there are at least 3 chambers in $S(Q)$ incident to $v$. On the other hand, 
    there are at least two chambers in $S(Q)$ incident to $x$, and at  least one of them 
is distinct from the three chambers incident to $v$. Similarly 
there is at least one chamber in 
$S(Q)$ incident to $y$ and  different  from the three chambers incident to $v$.
  Combining  with the observation  $n(Q)\le 5$, we see  $n(Q)=5$.    
   However, the union of these 5 chambers is either a pentagon or a hexagon. 
 The contradiction means $xy$ and $yz$ are edges in $\De$.

  Since the angle at $y$ is even, 
 $m(x)=m(z)$  holds.  It implies   that  the sum of angles at $x$, $y$ and $z$ is $\ge 2\beta_0+4\alpha_0$.
  Since the fourth angle of $Q$ is $\ge \alpha_0$ and
   $A_0= \pi-(\alpha_0+\beta_0+\gamma_0)\ge  2\pi/3-\alpha_0-\beta_0$,  Proposition \ref{quar11}
 (2)  implies
  $2\pi\ge 5\alpha_0+2 \beta_0+m(2\pi/3-\alpha_0-\beta_0)$.  
  By using $\alpha_0\le \pi/4$ and analyzing all  the   cases, one concludes that $n(Q)\le 3$. 
Here we have a contradiction, as one can check that the union of three chambers 
   can nerve be a quadrilateral
 with three even angles.

\end{proof}






\subsection{Geodesics   at different sides} \label{gedads}

In this section we discuss a condition (\lq\lq  being at different sides") on two geodesics 
 that  in most cases  guarantees
  the two geodesics intersect.   This notion  is defined in Definition \ref{atdsftg}.  
 In Section \ref{sufcond} we   give sufficient conditions for two 
 geodesics to be at different sides.  In Section \ref{inters} we show that in most cases
two geodesics at different sides 
 must intersect.

\subsubsection{Criterion for geodesics to be at different sides} \label{sufcond}


Recall for each vertex $v$ of $\De$, the link $\Link(\De,v)$ is a   rank two  spherical building. 

\b{Def}\label{localatsd}
{Let $\xi_1, \xi_2, \eta_1, \eta_2\in \p\De$  be  such that the geodesics 
 $\xi_1\xi_2$, $\eta_1\eta_2$  are contained in the 1-skeleton  and 
intersect at a single point $v$.  We say  \e{$\xi_1\xi_2$  and  $\eta_1\eta_2$
 are locally contained in an apartment}  if there is an apartment  in $\Link(\De, v)$
  that contains the four directions at $v$ induced by $\xi_1\xi_2$  and  $\eta_1\eta_2$.}

\end{Def}

\b{Prop}\label{ssufdas}
{Suppose  two  geodesics $\xi_1\xi_2$  and  $\eta_1\eta_2$
 are locally contained in an apartment and meet at a vertex $v$.  Then 
$\xi_1\xi_2$, $\eta_1\eta_2$ are  at different sides  if  any one of the following 
  conditions holds:\newline
\e{(1)} $R$ is not a right    triangle;\newline
\e{(2)} $R\not=(2,3,8)$  and $\xi_1\xi_2$  and  $\eta_1\eta_2$
 make a right angle at $v$.}

\end{Prop}

The proof of Proposition \ref{ssufdas}  is contained in the following two lemmas.

  Notice that if $p_1$, $p_2$ are two vertices in a generalized polygon opposite to the 
 same vertex $p$, then $p_1$ and $p_2$ are of the same type. It follows that 
if two geodesics $y_1y_2,y_1y_3\subset \De^{(1)}$  intersect in a nontrivial segment $y_1x$ ($x$ lies in the interior of 
 $y_1y_2$, $y_1y_3$), then  $xy_2$ and $xy_3$   make  a  nonzero  even angle at $x$.

\b{Le}\label{noritri}
{Suppose  two  geodesics $\xi_1\xi_2$  and  $\eta_1\eta_2$
 are locally contained in an apartment.  Then 
$\xi_1\xi_2$, $\eta_1\eta_2$ are  at different sides  if
  $R$ is not a right      triangle.}

\end{Le}

\b{proof} Suppose   $R$ is not a right      triangle.
 Let $v=\xi_1\xi_2\cap \eta_1\eta_2$.  
 Assume $\xi_i\eta_2\cap \eta_1\eta_2$ is a ray $x\eta_2$ ($x\in \eta_1\eta_2$) for some $i=1,2$. 
Then $\xi_i\eta_2\subset \De^{(1)}$.  Now both $\xi_j\xi_i$ ($j\not=i$) and $\eta_2\xi_i$  are contained 
 in $\De^{(1)}$,   and 
the  distance between  them  is asymptotically 0. It follows that $\xi_j\xi_i\cap \eta_2\xi_i $
  is  a  ray $y\xi_i$ ($y\in \xi_1\xi_2$),  and $(v, x,y)$ is a triangle with an  even angle  at 
 $x$,  contradicting to Proposition \ref{tri2}. Therefore
 $\xi_i\eta_2\cap \eta_1\eta_2=\phi$   for $i=1,2$.
  Similarly   for each $y\in  \xi_1\xi_2$, $y\not=v$,
  we have $y\eta_2\cap \eta_1\eta_2=\phi$.

  Continuity shows that  each $y\in  [\xi_1, \xi_2]:=\xi_1\xi_2\cup \{\xi_1, \xi_2\}$, $y\not=v$  has a neighborhood  $U$ in 
 $[\xi_1, \xi_2]-\{v\}$ such that  all the points in $U$ lie at the same side of $\eta_2$.  
 The connectivity of $(v,\xi_i]:=(v\xi_i-\{v\})\cup \{\xi_i\}$ implies that  all   the  points in 
$(v,\xi_i]$ lie at the same side of $\eta_2$.  
 On the other hand,  $\xi_1\xi_2$  and  $\eta_1\eta_2$
 are locally contained in an apartment.  
  It is not hard to see that,  if $y,z$    lie in different components of $\xi_1\xi_2-\{v\}$  and 
are sufficiently close to $v$, then $y$ and $z$ lie at different sides of $\eta_2$.  It follows that
  $\xi_1$ and $\xi_2$ lie at different sides of $\eta_2$.   Similarly 
 $\xi_1$ and $\xi_2$ lie at different sides of $\eta_1$,  and 
$\eta_1$ and $\eta_2$ lie at different sides of $\xi_i$ ($i=1,2$).

\end{proof}

\b{Le}\label{ritrdirigh}
{Suppose  two  geodesics $\xi_1\xi_2$  and  $\eta_1\eta_2$
 are locally contained in an apartment and meet at a vertex $v$.  Then 
$\xi_1\xi_2$, $\eta_1\eta_2$ are  at different sides  if  
$R\not=(2,3,8)$  and $\xi_1\xi_2$  and  $\eta_1\eta_2$
 make a right angle at $v$.}

\end{Le}

\b{proof} Suppose $R\not=(2,3,8)$  and   
$\xi_1\xi_2$  and  $\eta_1\eta_2$
 make a right angle at $v$.
  Assume   there is some $q\in \xi_1\xi_2$, $q\not=v$ such that $q\eta_2\cap \eta_1\eta_2$  is a ray
 $x\eta_2$.  Then $(v,  q,  x)\subset \De^{(1)}$ is a triangle      
 with a right angle at $v$ and an even angle at $x$, 
contradicting to  Proposition \ref{tri2}.  One continues to argue as in the proof of Lemma \ref{noritri}
 that $\xi_1\xi_2$, $\eta_1\eta_2$ are  at different sides.

\end{proof}

\b{remark}\label{insski}
{In the case when both $R_1$ and $R_2$ have 
 at least 5 edges it 
 is possible  to  understand the proof of Theorem \ref{main} without having to go through 
too much details, if one is willing to assume Propositions  \ref{tri1}, \ref{tri2} and \ref{quar11}.  
   Here is how to do so. 
Read Section \ref{pre}, Section \ref{crossr}, Section \ref{triqua1sta} (skip Corollaries  \ref{q3in2ev},
\ref{threeenenq}),  Section \ref{sufcond} (skip Lemma \ref{ritrdirigh}), Lemma \ref{onedif} of Section \ref{inters}, 
Section \ref{nonrightan} (skip Lemmas \ref{commint}, \ref{commintmis3}  and Proposition \ref{onetoneq})
    and  the following: 

Let $\De$ be a Fuchsian building whose   chamber $R$ is a $k$-gon with $k\ge 5$. 
Then Proposition \ref{tri2} and  Corollary \ref{shoq} imply that no triangle or quadrilateral 
 in $\De^{(1)}$ can be homeomorphic to a circle. 

We  claim  that  if  two geodesics   are  at different sides  then they  intersect.
Suppose 
$\xi_1\xi_2, \eta_1\eta_2\subset \De^{(1)}$  ($\xi_1, \xi_2, \eta_1, \eta_2\in \p\De$)
are disjoint  and  at different sides.    Then Lemma \ref{onedif}  implies that 
for each $i=1,2$ 
 there are  vertices $x_i\in \xi_1\xi_2$, $y_i\in \eta_1\eta_2$ such that
  $x_i\eta_i\cap \eta_1\eta_2=x'_i\eta_i$ and $y_i\xi_i\cap \xi_1\xi_2=y'_i\xi_i$ are rays.
 Since  no triangle or quadrilateral 
 in $\De^{(1)}$ can be homeomorphic to a circle, we have $x_1=x_2=y'_1=y'_2$ and $x'_1=x'_2=y_1=y_2$.
  It follows that $x_1x'_1$ makes an angle $\pi$ with the two rays $\xi_1\xi_2, \eta_1\eta_2$ 
  and      hence  $\xi_1\eta_1=x_1\xi_1\cup x_1x'_1\cup x'_1\eta_1\subset \De^{(1)}$,  
          contradicting to the definition of 
 being at different sides.

Proof of condition (1) in Lemma \ref{redaprtowh} when $R_1$ and $R_2$ have 
 at least 5 edges:  Let $c_1, c_2\subset A^{(1)}$ be two geodesics through $v$.  $c_1$, $c_2$ are clearly 
at different sides.  Lemma \ref{redaprtowh0.1}  implies $c'_1$, $c'_2$ are at different sides. The preceding paragraph shows 
 $c'_1$ and $c'_2$ intersect at a vertex $w\in \De_2$.  Now let $c_3\subset A^{(1)}$ be any other geodesic through $v$.
 Then $w_1:=c'_1\cap c'_3$ and $w_2:=c'_2\cap c'_3$ are vertices.  If $w\notin c'_3$, 
  then $(w, w_1, w_2)\subset \De_2^{(1)}$ is homeomorphic to a circle, contradicting to the above observation.

}
\end{remark}

\subsubsection{Intersection of geodesics at different sides}  \label{inters}

In this section we show, in most cases, that  geodesics at different sides must intersect with each other.

The proof of Lemma \ref{noritri}  also shows the following: 

\b{Le}\label{onedif}
{If $\xi_1\xi_2, \eta_1\eta_2\subset \De^{(1)}$  \e{($\xi_1, \xi_2, \eta_1, \eta_2\in \p\De$)}
are disjoint  and  at different sides,  then for each $i=1,2$ 
 there are  vertices $x_i\in \xi_1\xi_2$, $y_i\in \eta_1\eta_2$ such that
  $x_i\eta_i\cap \eta_1\eta_2$ and $y_i\xi_i\cap \xi_1\xi_2$ are rays.}

\end{Le}

\qed

\b{Prop}\label{samegeo}
{Suppose $\xi_1\xi_2, \eta_1\eta_2\subset \De^{(1)}$ \e{($\xi_1, \xi_2, \eta_1, \eta_2\in \p\De$)}
   are disjoint  and  at different sides. Let 
$y_i\in \eta_1\eta_2$  \e{($i=1,2$)}  be a vertex such that $y_i\xi_i\cap \xi_1\xi_2=y_i'\xi_i$
  is   a  ray   and  $y_i\xi_i\cap \eta_1\eta_2=\{y_i\}$.
  If $R\not=(2,3,8)$,  
  then one of the following  holds:\newline
\e{(1)}  $y_1'\in \interior(y_2'\xi_1)$ and $y_1\not=y_2$; \newline
\e{(2)} $y_1'=y_2'$ and  $y_1\not=y_2$. }

\end{Prop}

\b{proof}
 Since  $\De$ is a $\CAT(-1)$ space,  there is a unique geodesic segment between two points. 
We first notice $y_1'\in \interior(y_2'\xi_2)$ cannot happen, since otherwise
 $y_1y_1'\cup y_1'y_2'\cup y_2'y_2$ would be a  geodesic segment  connecting $y_1\in \eta_1\eta_2$
  and $y_2\in \eta_1\eta_2$,  but 
different from $y_1y_2\subset \eta_1\eta_2$.  Therefore $y_1'\in y_2'\xi_1$.

Suppose $y_1'\not=y_2'$, that is,  $y_1'\in \interior(y_2'\xi_1)$.  Notice  that  the two angles 
 $\angle_{y_1'}(y_1, \xi_2)$ and $\angle_{y_2'}(y_2, \xi_1)$  are even.  
  If $y_1=y_2$, then $(y_1, y_1', y_2')$ is a triangle with two even angles, contradicting
 to Proposition \ref{tri2}
  and the assumption that $R\not=(2,3,8)$.  Therefore  $y_1\not=y_2$
  and  (1)  holds.

Now suppose $y_1'=y_2'$.    We need to show  $y_1\not=y_2$.
 Suppose  $y_1=y_2$. Let $x_i\in \xi_1\xi_2$ ($i=1,2$)   be  a vertex such that
 $x_i\eta_i\cap \eta_1\eta_2=x'_i\eta_i$ is  a  geodesic ray  and  $x_i\eta_i\cap \xi_1\xi_2=\{x_i\}$. 
  Assume    $x'_i=y_1$ for some $i=1,2$. Then the uniqueness of geodesic implies
 $x_i=y'_1$. It follows that $\eta_i\xi_1=\eta_iy_1\cup y_1y'_1\cup y'_1\xi_1\subset \De^{(1)}$
  and $\eta_i\xi_1\cap \xi_2\xi_1$ is a ray, 
    contradicting to the assumption that 
$\xi_1\xi_2$ and $\eta_1\eta_2$ are at different sides.   Therefore
   $x'_i\not=y_1$ for  $i=1,2$. Notice $x'_i\in \interior(y_1\eta_i)$,
  otherwise $x_ix'_i\cup x'_iy_1$ and $x_iy'_1\cup y'_1y_1$ would be two distinct geodesic segments from $x_i$ to $y_1$.  At least one of the two angles 
$\angle_{y_1}(x'_1,y'_1)$, $\angle_{y_1}(x'_2,y'_1)$  is $\ge \pi/2$. We may assume
$\angle_{y_1}(x'_1,y'_1)\ge \pi/2$.  Then $x_1y_1=x_1y'_1\cup y'_1y_1$ and 
$(x_1, y_1, x'_1)$  is  a  triangle with an even angle  at $x'_1$  and   
   another  angle  $\angle_{y_1}(x'_1,x_1)=\angle_{y_1}(x'_1,y'_1)\ge \pi/2$.
Proposition \ref{tri2}
implies  $R=(2,3,8)$, contradicting to our assumption. 

\end{proof}

\b{Prop}\label{tgeoddif}
{Suppose $\xi_1\xi_2, \eta_1\eta_2\subset \De^{(1)}$  \e{($\xi_1, \xi_2, \eta_1, \eta_2\in \p\De$)}
  are disjoint  and at different sides. 
Let $x_1\in \xi_1\xi_2$, $y_1\in \eta_1\eta_2$ be vertices 
 such that $x_1\eta_1\cap \eta_1\eta_2=x_1'\eta_1$ and $y_1\xi_1\cap \xi_1\xi_2=y_1'\xi_1$ are rays
 and  $x_1\eta_1\cap \xi_1\xi_2=\{x_1\}$,  $y_1\xi_1\cap  \eta_1\eta_2=\{y_1\}$.
  Then after possibly switching  $\xi$ with $\eta$ and $x$ with $y$,
    one of the following   holds:\newline
\e{(1)}  $y_1'\in\interior(x_1\xi_1)$  and $y_1\in \interior(x_1'\eta_1)$;\newline
\e{(2)} $y_1'\in\interior(x_1\xi_1)$  and  $y_1=x_1'$;\newline
\e{(3)}  $y_1'\in\interior(x_1\xi_1)$  and  $x_1'\in \interior(y_1\eta_1)$.}
\end{Prop}

\b{proof}
  Assume     $x_1=y_1'$ and $x_1'=y_1$.  
 Then   $\xi_1\eta_1=\xi_1 y_1'\cup y_1'y_1\cup y_1\eta_1$, which 
  implies    $\xi_1\eta_1\cap  \xi_1\xi_2=y_1'\xi_1$ is a ray, contradicting to the assumption that 
$\eta_1$ and $\eta_2$  lie  at different sides of $\xi_1$. Hence  either  $x_1\not=y_1'$ or  $x_1'\not=y_1$. 
We assume $x_1\not=y_1'$. The case $x_1'\not=y_1$ can be handled similarly.

Suppose
 $x_1\not=y_1'$ and $x_1'=y_1$.  Then either 
$y_1'\in\interior(x_1\xi_1)$  or $x_1\in\interior(y'_1\xi_1)$. 
 Since $y_1\xi_1\cap \xi_1\xi_2=y_1'\xi_1$, the uniqueness of geodesic implies that 
$x_1\in\interior(y'_1\xi_1)$ cannot happen and we have (2). 

Suppose  $x_1\not=y_1'$ and  $x_1'\not=y_1$.
There are 4 possibilities:\newline
 (a) $y_1'\in\interior(x_1\xi_1)$ and $y_1\in \interior(x_1'\eta_1)$;\newline
(b) $y_1'\in\interior(x_1\xi_1)$ and $x'_1\in \interior(y_1\eta_1)$;\newline
(c) $x_1\in\interior(y'_1\xi_1)$ and $x'_1\in \interior(y_1\eta_1)$;\newline
(d) $x_1\in\interior(y'_1\xi_1)$ and  $y_1\in \interior(x_1'\eta_1)$. \newline
Case (d) cannot happen since otherwise there are two different geodesic segments 
$y_1y_1'\cup y_1'x_1$ and $y_1x_1'\cup x_1'x_1$ from $y_1$ to $x_1$, a contradiction. 
(a) corresponds to (1) and (b) corresponds to (3).  (c) corresponds to (1) after switching  
 $x$ with $y$ and  $\xi$ with $\eta$.

\end{proof}

\b{Prop}\label{klarg4}
{Let $\xi_1\xi_2, \eta_1\eta_2\subset \De^{(1)}$  \e{($\xi_1, \xi_2, \eta_1, \eta_2\in \p\De$)}
   be   two   disjoint geodesics.  If  
$\xi_1\xi_2, \eta_1\eta_2$   are at different sides,  then $R$ 
 is  a right   triangle.}

\end{Prop}

We need some lemmas. 

\b{Le}\label{qandtrinori}
{Suppose $\xi_1\xi_2, \eta_1\eta_2$  are disjoint and  at different sides. 
Let $x_1\in \xi_1\xi_2$, $y_1\in \eta_1\eta_2$ be vertices 
 such that $x_1\eta_1\cap \eta_1\eta_2=x_1'\eta_1$ and $y_1\xi_1\cap \xi_1\xi_2=y_1'\xi_1$ are rays
  and  $x_1\eta_1\cap \xi_1\xi_2=\{x_1\}$,  $y_1\xi_1\cap  \eta_1\eta_2=\{y_1\}$.
  If $R$ is not a right  triangle,  
   then   Proposition \ref{tgeoddif} \e{(3)}  holds.}  

\end{Le}

\b{proof}
We prove by contradiction that Proposition \ref{tgeoddif} (1)
 and (2) cannot happen.  Since $R$ is not a right  triangle,  Proposition \ref{tri2}
   implies there is no triangle contained in the 1-skeleton 
  that has a  nonzero  even angle.
If Proposition \ref{tgeoddif} (2)
holds, then $(x_1, y_1', y_1)$ is a triangle with a  nonzero 
    even angle $\angle_{y_1'}(x_1,y_1)$, a contradiction.
 Suppose  Proposition \ref{tgeoddif} (1)
holds.   Then  $x_1y_1=x_1x_1'\cup x_1'y_1$ and  again 
$(x_1, y_1', y_1)$ is a triangle with a  nonzero   even angle $\angle_{y_1'}(x_1,y_1)$, a contradiction.
Finally we notice that Proposition \ref{tgeoddif} {(3)}  remains the same after we switch 
$x$ with  $y$ and $\xi$ with $\eta$.

\end{proof}

\b{Le}\label{hau}
{Suppose $\xi_1\xi_2, \eta_1\eta_2$  are  disjoint and at different sides. Let 
$y_i\in \eta_1\eta_2$  \e{($i=1,2$)} be a vertex such that $y_i\xi_i\cap \xi_1\xi_2=y_i'\xi_i$
  is   a  ray and  $y_i\xi_i\cap \eta_1\eta_2=\{y_i\}$.  If $R$ is not a right  triangle, then 
Proposition \ref{samegeo}  \e{(1)}   holds.}

\end{Le}

\b{proof}
  Assume  Proposition \ref{samegeo}  (2) holds, that is, $y_1'=y_2'$ and  $y_1\not=y_2$.
By Lemma \ref{qandtrinori},  $y_1'\in\interior(x_1\xi_1)$  and  $x_1'\in \interior(y_1\eta_1)$.
 If $y_2\in \interior(x_1'\eta_2)$, then $x_1y_2=x_1y_1'\cup y_1'y_2$ and 
$(x_1, y_2, x_1')$ is a triangle with a  nonzero  even angle $\angle_{x_1'}(x_1, y_2)$,
  contradicting to Proposition \ref{tri2}.
 If $y_2\in x_1'\eta_1$,  then $x_1y_1'\cup y_1'y_2$  and $x_1x_1'\cup x_1'y_2$
  are two distinct   geodesic segments from $x_1$ to $y_2$,   a contradiction.


\end{proof}

\noindent
{\bf{Proof of Proposition \ref{klarg4}.}}
Suppose $R$ is not a right  triangle.  By  Lemma \ref{onedif}   there are vertices 
 $y_i\in \eta_1\eta_2$ ($i=1,2$) such that $y_1\xi_1\cap \xi_1\xi_2=y'_1\xi_1$ and
$y_2\xi_2\cap \xi_1\xi_2=y'_2\xi_2$ are rays.  We may assume 
   $y_1\xi_1\cap \eta_1\eta_2=\{y_1\}$,  
  $y_2\xi_2\cap \eta_1\eta_2=\{y_2\}$.  Lemma \ref{hau} implies 
$y'_1\in \interior(y'_2\xi_1)$  and $y_1\not=y_2$.   Let $\bar{\eta}\not=\bar{\bar{\eta}}\in \{\eta_1, \eta_2\}$
   such that  $y_1\in interior(y_2\bar{\eta})$.   Let $\bar{x}\in \xi_1\xi_2$ be  a vertex such that 
 $\bar{x}\bar{\eta}\cap \bar{\eta}\bar{\bar{\eta}}=\bar{x}'\bar{\eta}$ is  a ray
  and  $\bar{x}\bar{\eta}\cap \xi_1\xi_2=\{\bar{x}\}$. 
   Lemma \ref{qandtrinori}
  implies  the six points $\{y_1, y_1', y_2, y_2', \bar{x}, \bar{x}'\}$ are pairwise distinct. 
 Applying Lemma \ref{qandtrinori} to the rays $y_1\xi_1$ and $\bar{x}\bar{\eta}$ we see 
   $\bar{x}\in \interior(y'_1\xi_2)$. 
  Applying Lemma \ref{qandtrinori}  again  to the rays $y_2\xi_2$ and $\bar{x}\bar{\eta}$ we see 
$\bar{x}\in \interior(y'_2\xi_1)$.  It follows that $\bar{x}\in \interior(y'_1y'_2)$.
 Similarly,  if  $\bar{\bar{x}}\in \xi_1\xi_2$ is    a vertex such that 
 $\bar{\bar{x}}\bar{\bar{\eta}}\cap \bar{\eta}\bar{\bar{\eta}}=\bar{\bar{x}}'\bar{\bar{\eta}}$ is  a ray,
  then $\bar{\bar{x}}\in \interior(y'_1y'_2)$.  Lemma \ref{hau} implies 
$\bar{\bar{x}}\not=\bar{x}$.  Therefore $y'_1y'_2$ contains   at least two vertices 
in its interior.   Notice  $(y'_1, y_1,  y_2,  y'_2)\subset \De^{(1)}$ 
  is  a  quadrilateral with  even angles 
  at two  adjacent corners  $y'_1$  and   $y'_2$.
  Now   Corollary   \ref{q3in2ev}
  implies that $R=(2,3,8)$, contradicting to our assumption.

\qed

 The following  lemma follows easily from the definition.

\b{Le}\label{diinet1}
{Let $\xi_1, \xi_2, \eta_1, \eta_2\in \p\De$ with  $\xi_1\xi_2, \eta_1\eta_2\subset \De^{(1)}$.
Suppose $\xi_1\xi_2, \eta_1\eta_2$ are at different sides and $\xi_1\xi_2\cap \eta_1\eta_2\not=\phi$,
  then  $\xi_1\xi_2\cap \eta_1\eta_2$ is  a vertex.}

\end{Le}

\subsection{When $R_1$, $R_2$ are not right  triangles}\label{nonrightan}

In this section we prove   Theorem \ref{main}
in the case when  the chambers  are not right  triangles.

\subsubsection{Criterion for extension}\label{criexten}

In this section we find  conditions that  ensure  a homeomorphism 
$h:\p\De_1\ra \p\De_2$  extends  to  an isomorphism $\De_1\ra \De_2$.  

Let $\De_1$  and   $\De_2$  be Fuchsian buildings with chambers $R_1$ and $R_2$ respectively,  
and $h:\p\De_1\ra \p\De_2$  a homeomorphism that preserves the combinatorial cross ratio almost everywhere.
For any $\xi\in \p\De_1$,  let $\xi'=h(\xi)\in \p\De_2$. For $\xi, \eta\in \p\De_1$,
  we call $\xi'\eta'$ the   \e{image}   of $\xi\eta$.  


Since nonempty open subsets in $\p\De_1$  and   $\p\De_2$  have positive measures  and 
 $h$ preserves the combinatorial cross ratio a.e.,  
    Lemma \ref{geo1ske}
  implies    the following  result.

\b{Le}\label{redaprtowh0}
{For any  $\xi, \eta\in \p \De_1$,   
$\xi\eta$ lies in the 1-skeleton  of $\De_1$ 
  if and only if  its image 
 lies in the 1-skeleton of $\De_2$.
  Furthermore,   any  edge in $\xi\eta$ is  contained in exactly $q+1$ chambers
 if and only   if  any edge in $\xi'\eta'$ is contained in exactly $q+1$ chambers.}

\end{Le}

Lemma \ref{redaprtowh0}  generalizes Lemma 2.4.8 of \cite{Bo1}. It implies that $h(B_1)=B_2$,
 where $B_i$ ($i=1,2$) is the set of regular points in $\p\De_i$.   
The following crucial lemma follows  easily from Lemma \ref{redaprtowh0}
  and Lemma \ref{difes}.

\b{Le}\label{redaprtowh0.1}
{Let $\xi_1, \xi_2, \eta_1, \eta_2\in \p\De_1$ with  $\xi_1\xi_2, \eta_1\eta_2\subset \De_1^{(1)}$.
 Then $\xi_1\xi_2,   \eta_1\eta_2\subset \De_1$ are at different sides if and only if 
$\xi'_1\xi'_2,   \eta'_1\eta'_2\subset \De_2$ are at different sides.}

\end{Le}

For any vertex $v\in \De_1$, let ${\mathcal{D}}_{v}$  be the family of geodesics passing through $v$ and contained in
   $\De_1^{(1)}$,  and ${\mathcal{D}}'_{v}$  the family of images of the  geodesics in ${\mathcal{D}}_{v}$.  
  For each apartment $A$ containing $v$, let ${\mathcal{D}}_{A,v}$  be the set of geodesics passing through $v$ 
  and contained 
 in $A^{(1)}$,  and ${\mathcal{D}}'_{A, v}$  the set  of images of the geodesics in ${\mathcal{D}}_{A, v}$.
 We shall prove that for any vertex $v\in \De_1$, the  geodesics in ${\mathcal{D}}'_{v}$
  intersect in a unique vertex of  $\De_2$. 

\b{Le}\label{redaprtowh}
{Let $v\in \De_1$  be  a vertex. 
Suppose the following conditions hold:\newline
\e{(1)} for  any apartment $A$ containing $v$,  the geodesics in 
 ${\mathcal{D}}'_{A,v}$   intersect in a unique vertex  $v_A\in \De_2$;\newline
\e{(2)} $v_{A_1}=v_{A_2}$  holds for any two apartments 
 $A_1, A_2\subset \De_1$    containing   $v$   with  
 $\Link(A_1, v)\cap \Link(A_{2}, v)$     a  half apartment in $\Link(\De_1, v)$.\newline
Then   the  geodesics in ${\mathcal{D}}'_{v}$
  intersect in a unique vertex.}

\end{Le}

\b{proof}
Let $c, c'\subset \De_1^{(1)}$ be two geodesics  that pass through $v$.  
 Let $A$  and  $A'$ be two apartments that contain $c$ and $c'$ respectively.  
Lemma \ref{halfa}  implies there is  a sequence of apartments  containing $v$:  $A_0=A,  \cdots, A_n=A'$
  such that  $\Link(A_i, v)\cap \Link(A_{i+1}, v)$ ($0\le i\le n$) is   a  half apartment in $\Link(\De_1, v)$.
Now the lemma follows from the assumptions.

\end{proof}

\b{Le}\label{allha}
{Suppose for any vertex $v\in \De_1$,  the  geodesics in ${\mathcal{D}}'_{v}$
  intersect in a unique vertex of  $\De_2$  denoted by   $f(v)$,
  and  for any vertex $w\in \De_2$,  the  family  of  geodesics 
   $\{h^{-1}(\xi)h^{-1}(\eta):  \xi, \eta\in \p\De_2,\;  w\in \xi\eta\subset \De_2^{(1)}\}$  
  also intersect in a unique vertex of $\De_1$.  
     Then the map $f: \De_1^{(0)}\ra   \De_2^{(0)}$
  is a bijection that extends to an isomorphism from $\De_1$ to    $\De_2$.}

\end{Le}

\b{proof} The map $f$ is clearly bijective. 
We   claim  for any  two vertices $v_1, v_2\in \De_1$,  
  $v_1$ and $v_2$ are adjacent if and only if $f(v_1)$ and $f(v_2)$ are adjacent. 
 The claim implies $f$ extends to an isomorphism between the 1-skeletons of $\De_1$ and 
 $\De_2$.  It is an exercise to show that $f$ maps the boundary of a chamber in $\De_1$ to 
 the boundary of a chamber in $\De_2$.   Therefore $f$ extends to an isomorphism
   from $\De_1$ to    $\De_2$. 
 Next we prove the claim.

Let $\xi, \eta\in \p\De_1$ with $\xi\eta\subset \De_1^{(1)}$.  
 By the definition of $f$,    $f(v)\in \xi'\eta'$  
for any vertex $v\in \xi\eta$.  It suffices to show that if $x, y\in   \xi\eta$ are two vertices 
 such that $x\in \interior(y\xi)$, then 
 $f(x)\in \interior(f(y)\xi')$.
 Since the link $\Link(\De_1, x)$ is a thick spherical building,  there is a an edge $e=xz$
  such that $zx\cup x \xi$ is a geodesic ray and $xz\cap xy=\{x\}$.  We extend the ray
 $zx\cup x \xi$ to obtain a complete geodesic $\xi\xi_1\subset \De_1^{(1)}$.
  We note $x\in \xi\xi_1$ and $y\notin \xi\xi_1$.  The images $\xi'\eta'$ and $\xi'\xi_1'$ of 
$\xi\eta$ and $\xi\xi_1$ are contained in the 1-skeleton of $\De_2$ and
 have a common point at infinity. Since $\De_2$ is $\CAT(-1)$, the  intersection $\xi'\eta'\cap\xi'\xi_1'$
  is  a ray  asymptotic to $\xi'$.  The assumption and the fact that $x\in \xi\eta\cap \xi\xi_1$ imply 
$f(x)\in \xi'\eta'\cap\xi'\xi_1'$.  
 If $f(y)\in \interior(f(x)\xi')$, then $f(y)\in \xi'\eta'\cap\xi'\xi_1'$,
which in turn implies $y\in \xi\eta\cap\xi\xi_1$, contradicting to $y\notin \xi\xi_1$.

\end{proof}

By Lemma  \ref{redaprtowh} and  Lemma  \ref{allha}, to prove Theorem \ref{main} we only need to verify the two 
 conditions  in Lemma \ref{redaprtowh}.

\subsubsection{Proof of Theorem \ref{main}  when $R_1$,  $R_2$ are not right triangles}\label{crieproo}

In this section we shall show that  the conditions in Lemma  \ref{redaprtowh} are satisfied when 
$R_1$,  $R_2$ are not right triangles.

\b{Le}\label{threeev}
{Let $\De$  be a Fuchsian building,  and  
$\xi_1, \xi_2, \xi_3\in \p\De$ be three points such that $\xi_1\xi_2$,  $\xi_2\xi_3$  and   $\xi_3\xi_1$
   all lie in     $ \De^{(1)}$. 
    Then exactly one of the following occurs:\newline
\e{(1)} there exists  a vertex  $v\in \De$
 such that  $\xi_1\xi_2=v\xi_1\cup v\xi_2$,   $\xi_1\xi_3=v\xi_1\cup v\xi_3$  and 
 $\xi_2\xi_3=v\xi_2\cup v\xi_3$;\newline
\e{(2)} $R=(2,3,8)$  and there are  vertices $v_1, v_2,  v_3$ such that 
$\xi_1\xi_2\cap \xi_1\xi_3=v_1\xi_1$, $\xi_2\xi_1\cap \xi_2\xi_3=v_2\xi_2$, 
$\xi_1\xi_3\cap \xi_2\xi_3=v_3\xi_3$;   furthermore,  
 $T:=(v_1, v_2,v_3)$ has three even angles and $S(T)$
is  as shown in   Figure \ref{238} \e{(h)}.}

\end{Le}

\b{proof}
Since $\De$ is $\CAT(-1)$ and $\xi_1\xi_2$, $\xi_2\xi_3$ $\xi_3\xi_1$ are all
 contained in the 1-skeleton,  the intersections 
$\xi_1\xi_2\cap \xi_1\xi_3$, $\xi_1\xi_2\cap \xi_2\xi_3$  and 
$\xi_1\xi_3\cap \xi_2\xi_3$  are all rays. 
Let $v$ be the vertex  with  $\xi_1\xi_2\cap \xi_1\xi_3=v\xi_1$. 
If $v\in \xi_2\xi_3$, then the convexity of geodesic implies 
 $\xi_2\xi_3=v\xi_2\cup v\xi_3$ and so (1) holds. Suppose $v\notin \xi_2\xi_3$. 
    Then there are vertices $v_2\in \interior(v\xi_2)$, $v_3\in \interior(v\xi_3)$ such that 
$\xi_1\xi_2\cap \xi_2\xi_3=v_2\xi_2$  and 
$\xi_1\xi_3\cap \xi_2\xi_3=v_3\xi_3$.   Now $T=(v, v_2, v_3)\subset \De^{(1)}$ is a triangle  with  three even 
 angles.   Proposition  \ref{tri2}
  implies that  $R=(2,3,8)$ and $S(T)$
is as shown in   Figure \ref{238} (h), hence (2) holds.

\end{proof}

\b{remark}\label{prop5.2e}
{In the above proof we came up with a triangle $T\subset \De^{(1)}$ which has three even angles. We then check 
 the triangles listed in Proposition \ref{tri2}  and find out that  only one  
   of them has this property and the chamber has to be $(2,3,8)$.
  We shall use Proposition \ref{tri2}  very often in this paper.  Usually a triangle $T\subset \De^{(1)}$
  is present  which has  a  certain property, for example, \lq\lq has an even angle", \lq\lq one side contains 
 one or two vertices in the interior", \lq\lq has an angle  $>\pi/2$".  When we apply 
Proposition \ref{tri2},  we  check the triangles listed in Proposition \ref{tri2} to see whether any of them 
  or which ones of them have  the stated  property.  This is how we apply Proposition \ref{tri2}.}

\end{remark}

Recall for each vertex $v$ in a Fuchsian building $\De$,  
the link $\Link(\De, v)$  is a generalized $m$-gon for 
 some $m\in \{2,3,4,6,8\}$.  We  define    $m(v)=m$
 if $\Link(\De, v)$  is a generalized $m$-gon.

\b{Le}\label{commint}
{Suppose the following condition holds:\newline
 for  any two geodesics $\xi_1\xi_2, \eta_1\eta_2\subset \De_1$ \e{($\xi_1, \xi_2, \eta_1, \eta_2\in \p\De_1$)}
 that are locally contained in an apartment and are at different sides, their images $\xi'_1\xi'_2$,
 $\eta'_1\eta'_2$ have nonempty intersection. \newline
 If  $R_2$  is  not  a  right   triangle,  then for   any vertex $v\in \De_1$ with $m(v)\not=3$  and any apartment $A\owns v$, 
 the geodesics in ${\mathcal{D}}'_{A, v}$  intersect in a unique vertex $v_A$.}


\end{Le}



\b{proof}We observe  that 
for  any two geodesics $\xi_1\xi_2, \eta_1\eta_2\subset \De_1$ {($\xi_1, \xi_2, \eta_1, \eta_2\in \p\De_1$)}
 that are locally contained in an apartment and are at different sides, 
 the intersection
 $\xi'_1\xi'_2\cap \eta'_1\eta'_2$  is  a single point.  Otherwise $\xi'_1\xi'_2\cap \xi'_1\eta'_i$ is a ray for some $i=1, 2$,
which implies that $\xi_1\xi_2\cap \xi_1\eta_i$ is a ray,  contradicting to the fact  that 
$\xi_1\xi_2$  and   $\eta_1\eta_2$  are at different sides.

  The case  $m(v)=2$ follows from the assumption. 
   Assume $m(v)\ge 4$. 
First suppose there  are  three members  $c_1$, $c_2$, $c_3$ of ${\mathcal{D}}'_{A,v}$ intersecting at one point $w$. 
  Suppose  $c$ is  a member of ${\mathcal{D}}'_{A,v}$  that does not contain $w$.  Then the three intersection points
 $p_i=c_i\cap c$($i=1,2,3$)  are  pairwise distinct. We may assume $p_2$ lies between $p_1$ and $p_3$. 
   Then  $(w, p_1, p_3)\subset \De_2^{(1)}$ is a triangle that is  not the boundary of a chamber, contradicting to 
Proposition \ref{tri2}.

Now suppose  the intersection of any   three members of   ${\mathcal{D}}'_{A,v}$  is empty. 
 Fix some  $c\in {\mathcal{D}}'_{A,v}$. Since $m(v)\ge 4$, there are at least three members $c_1$,
  $c_2$, $c_3$ of   ${\mathcal{D}}'_{A,v}$  that are  distinct from $c$. Let $p_i=c\cap c_i$ ($i=1,2,3$). We may assume $p_2$ lies between 
 $p_1$ and $p_3$. Let $w=c_1\cap c_3$. We obtain a contradiction as above.

\end{proof}

\b{Le}\label{commintmis3}
{Suppose the following condition holds:\newline
 for  any two geodesics $\xi_1\xi_2, \eta_1\eta_2\subset \De_1$ \e{($\xi_1, \xi_2, \eta_1, \eta_2\in \p\De_1$)}
 that are locally contained in an apartment and are at different sides, their images $\xi'_1\xi'_2$,
 $\eta'_1\eta'_2$ have nonempty intersection. \newline
If  $R_2$  is  not  a  right   triangle,  then for   any vertex $v\in \De_1$ with $m(v)=3$  and any apartment $A\owns v$, 
 the geodesics in ${\mathcal{D}}'_{A, v}$  intersect in a unique vertex $v_A$.}


\end{Le}

\b{proof}  Let $\xi_i\eta_i\subset A^{(1)}$ ($\xi_i, \eta_i\in \p A$,  $i=1,2,3$)  be the 
three geodesics    that  pass  through $v$. We may assume 
 $\angle_v(\xi_i,\xi_j)=2\pi/3$ for $i\not=j$. 
Note $\Link(\De_1, v)$ is a thick generalized 3-gon,   $\Link(A,v)$ is an apartment in 
$\Link(\De_1, v)$  and the initial   directions of $v\xi_1$, $v\xi_2$, $v\xi_3$  are of the same type.
 Proposition \ref{oppoti} implies that there is an edge $vv'\subset \De_1$
  such that $v'v\cup v\xi_i$ ($i=1,2,3$) is still a   geodesic. 
We extend $vv'$ to a ray $v\xi\supset vv'$ ($\xi\in \p\De_1$).  Then 
 $\xi\xi_i=v\xi\cup v\xi_i$ for $i=1,2,3$. 

By assumption, $\xi'_1\eta'_1$ and $\xi'_2\eta'_2$ intersect at  some  vertex $w$. 
   Assume  exactly one of $\xi'_1\xi'$,  $\xi'_2\xi'$  contains $w$, say 
$w\in \xi'_1\xi'$ and $w\notin \xi'_2\xi'$.  Then $\xi'\xi'_1=w\xi'\cup  w\xi'_1$,  and  
  $\xi'_2\xi'\cap \xi'_2\eta'_2=x\xi'_2$ 
  and  $\xi'\xi'_2\cap \xi'\xi'_1=y\xi'$  for some $x\in\interior(w\xi'_2)$, $y\in \interior(w\xi')$.
 It follows that $(w, x, y)\subset \De_2^{(1)}$ is a triangle with  two even angles
$\angle_x(w,y)$, $\angle_y(w,x)$.  Proposition \ref{tri2}  implies  $R_2=(2,3,8)$, contradicting to our assumption.
   Therefore either  $w\in\xi'_1\xi'\cap \xi'_2\xi'$ or $w\notin\xi'_1\xi'$,  $w\notin\xi'_2\xi'$.

  Assume  none of $\xi'_1\xi'$,  $\xi'_2\xi'$ contains $w$.   Then there are $x\in \interior(w\xi'_2)$,
  $y\in \interior(w\xi'_1)$ such that $\xi'\xi'_2\cap \xi'_2\eta'_2=x\xi'_2$,
  $\xi'\xi'_1\cap \xi'_1\eta'_1=y\xi'_1$.   Note $y\xi'\cap \xi'_2\eta'_2=\phi$, otherwise 
 if  $p= y\xi'\cap \xi'_2\eta'_2$ then $(w, y,p)\subset \De_2^{(1)}$ is a triangle with an even angle $\angle_y(w,p)$,
  contradicting to Proposition \ref{tri2}.
  In particular,  $x\notin y\xi'\subset \xi'_1\xi'$.
 It follows that   $\xi'\xi'_1\cap \xi'\xi'_2=z\xi'$  for some $z\in \interior(x\xi')$.
Similar argument shows   $z\in \interior(y\xi')$.  Now 
$(z,  y,   w,  x)\subset \De_2^{(1)}$ is a quadrilateral with
  three even angles $\angle_x(w,z)$, $\angle_y(w,z)$ and $\angle_z(x,y)$,  contradicting to 
  Corollary  \ref{threeenenq}.   Therefore   both $\xi'_1\xi'$  and   $\xi'_2\xi'$ contain $w$.

  Assume  $w\notin \xi'_3\eta'_3$.  Then 
$\xi'_1\eta'_1$ and $\xi'_3\eta'_3$ intersect at  some  vertex $w'\not=w$. Then the above argument shows that 
 both $\xi'_1\xi'$  and   $\xi'_3\xi'$ contain $w'$.   In particular, $\xi'_1\xi'$  contains both $w$ and $w'$. 
We have either $w'\in \interior(w\xi')$  or  $w\in \interior(w'\xi')$.
 First suppose $w'\in \interior(w\xi')$.  Then $w'\in\interior(w\eta'_1)$.  It follows that
 $\angle_w(\xi'_2, \eta'_1)=\angle_w(\xi'_2, w')=\angle_w(\xi'_2, \xi')$, which is $\pi$ since
 $\xi'_2\xi'$ contains   $w$. Hence  $\xi'_2\eta'_1=  w \xi'_2\cup w\eta'_1\subset \De_2^{(1)}$.  
    Then $\xi_2\eta_1$ must be a geodesic contained in the 1-skeleton
 of $\De_1$, which is not true since  $\xi_1\eta_1$ and $\xi_2\eta_2$ are contained in the apartment
 $A$.   Similarly one obtains a contradiction   if  $w\in \interior(w'\xi')$.

\end{proof}

Lemma  \ref{redaprtowh0.1}  and  Proposition \ref{klarg4}  imply that  the condition in Lemma \ref{commint}
  and  Lemma \ref{commintmis3} is satisfied when $R_2$ is  not  a right  triangle.   
It follows that condition (1) of Lemma  \ref{redaprtowh}  holds  when $R_2$ is  not  a right  triangle. 
 When  $R_1$, $R_2$ are not right triangles,  Theorem  \ref{main} follows from 
Proposition \ref{twoapi}, Lemma \ref{redaprtowh}  and  Lemma  \ref{allha}.

\b{Prop}\label{twoapi}
{Suppose  condition \e{(1)}  of  Lemma \ref{redaprtowh}
  is satisfied and $R_1,  R_2$ are  not  right  triangles.   Then
condition \e{(2)} of Lemma \ref{redaprtowh} is also satisfied.}

\end{Prop}

 We prove Proposition \ref{twoapi}
in   two lemmas depending on  whether $m(v)=3$.

\b{Le}\label{twoapiev}
{Suppose  condition \e{(1)} of  Lemma \ref{redaprtowh}  is satisfied and $R_1,  R_2$   are  not  right  triangles.   Then
condition  \e{(2)} of  Lemma \ref{redaprtowh}  holds for those $v$ with $m(v)\not=3$.}

\end{Le}

\b{proof}
Let  $v\in \De_1$ be  a  vertex with $m(v)\not=3$  and 
$A_1, A_2\subset \De_1$      two  apartments 
 containing   $v$   with  
 $\Link(A_1, v)\cap \Link(A_{2}, v)$     a  half apartment in $\Link(\De_1, v)$.
Let $a, b \in \Link(\De_1, v)$  be the  two opposite vertices of the half apartment
 $\Link(A_1, v)\cap \Link(A_{2}, v)$,  and $c$ the midpoint of $\Link(A_1, v)\cap \Link(A_{2}, v)$.
 Note $c$ is a  vertex as $m(v)\not=3$. Denote by $c_i$ ($i=1,2$) the vertex in $\Link(A_i,v)$
 opposite to $c$.  The points $\xi_i, \eta_i, \omega_i, \zeta_i\in \p A_i$($i=1,2$)
  are defined as follows:  $a$ is the initial direction of $v\xi_i$,
 $b$ is the initial direction of $v\eta_i$, $c$ is the initial direction of $v\omega_i$ 
 and $c_i$ is the initial direction of $v\zeta_i$.  
Note $\omega_i\zeta_j=v\omega_i\cup v\zeta_j$ and $\zeta_1\zeta_2=v\zeta_1\cup v\zeta_2$.


Note    $\zeta_1\zeta_2$ and $\xi_1\eta_1$ are locally contained in an apartment. 
Since $R_1$ is not  a right   triangle,
 Lemma \ref{noritri}  implies    $\zeta_1\zeta_2$ and $\xi_1\eta_1$  are at different sides. 
 It follows that $\zeta'_1\zeta'_2$ and $\xi'_1\eta'_1$
 are at different sides.  Since $R_2$ is not a right  triangle,  Proposition \ref{klarg4}
  implies  
$\zeta'_1\zeta'_2\cap \xi'_1\eta'_1\not=\phi$.   
Similarly the following holds:  $\zeta'_1\zeta'_2\cap \xi'_2\eta'_2\not=\phi$,  
$\zeta'_2\omega'_i\cap \xi'_1\eta'_1\not=\phi$  and 
$\xi'_2\eta'_2\cap \omega'_1\zeta'_i\not=\phi$ ($i=1,2$).

Let $w=\xi'_1\eta'_1\cap \zeta'_1\omega'_1$.  We notice  $w=v_{A_1}$.  Since $R_2\not=(2,3,8)$, 
Lemma \ref{threeev}  applied to  the three points $\omega'_1,  \zeta'_1,  \zeta'_2$
  shows that  there is some $w'\in \omega'_1\zeta'_1$ such that 
 $\zeta'_1\zeta'_2=w'\zeta'_1\cup w'\zeta'_2$ and $\omega'_1\zeta'_2=w'\omega'_1\cup w'\zeta'_2$.
   Assume    $w'\in \interior(w\zeta'_1)$.   Then $w$ is the only intersection point of 
$\zeta'_i\omega'_1$($i=1,2$) with   $\xi'_1\eta'_1$. It follows that 
$\zeta'_1\zeta'_2$ and $\xi'_1\eta'_1$  are disjoint, contradicting to the conclusion in the second paragraph.
 Similarly $w'\in \interior(w\omega'_1)$  is impossible and we conclude $w'=w$.
In particular, $w\in \zeta'_2\omega'_1$.

We need to prove $w\in \zeta'_2\omega'_2$. If $\omega_2=\omega_1$, then we are done. 
 Suppose $\omega_2\not=\omega_1$, then $\zeta_2\omega_2\cap \zeta_2\omega_1$  is  a ray. 
 It follows that $\zeta'_2\omega'_2\cap \zeta'_2\omega'_1=w_2\zeta'_2$ (for some $w_2\in \zeta'_2\omega'_1$)     is also a ray.
    Assume  $w_2\in \interior(w\zeta'_2)$. 
By the second paragraph $\zeta'_2\omega'_2$  and $ \xi'_1\eta'_1$  intersect at some 
point $w_1\in \xi'_1\eta'_1$.   Now we have a triangle $(w, w_1, w_2)$ with an 
 even angle $\angle_{w_2}(w, w_1)$, contradicting to Proposition \ref{tri2}.
   Therefore $w\in w_2\zeta'_2\subset \zeta'_2\omega'_2$.

Next we prove $w\in \xi'_2\eta'_2$.  By the second paragraph, $\xi'_2\eta'_2$
  and $\omega'_1\zeta'_1$ intersect at some point $x$.  If  $x\in\interior(w\omega'_1)$, 
   then $\xi'_2\eta'_2\cap \zeta'_1\zeta'_2=\phi$, contradicting to the second paragraph.
 Similarly $x\notin\interior(w\zeta'_1)$ and  we conclude $w=x\in \xi'_2\eta'_2$.

\end{proof}

\b{Le}\label{twoapi3}
{Suppose  condition \e{(1)} of  Lemma \ref{redaprtowh}  is satisfied and $R_1,   R_2$   are  not  right  triangles.   Then
condition \e{(2)} of   Lemma \ref{redaprtowh} holds for those $v$ with $m(v)=3$.}

\end{Le}

\b{proof}
Let $v\in \De_1$ be a vertex with $m(v)=3$,  
$A_1, A_2\subset \De_1$   be     two  apartments 
 containing   $v$   with  
 $\Link(A_1, v)\cap \Link(A_{2}, v)$     a  half apartment in $\Link(\De_1, v)$.
Let $a, b \in \Link(\De_1, v)$  be the  two opposite vertices of the half apartment
 $\Link(A_1, v)\cap \Link(A_{2}, v)$,  and   $c, d$ the other two vertices on  
$\Link(A_1, v)\cap \Link(A_{2}, v)$ such that $a$ and $c$ are adjacent.
 Let $c_i, d_i\in \Link(A_i,v)$ ($i=1,2$)  be vertices such that 
$c_i$ and $c$ are opposite, and $d_i$ and $d$  are opposite. 
  The points  $\xi_i,\eta_i, \omega_i, \tau_i, \zeta_i, \phi_i\in \p A_i$ ($i=1,2$) are defined as follows:
    $a$ is the initial direction of $v\xi_i$,
$b$ is the initial direction of $v\eta_i$,
 $c$ is the initial direction of $v\omega_i$,
  $d$ is the initial direction of $v\tau_i$,
   $c_i$  is the initial direction of $v\zeta_i$  and 
 $d_i$ is the initial direction of $v\phi_i$.

Note $d_2$ is opposite to both $d$ and $c_1$.  Hence $\tau_1\phi_2=v\tau_1\cup v\phi_2$
  and $\zeta_1\phi_2=v\zeta_1\cup v\phi_2$.  Let   $w=\omega'_1\zeta'_1\cap \tau'_1\phi'_1$.      
   Then the proof of Lemma \ref{commintmis3}
shows that $\tau'_1\phi'_2=w\tau'_1\cup w\phi'_2$
  and $\zeta'_1\phi'_2=w\zeta'_1\cup w\phi'_2$.
  It follows  that  $\zeta'_1\phi'_2\cap \tau'_1\phi'_1=w$.  
 On the other hand,  $d$ is opposite to both $d_1$ and $d_2$ and we have 
 $\tau_2\phi_i=v\tau_2\cup v \phi_i$.   Apply the arguments in Lemma \ref{commintmis3} 
  again we see $\tau'_2\phi'_i=w\tau'_2\cup w \phi'_i$.  In particular,   $w\in \tau'_2\phi'_2$.
   Similarly   $w\in  \zeta'_2\omega'_2$. It follows that
 $v_{A_1}=w=\zeta'_2\omega'_2\cap \tau'_2\phi'_2=v_{A_2}$.

\end{proof}


\b{Prop}\label{onetoneq}
{Let $\De_1$ and $\De_2$ be two Fuchsian buildings,  and $h:\p\De_1\ra \p\De_2$  a homeomorphism that preserves the combinatorial cross rario almost everywhere.
  If one of $R_1$, $R_2$ is a right  triangle, then so is the other.}

\end{Prop}

\b{proof}
We  assume   $R_1$ is  a right  triangle and $R_2$ is not,  and shall  derive a contradiction from this.
Let $A$ be  an apartment in $\De_1$  and $C$ any chamber in $A$.  Denote the vertices of $C$ by 
 $x, y, z$ such that the angle at $y$ is $\pi/2$.  Denote the geodesics in $A$ that  contain $xz$, 
  $xy$ and $yz$ respectively by $c_1$, $c_2$ and $c$. Let $C'$ be the chamber in $A$ that shares the edge
 $xy$ with $C$. Denote the third vertex of $C'$ by $z'$ and the geodesic in $A$ containing $xz'$ by $c_3$.

It is clear that any two of the 4 geodesics $c$, $c_1$, $c_2$, $c_3$  are at different sides.
Let $c'$ and $c'_i$   ($i=1,2,3$)  be the images of $c$ and $c_i$ respectively.
Lemma \ref{redaprtowh0.1}
implies that any two of the 4 geodesics $c'$, $c'_1$, $c'_2$, $c'_3$
 are  also  at different sides.
Since $R_2$ is not a right   triangle,   Proposition \ref{klarg4}
 implies  any two of the 4 geodesics $c'$, $c'_1$, $c'_2$, $c'_3$
have nonempty intersection.  Then the argument in the proof of Lemma \ref{commint}
  shows that these 4 geodesics have a common point. In particular, the three geodesics
 $c'_1$, $c'_2$ and $c'$,   which  are the images of the three geodesics containing the three edges 
 of the   chamber  $C$, have  a  unique common point $w_C$.  We have also shown 
 the following:   if   two chambers $C_1, C_2\subset A$ share an edge  $e$  and
  the angle  of $C_1$ at one of the two vertices of   $e$  is $\pi/2$, then 
$w_{C_1}=w_{C_2}$.

Suppose two chambers $C_1, C_2\subset A$ share an edge  $e=v_1v_2$  and
  both  angles  of $C_1$ at  the two vertices of   $e$    are $\not=\pi/2$. 
Then at least one of these two angles is $\not=\pi/3$, say the angle at $v_1$ is $\not=\pi/3$.
  There are   $m(v_1)\ge 4$  geodesics contained in $A^{(1)}$ and passing through $v_1$. 
Then the argument in the proof of Lemma \ref{commint}  again 
  shows that  their images in $\De_2$ have a common point.  Since the 
images of the two geodesics containing the  two edges of $C_i$ incident to $v_1$  intersect at 
 $w_{C_i}$, we have $w_{C_1}=w_{C_2}$.   Together with the previous paragraph we conclude
 $w_{C_1}=w_{C_2}$ for any two adjacent chambers $C_1, C_2\subset A$. Since $A$ is chamber connected, 
 $w_{C_1}=w_{C_2}$ for any two  chambers $C_1$, $C_2$  in $A$.  
 It follows that there is a  vertex $w\in \De_2$ such that $\xi'\eta'$ contains $w$ for any
  geodesic $\xi\eta$ contained in the 1-skeleton of $A$.

  Since the map
 $h: \p\De_1\ra \p\De_2$  is a homeomorphism,  it induces 
 a homeomorphism $h\times h: \p\De_1\times \p\De_1\ra 
\p\De_2\times \p\De_2$  satisfying  $h\times h (\p\De_1\times \p\De_1-\{(\xi, \xi): \xi\in \p\De_1\})=
\p\De_2\times \p\De_2-\{(\xi, \xi): \xi\in \p\De_2\}$. 
Let  $E_1=\{(\xi,\eta): \xi,\eta\in \p A \;\text{and}\; \xi\eta\subset A^{(1)}\}$
  and $E_2=\{(\xi', \eta'): \xi',\eta'\in \p\De_2 \; \text{and} \; \xi'\eta'\owns w\}$.
Then $h\times h(E_1)\subset E_2$.  On the other hand, 
  $\overline{E_2}\subset \p\De_2\times \p\De_2-\{(\xi, \xi): \xi\in \p\De_2\}$  and 
$\overline{E_1}\not \subset  \p\De_1\times \p\De_1-\{(\xi, \xi): \xi\in \p\De_1\}$, 
  contradicting to the fact that $h\times h$ is a homeomorphism.


\end{proof}

\subsection{The right  triangle case}\label{rightan}

In this section we prove   Theorem  \ref{main}
    when both $R_1$ and $R_2$ are right  triangles 
but none of them is $(2,3,8)$.  The exceptional case $(2,3,8)$ will be considered 
 in the next section.  We remind the reader that Proposition \ref{tri2} shall be used in the manner 
indicated in Remark \ref{prop5.2e}.

\subsubsection{Types of geodesics}  \label{typeog}

Let $\De$ be a Fuchsian building.     
   Each complete geodesic contained  in $\De^{(1)}$ is a labeled  complex with the 
 induced labeling.
Two complete geodesics contained in $\De^{(1)}$  have the  \e{same type}  if  they are isomorphic as labeled
complexes.    Let us   first  find  all possible types of geodesics.  Since geodesics are contained in apartments and 
 all apartments are isomorphic as labeled complexes, we only need to look at  one apartment.  Let $A$ be a fixed apartment  and 
  $C\subset A$ a fixed chamber.    
Since the action of  Coxeter group preserves the labeling,
  for each complete geodesic $c_1\subset A^{(1)}$, there is a complete geodesic $c_2\subset A^{(1)}$ containing an edge of $C$
   such that $c_1$ and $c_2$ are of the same type.   So there are at most $k$ different types of geodesics.  
Let $c\subset A^{(1)}$ be a complete geodesic   and $v\in c$ a vertex  labeled by $\{i, i+1\}$. 
  Note that if $m(v)$ is even, then the two edges of $c$ incident to $v$  have the same labeling;
 and if $m(v)=3$, then these two edges are labeled by  $\{i\}$ and $\{i+1\}$ respectively. 
From this we obtain a description of all different types of geodesics in $\De$.
If $m(v)=3$ for all vertices, then there is only one type of geodesics  and the edges on each complete geodesic are 
 periodically labeled by $\{1\}, \{2\}$, $\cdots$, $\{k\}$.  Assume there is at least one vertex $v$ with $m(v)\not=3$.  
 We consider the boundary  $\p R$  of $R$ and let  $\p R'$  be the complement in  $\p R$ of 
   all vertices  $v$ with $m(v)\not=3$.
  Let $L$  be a component of  $\p R'$.  Then the edges of $L$ are linearly labeled by $\{i\}, \{i+1\},  \cdots, \{j\}$. 
 Let  $c$ be  a  complete geodesic containing an edge labeled by some $\{l\}$ 
  satisfying    $i\le l\le j$. Then the edges of $c$ are periodically 
    labeled by $\{i\}, \{i+1\}, \cdots, \{j\}, \{j\}, \cdots, \{i+1\}, \{i\}$. 
    It is clear that this establishes  a  1-to-1 correspondence between 
 the types of geodesics and  the components of $\p R'$.

Suppose $\p R'$  has at least two components, or equivalently,  there are at least two vertices  $w$  of $R$ with
  $m(w)\not=3$. 
Now let $v\in \De$ be a fixed vertex. Consider the set  $\mathcal{D}_v$   of 
all  complete geodesics through $v$   that  are 
contained in $\De^{(1)}$.  The above 1-to-1 correspondence between the types of geodesics and  the components of $\p R'$ shows the following: if $m(v)=3$, then all the geodesics in $\mathcal{D}_v$ have the same type;
  if $m(v)\not=3$, then $c_1, c_2\in \mathcal{D}_v$  have the same type  if and only if they make an even angle at $v$. 
In general, if two geodesics $c_1, c_2\subset \De^{(1)}$ intersect at a vertex $v$, they actually make 4 angles at $v$. 
 If $m(v)$ is even, then either all 4 angles are even, or all 4 angles are odd; if $m(v)=3$, then 
 2 of the angles are even and two are odd.   Hence   two intersecting geodesics $c_1, c_2\subset \De^{(1)}$ 
     have the same type if and only if at least one of the angles they make is even.

Now suppose $R$  is  a right  triangle.  If $R\not=(2,3,8)$, then there are three different types of geodesics;
  if $R=(2,3,8)$, then there are only two different types of geodesics.

\subsubsection{Two  disjoint geodesics  at different sides}  \label{twogeod}

Recall that  there are only 6 right  triangles (up to  isometry) 
    that can appear as the chamber of a Fuchsian building:
$(2,8,8)$, $(2,6,8)$,   $(2,6,6)$,  $(2,4,6)$,   $(2,4,8)$,  $(2,3,8)$.


In Lemma \ref{onedif} 
  and  Propositions  \ref{samegeo}  and  \ref{tgeoddif}
we have done some preliminary studies on disjoint geodesics that are at different sides. 
 Here we continue this  analysis.

\begin{figure}[h]
\centering\epsfig{file=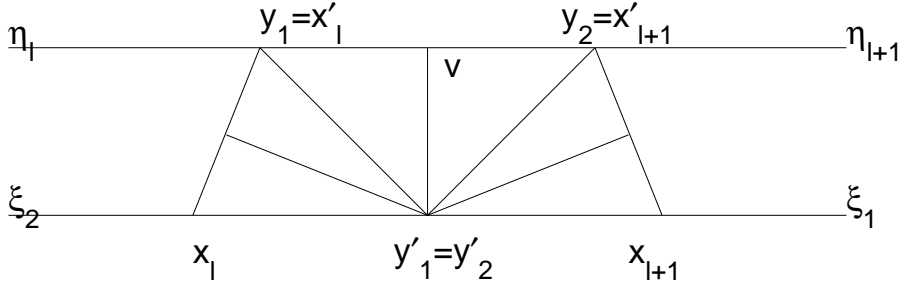,  height=4cm,  width=12cm}\caption{Two disjoint geodesics at different sides}\label{8l1}
\end{figure}

\b{Le}\label{riasamegeo}
{Suppose $R\not=(2,3,8)$ is a right  triangle.
Let  $\xi_1\xi_2, \eta_1\eta_2\subset \De^{(1)}$ \e{($\xi_1$,  $\xi_2$, $\eta_1$, $\eta_2\in \p\De$)}
   be  geodesics that
   are  disjoint and are
   at different sides. 
Let $x_i\in \xi_1\xi_2$, $y_i\in \eta_1\eta_2$  be vertices   such that
  $x_i\eta_i\cap \eta_1\eta_2$,  $y_i\xi_i\cap \xi_1\xi_2$ are rays
  and  $x_i\eta_i\cap  \xi_1\xi_2=\{x_i\}$,   $y_i\xi_i\cap   \eta_1\eta_2=\{y_i\}$.
   If  $y_1'=y_2'$ and  $y_1\not=y_2$,   
  then $R=(2,4,6)$ or $(2,4,8)$, $m(y'_1)=4$, $m(y_1)=m(y_2)=6$ or 8, and 
      there is a  configuration 
        as shown in  Figure \ref{8l1}.}

\end{Le}

\b{proof}
Recall in a generalized  $m$-gon  opposite vertices have the same type if  $m$
  is   even. Since $\Link(\De, y'_1)$ is a generalized  $m(y'_1)$-gon
  and $m(y'_1)$ is even, it follows  from the assumptions that the angle 
 $\angle_{y'_1}(y_1,y_2)$ is even.      Proposition \ref{tri2}  applied to $(y_1, y_2, y'_1)$ shows that 
     $y_1y_2$ contains exactly one vertex
 in the interior,    $m(y_1)=m(y_2)\ge 4$  and 
    $\angle_{y_1}(y'_1,y_2)=\angle_{y_2}(y'_1,y_1)=\frac{\pi}{m(y_1)}$.

First we assume one of $x'_1, x'_2  $ does not lie on $y_1y_2$.  Then there are   $i$ and $j$  with 
$y_{j+1}\in \interior(y_jx'_i)$.  Here,  below  and in Figure  \ref{8l1} the indices $j+1, i+1$, $l+1$ are taken 
  mod 2.  
Notice $\angle_{y_{j+1}}(y'_1, x'_i)\ge \pi-\angle_{y_{j+1}}(y_j, y'_1)\ge 3\pi/4$. 
If $x_i\in y'_1\xi_{j+1}$, then $x_iy_{j+1}=x_iy'_1\cup y'_1 y_{j+1}$ and
 $(x_i, y_{j+1}, x'_i)$ is a triangle with sum of angles $\ge \angle_{y_{j+1}}(y'_1, x'_i)+\pi/8+\pi/8=\pi$, which is 
 impossible. Hence $x_i\in y'_1\xi_j$. In this case $x_iy_j=x_iy'_1\cup y'_1y_j$ and $(x_i, y_j, x'_i)$ is a 
 triangle.  By the definition of $x_i$,  $\angle_{x'_i}(x_i, y_j)$  either $=\pi$ or is an even angle. 
 $\angle_{x'_i}(x_i, y_j)=\pi$ is impossible since otherwise $x_ix'_i\cup x'_iy_j$ and $x_iy'_1\cup y'_1y_j$ 
 would be two distinct geodesic segments from $x_i$ to $y_j$.  Therefore   $(x_i, y_j, x'_i)$
 has an even angle at $x'_i$. Since $y_{j+1}\in \interior(y_jx'_i)$, Proposition \ref{tri2}
  implies $m(y_{j+1})=2$, contradicting to the above observation that $m(y_1)=m(y_2)\ge 4$.  
The contradiction shows 
 $x'_1, x'_2\in y_1y_2$.

Assume   $x'_i\in \interior(y_1y_2)$ for some $i=1,2$.  Then $x'_i$ is the only vertex in the interior
of $y_1y_2$.  There is some $j=1,2$  with 
   $x_i\in y'_1\xi_j$. 
 The uniqueness of the  geodesic  $x_iy_j$ shows $\angle_{x'_i}(x_i, y_j)\not=\pi$.
Hence $\angle_{x'_i}(x_i, y_j)$  is  an even angle.  
 Since $x'_iy_j$ contains no vertex in the interior, Proposition \ref{tri2}  applied to $(x_i, x'_i, y_j)$ 
    implies that $m(y'_1)=2$,  contradicting to the fact that $(y'_1, y_1, y_2)$ has an even angle at $y'_1$. 
 Therefore we conclude $\{x'_1, x'_2\}\subset \{y_1,y_2\}$.

There is some $i=1,2$ with $x'_1=y_i$. 
 The uniqueness of geodesic implies  $x_1\in y'_1\xi_{i+1}$  and 
$\angle_{x'_1}(x_1, y_{i+1})\not=\pi$.  We have $x_1y_{i+1}=x_1y'_1\cup y'_1y_{i+1}$ and  $(x'_1, x_1, y_{i+1})$
has an even angle $\angle_{x'_1}(x_1, y_{i+1})$. Since  $y_iy_{i+1}$ has a vertex $v$ in the interior,
   Proposition \ref{tri2}
       implies the following:
 $R=(2,4,6)$ or $(2,4,8)$;    $m(v)=2$;  $m(y'_1)=4$;  there is a vertex $v_1\in \interior(x_1x'_1)$ with
 $m(v_1)=2$; $(y'_1, x_1, v_1)$, $(y'_1, x'_1, v_1)$, $(y'_1, x'_1, v)$,
  $(y'_1, y_{i+1}, v)$ are chambers and $S(T)$ is the union of these 4 chambers.
 Applying  the same argument to $x'_2$  we obtain  the following:
 $x'_2=y_{i+1}$, $x_2\in \interior(y'_1\xi_i)$; there is a vertex $v_2\in \interior(x_2x'_2)$ 
 such that $(y'_1, x_2, v_2)$ and $(y'_1, x'_2, v_2)$  are chambers.

\end{proof}

\b{Le}\label{samethffl}
{Suppose $R\not=(2,3,8)$ is a right  triangle.
Let  $\xi_1\xi_2, \eta_1\eta_2\subset \De^{(1)}$ \e{($\xi_1$,  $\xi_2$, $\eta_1$, $\eta_2\in \p\De$)}
   be  geodesics that
   are  disjoint and are
   at different sides. 
Let $x_i\in \xi_1\xi_2$, $y_i\in \eta_1\eta_2$  be vertices   such that
  $x_i\eta_i\cap \eta_1\eta_2$,  $y_i\xi_i\cap \xi_1\xi_2$ are rays
  and  $x_i\eta_i\cap  \xi_1\xi_2=\{x_i\}$,   $y_i\xi_i\cap   \eta_1\eta_2=\{y_i\}$.
   Suppose   $y_1'\in\interior(x_1\xi_1)$  and $y_1\in \interior(x_1'\eta_1)$.
Then $R=(2,4,6)$ or $(2,4,8)$,  and there must be a configuration   isomorphic to the one in Figure \ref{8l1}.}

\end{Le}

\b{proof}
   If $x'_2=x'_1$, then   we  are done  by Lemma \ref{riasamegeo}.
   We suppose $x'_2\not=x'_1$.
 By Proposition \ref{samegeo},  $x'_2\in \interior(x'_1\eta_2)$.
Note $m(x'_1)\not=2$, otherwise $\angle_{x'_1}(x_1,y_1)=\angle_{x'_1}(x_1,x'_2)=\pi$, contradicting 
   to  Proposition \ref{samegeo}.
First assume  $x_2\in y'_1\xi_1$. In this case $x_2y_1=x_2y'_1\cup y'_1y_1$ and 
 $(x_2, y_1, x'_2)$  has an even angle $\angle_{x'_2}(x_2, y_1)$. Since 
 $y_1x'_2$ contains the vertex $x'_1$ in the interior, Proposition \ref{tri2}
  implies  $m(x'_1)=2$,  a contradiction.

Note Proposition \ref{samegeo}  implies $x_2\not=x_1$. Next we assume $x_2\in \interior(x_1y'_1)$. 
 Since $(y'_1, x_1, y_1)$  has  an even angle  at $y'_1$,  Proposition \ref{tri2}
 implies   $m(x_2)=2$  and $\angle_{x_2}(x_1, x'_1)=\pi/2$.
  If $\angle_{x_2}(x_1, x'_2)=\pi$, then $\angle_{x_2}(y'_1, x'_2)=\pi$, contradicting to 
 Proposition \ref{samegeo}.  Hence  $\angle_{x_2}(x_1, x'_2)=\pi/2$,
  which implies $\angle_{x_2}(x'_2, x'_1)=\pi$ and so $x'_1x'_2=x'_1x_2\cup x_2x'_2$ is not contained in
   $\eta_1\eta_2$,
 a contradiction. We conclude $x_2\in \interior(x_1\xi_2)$.

Since $(y'_1, y_1, x_1)$ is a triangle with an even angle
 $\angle_{y'_1}(x_1,y_1)$    and the side $x_1y_1$ contains
  the   vertex $x'_1$ with $m(x'_1)\not=2$ in the interior, 
Proposition \ref{tri2}
  implies   $R=(2,4,6)$ or $(2,4,8)$,
 $m(x'_1)=4$ and $\angle_{x_1}(x'_1, y'_1)=\pi/6$ or $\pi/8$.
 It follows  that  $\angle_{x_1}(x_2, x'_1)=5\pi/6$ or $7\pi/8$. 
 Consider the quadrilateral $Q=(x_1, x'_1,x'_2,x_2)$:
 $\angle_{x'_1}(x_1, x'_2)\ge \pi/2$ since it is even and $m(x'_1)=4$; 
  $\angle_{x'_2}(x_2, x'_1)\ge \pi/4$   since it is  even. 
  Proposition \ref{quar11}  implies   $n(Q)\le 2$, 
   which is impossible since there are at least 
 5 or 7 chambers in $ S(Q)$ incident to $x_1$.

\end{proof}

\b{Le}\label{tanej}
{Suppose $R\not=(2,3,8)$ is a right  triangle.
Let  $\xi_1\xi_2, \eta_1\eta_2\subset \De^{(1)}$ \e{($\xi_1$,  $\xi_2$, $\eta_1$, $\eta_2\in \p\De$)}
   be  geodesics that
   are  disjoint and are
   at different sides.    
  Let $x_i\in \xi_1\xi_2$, $y_i\in \eta_1\eta_2$  be vertices   such that
  $x_i\eta_i\cap \eta_1\eta_2$,  $y_i\xi_i\cap \xi_1\xi_2$ are rays
  and  $x_i\eta_i\cap  \xi_1\xi_2=\{x_i\}$,   $y_i\xi_i\cap   \eta_1\eta_2=\{y_i\}$.
Suppose   $y_1'\in\interior(x_1\xi_1)$  and $y_1=x'_1$.
Then $R=(2,4,6)$ or $(2,4,8)$,  and there must be a configuration   isomorphic to the one  in 
   Figure \ref{8l1}.}

\end{Le}

\b{proof}
The  lemma   follows from
  Lemma  \ref{riasamegeo}
  if $y'_2=y'_1$. We assume $y'_2\not=y'_1$.
Proposition \ref{samegeo}  implies $y'_2\in \interior(y'_1\xi_2)$. 
We first assume $y'_2\in \interior(y'_1x_1)$.  Note $(y'_1, x_1, y_1)$ is a triangle with 
 an even angle at $y'_1$ and the side $y'_1x_1$ contains the  vertex  $y'_2$ in 
 the interior.   Proposition \ref{tri2}   implies   $m(y'_2)=2$, which in turn implies
 $\angle_{y'_2}(y'_1, y_2)=\angle_{y'_2}(x_1, y_2)=\pi$, contradicting to Proposition \ref{samegeo}.
 Hence $y'_2\in x_1\xi_2$.

Now assume   $y'_2=x_1$.  The uniqueness of geodesic implies that $y_2\in y_1\eta_2$.
 Proposition \ref{samegeo}  implies $y_2\not=y_1$. 
The segment $x_1y_1$ is the side of $(y'_1, x_1, y_1)$  opposite to an  even angle.
Proposition \ref{tri2}
     implies $x_1y_1$ contains exactly one vertex in the interior.
 On the other hand,  $(x_1, y_1, y_2)$ has an even angle at $y_1$. 
By   Proposition \ref{tri2},  $R=(2,4,6)$ or $(2,4,8)$,
 $y_1y_2$ contains exactly one vertex $v$ in the interior and  $m(v)=2$, and 
$x_1y_2$  also  contains exactly one vertex $v'$ in the  interior  and  $m(v')=4$.
 Now let $x_2\in \xi_1\xi_2$ be a vertex such that $x_2\eta_2\cap \eta_1\eta_2=x'_2\eta_2$ is  a ray
  and  $x_2\eta_2\cap \xi_1\xi_2=\{x_2\}$.
 By Proposition \ref{samegeo}   and Lemma  \ref{riasamegeo} we may assume $x'_2\in \interior(y_1\eta_2)$. 
The argument in the   first  paragraph shows $x'_2\not=v$. 
   Assume  $x'_2=y_2$.  The uniqueness of geodesic shows $x_2\in \interior(x_1\xi_1)$.  
 The triangle $(x_1, x_2, y_2)$ has an even angle at $x_1$ and the side $x_1y_2$ contains 
  the  vertex $v'$ with $m(v')=4$, contradicting to Proposition \ref{tri2}.
 Therefore we must have $x'_2\in \interior(y_2\eta_2)$.  
The lemma follows from Lemma \ref{samethffl}
if  $x_2\in \interior(x_1\xi_2)$.  Proposition \ref{samegeo}  implies $x_2\not=x_1$. 
Therefore $x_2\in y'_1\xi_1$.  In this case $x_2y_1=x_2y'_1\cup y'_1y_1$,  
    $(x_2, y_1, x'_2)$  has an even angle at $x'_2$ and the side $y_1x'_2$ contains two vertices 
 $v$, $y_2$ in the interior, contradicting to Proposition \ref{tri2}.

Finally we assume  $y'_2\in \interior(x_1\xi_2)$. Proposition \ref{samegeo}
  implies $y_2\not=y_1$. The lemma follows from Lemma \ref{samethffl} if $y_2\in \interior(y_1\eta_1)$.
 We now assume $y_2\in \interior(y_1\eta_2)$.   Then we have a quadrilateral $Q=(x_1, y_1,y_2, y'_2)$.
We shall show such a quadrilateral does  not exist. We only give the proof for $R=(2,4,6)$.
 The other cases can be handled
similarly. Recall $(y'_1, x_1, y_1)$ has an even angle at $y'_1$. 
Proposition \ref{tri2} 
 implies that $m(x_1)=m(y_1)$  and 
$\angle_{y_1}(x_1, y'_1)=\angle_{x_1}(y_1, y'_1)=\pi/4$  or $\pi/6$.
The quadrilateral
  $Q$   has  even angles at $y_1$ and $y'_2$.   Notice  
    $\angle_{x_1}(y_1, y'_2)\ge \pi-\angle_{x_1}(y_1, y'_1)$
  and $\angle_{y_1}(x_1,y_2)\ge 2\angle_{y_1}(x_1, y'_1)=2\angle_{x_1}(y_1, y'_1)$.  If $\angle_{x_1}(y_1, y'_1)=\pi/6$,
   then  $\angle_{x_1}(y_1, y'_2)=5\pi/6$ and  Proposition \ref{quar11}  implies
 $n(Q)\le 4$, which is impossible because there 
 are at least 5 chambers in $S(Q)$ incident to $x_1$. If $\angle_{x_1}(y_1, y'_1)=\pi/4$, then 
  $\angle_{x_1}(y_1, y'_2)=3\pi/4$ and  Proposition \ref{quar11}  implies
  $n(Q)\le 3$. It follows that 
$S(Q)$  is the union of the 3 chambers in $S(Q)$
incident to $x_1$.  But such a union has only one even angle, contradiction.

\end{proof}


\b{Prop}\label{twogdisatd}
{Suppose $R\not=(2,3,8)$ is a right  triangle. 
Let  $\xi_1\xi_2, \eta_1\eta_2\subset \De^{(1)}$ \e{($\xi_1,  \xi_2, \eta_1, \eta_2\in \p\De$)}
   be  geodesics that
   are  disjoint and are
   at different sides.  Then $R=(2,4,6)$ or $(2,4,8)$  and 
there must  be a configuration   isomorphic to the one  in Figure \ref{8l1}.    
   In particular,
    $\xi_1\xi_2$  and  $\eta_1\eta_2$  are of different types.}

\end{Prop}

\b{proof}
The proof is similar to that  of Proposition \ref{klarg4}.  Lemma \ref{onedif}
 implies there are vertices 
$x_i\in \xi_1\xi_2$, $y_i\in \eta_1\eta_2$   such that
  $x_i\eta_i\cap \eta_1\eta_2$,  $y_i\xi_i\cap \xi_1\xi_2$ are rays
  and  $x_i\eta_i\cap  \xi_1\xi_2=\{x_i\}$,   $y_i\xi_i\cap   \eta_1\eta_2=\{y_i\}$.
By Lemmas \ref{riasamegeo}, \ref{samethffl}  and \ref{tanej},   $R=(2,4,6)$ or $(2,4,8)$, and  
there is a configuration   isomorphic to the one  in Figure \ref{8l1}  if any one of the following occurs:  
Proposition \ref{samegeo} (2), Proposition \ref{tgeoddif} (1), (2). 
 So we can assume only  Proposition \ref{samegeo} (1)
 and Proposition \ref{tgeoddif} (3)  can occur.  Then the proof of Proposition \ref{klarg4}
  shows  that 
there is a quadrilateral $Q$ with one side $xy$ such that the angles at $x$ and $y$ are even
  and $xy$ contains two vertices in the interior.  Corollary \ref{q3in2ev}
implies $R=(2,3,8)$,    contradicting to our assumption.

\end{proof}

Let us make some  observations  about the configuration in Figure \ref{8l1}.
Notice $\angle_v(\eta_1,y'_1)=\angle_v( \eta_2,y'_1)=\pi/2$ and 
$\angle_{y'_1}(\xi_1, v)=\angle_{y'_1}(\xi_2,v)=3\pi/4$.  It follows that $vy'_1$ is the shortest geodesic between 
 $\xi_1\xi_2$ and $\eta_1\eta_2$.  Consider the quadrilateral $Q=(x_l, x'_l,  x'_{l+1}, x_{l+1})$. 
     Notice $Q$ is uniquely determined by 
 $\xi_1\xi_2$ and the segment $vy'_1$:  the geodesic segment in $\Link(\De, y'_1)$ from the initial direction of 
 $y'_1v$ to the initial direction of $y'_1\xi_1$ consists of 3 edges, and these 3 edges determine 3 chambers incident to $y'_1$;
 similarly  the geodesic segment in $\Link(\De, y'_1)$ from the initial direction of 
 $y'_1v$ to the initial direction of $y'_1\xi_2$  determine 3 chambers; the union of these 6 chambers is a  disk whose boundary is 
 $Q$.  We denote this $Q$ by $Q(   \xi_1\xi_2,   \eta_1\eta_2)$.

\subsubsection{Geodesics in the same apartment}\label{geitsaap}

The eventual  goal of  this section is to show

\b{Prop}\label{serk2po2}
{Let $\De_1$,  $\De_2$ be two Fuchsian buildings, and $h:\p\De_1\ra \p\De_2$  a homeomorphism that preserves the combinatorial cross ratio almost everywhere.
  Let $A$ be an apartment in $\De_1$ and $v\in A$ a vertex.
 If both $R_1$ and $R_2$ are right  triangles different from $(2,3,8)$, then 
the geodesics in ${\mathcal{D}}'_{A,v}$   intersect in a unique vertex of $\De_2$.}

\end{Prop}

  One  first has  to show that any two geodesics in ${\mathcal{D}}'_{A,v}$ have 
  nonempty intersection (Proposition \ref{sec8la4}).




\b{Le}\label{sec8la1}
{Suppose $R\not=(2,3,8)$ is a right  triangle.
 Let $c, c_1, c_2\subset \De^{(1)}$ be complete geodesics such that  
  any two of them   are at different sides.
  If $c_1$ and $c_2$ have the same type  and $c\cap c_1=\phi$, then $c\cap c_2=\phi$.}

\end{Le}

\b{proof}  Assume  $c\cap c_2\not=\phi$.
Since $c\cap c_1=\phi$, Proposition \ref{twogdisatd}   implies  $R=(2,4,6)$ or $(2,4,8)$ and 
there is a configuration as shown in Figure \ref{8l1},
  where either $c_1=\xi_1\xi_2$,  $c=\eta_1\eta_2$ or $c=\xi_1\xi_2$,  $c_1=\eta_1\eta_2$.
 Proposition \ref{twogdisatd}  also  implies $c_1\cap c_2\not=\phi$ because 
$c_1$ and $c_2$ have the same type and are at different sides.
Let  $p=c_2\cap \eta_1\eta_2$,    $q=c_2\cap \xi_1\xi_2$.  Then $pq\subset c_2\subset \De^{(1)}$.
  We will show $p\in v\eta_{l+1}$ is impossible. Since the configuration is symmetric about $vy'_1$, $p\in v\eta_l$
  is also impossible,  and we have a contradiction. We shall only write down the proof  for  $R=(2,4,6)$,
 the case  $R=(2,4,8)$  can be similarly handled. 

First assume  $p\in \interior(y_2\eta_{l+1})$.   
Then $q\in \interior(x_{l+1}\xi_1)$ is impossible, since otherwise $px_{l+1}=py_2\cup y_2x_{l+1}$ and $(p,x_{l+1}, q)$ has 
  an angle $\angle_{x_{l+1}}(p, q)=\pi-\angle_{x_{l+1}}(y_2, y'_1)=5\pi/6$, which is impossible.
  Note $q=x_{l+1}$ cannot occur since otherwise $c_2\cap \eta_1\eta_2$ contains a nontrivial segment $py_2$.
 Note $\angle_{y'_1}(y_2,x_l)=\pi$. 
 If $q\in y'_1\xi_2$, then $y_2q=y_2y'_1\cup y'_1q$ and  
$(p,y_2, q)$ has an angle $\angle_{y_2}(p,q)\ge \angle_{y_2}(p,x_{l+1})-\angle_{y_2}(x_{l+1},y'_1)=\pi-\pi/6=5\pi/6$,
 a contradiction.

Now assume $p=y_2$. 
The angle $\angle_{x_{l+1}}(y_2, \xi_1)=5\pi/6$ implies $q\notin   \interior(x_{l+1}\xi_1)$. 
Note that if two intersecting geodesics  are at different sides, then 
none of the 4 angles they 
  make   can be  $\pi$. 
It follows that $q\not=x_{l+1}, y'_1$ since $\angle_{y_2}(x_{l+1}, \eta_{l+1})=\angle_{y'_1}(y_2, \xi_2)=\pi$. 
  The fact $y'_1\in y_2\xi_2$ implies $q\notin \interior(y'_1\xi_2)$ since 
  $c_2\cap \xi_1\xi_2$ is a single point.

Finally we assume $p=v$.  Note $m(v)=2$ and $m(y'_1)=4$.
If $q\not=y'_1$, then $(v, y'_1, q)$ is a triangle
 with $\angle_v(y'_1,q)=\pi/2$ and $\angle_{y'_1}(v,q)=3\pi/4$, which is impossible.
 Hence $q=y'_1$.   It  follows  that  the type of $c_2\supset vy'_1$ is different from the types of
 $\xi_1\xi_2$ and $\eta_1\eta_2$, contradicting to the assumption.

\end{proof}

\b{Le}\label{sec8la2.5}
{Let $\De$  be a Fuchsian building with chamber $R\not=(2,3,8)$,  
 $v_0\in \De$   a  vertex  and $\xi_1, \xi_2, \eta_1, \eta_2, \eta_3\in \p \De$.
Suppose the following conditions are satisfied:\newline
\e{(1)} $\xi_1\xi_2\subset \De^{(1)}$;\newline
\e{(2)} $\eta_i\eta_j=v_0\eta_i\cup v_0\eta_j\subset \De^{(1)}$  
for all   $1\le i,j\le 3$,  $i\not=j$;\newline
\e{(3)}  $\xi_1\xi_2$ and $\eta_i\eta_j$ are at different sides for 
all   $1\le i,j\le 3$,  $i\not=j$.\newline
 Then  $v_0\in  \xi_1\xi_2$. }

\end{Le}

\b{proof}
First assume $\xi_1\xi_2\cap \eta_i\eta_j\not=\phi$ for some $i\not=j$, say,
 $\xi_1\xi_2\cap \eta_1\eta_2\not=\phi$. Lemma \ref{diinet1}
implies  $\xi_1\xi_2\cap \eta_1\eta_2$ is some vertex $v'\in \De$. If $v'=v_0$, then we are done. 
 Suppose  $v'\not=v_0$. We may assume $v'\in \interior(v_0\eta_1)$. 
 Since $\eta_1\eta_i=v_0\eta_1\cup v_0\eta_i$ for $i=2,3$,  Lemma \ref{diinet1}
 also implies $\xi_1\xi_2\cap v_0\eta_i=\phi$ for $i=2,3$. It follows that
 $\xi_1\xi_2\cap \eta_2\eta_3=\phi$.  
Now $\xi_1\xi_2$ and $\eta_2\eta_3$ are disjoint  and are  at different sides,  and 
$v'\eta_i\cap \eta_2\eta_3=v_0\eta_i$ ($i=2,3$) with $v'\in \xi_1\xi_2$,  contradicting to 
Proposition \ref{samegeo}.

 Now we assume $\xi_1\xi_2\cap \eta_i\eta_j=\phi$ for all  $i\not=j$.
 Proposition \ref{twogdisatd}  applied 
to $\xi_1\xi_2$ and $\eta_1\eta_2$  shows that   $R=(2,4,6)$ or $(2,4,8)$, 
  and there is a configuration   isomorphic to the one  in Figure \ref{8l1}. 
 Consider the perpendicular $xy$ ($x\in \eta_1\eta_2$, $y\in \xi_1\xi_2$)  
between $\eta_1\eta_2$
  and $\xi_1\xi_2$.
We claim $v_0=x$.  Suppose $v_0\not=x$, say $v_0\in \interior(x\eta_2)$.  
 Since $\eta_1\eta_i=v_0\eta_1\cup v_0\eta_i$ for $i=2,3$  and  $\De$ is  $\CAT(-1)$, $d(z, \xi_1\xi_2)>d(x,y)$  holds 
   for all $z\in v_0\eta_i$, $i=2,3$.  It  follows  that   $d(\eta_2\eta_3, \xi_1\xi_2)>d(x,y)$.
Here we have a contradiction since  Proposition \ref{twogdisatd} 
  applied to $\xi_1\xi_2$ and $\eta_2\eta_3$ shows that $d(\eta_2\eta_3, \xi_1\xi_2)=d(x,y)$.
  Hence $v_0=x$.   
Notice $xy$ is the perpendicular between $\xi_1\xi_2$ and $\eta_i\eta_j$ for all $i\not=j$.
By Proposition \ref{twogdisatd}  $m(v_0)=2$ or $4$. We consider these two cases separately.

First assume $m(v_0)=2$. Then $m(y)=4$.  Consider 
the quadrilaterals $Q(\xi_1\xi_2,   \eta_i\eta_j)$, $i\not=j$.  
 By the observation made after Proposition \ref{twogdisatd},  
$Q(\xi_1\xi_2,   \eta_i\eta_j)$
   is uniquely determined by $\xi_1\xi_2$ and $yx$. 
 It  follows  that    $Q(\xi_1\xi_2, \eta_1\eta_2)=Q(\xi_1\xi_2,  \eta_1\eta_3)$,
  which is not true since $Q(\xi_1\xi_2,  \eta_1\eta_3)$ contains an edge   from   $v_0\eta_3$ 
  while $Q(\xi_1\xi_2,  \eta_1\eta_2)$ does not.

Now assume $m(v_0)=4$. Then $m(y)=2$.  Let $y_1, y_2\in \xi_1\xi_2$ be the two vertices on $\xi_1\xi_2$ adjacent to $y$,
  and $\sigma,  \sigma_i\in \Link(\De, v_0)$ ($i=1,2$) the initial directions of $v_0y$, $v_0y_i$ respectively.
 Also let $\omega_j\in  \Link(\De, v_0)$ be the  initial direction  of $v_0\eta_j$. 
  Proposition \ref{twogdisatd}  implies that for each $j$, $1\le j\le 3$,  the segment 
 $\sigma\omega_j\subset \Link(\De, v_0)$  consists of three edges and the initial edge is $\sigma\sigma_1$ or 
$\sigma\sigma_2$.   Hence there is a map $f: \{\eta_1, \eta_2, \eta_3\}\ra \{\sigma\sigma_1, \sigma\sigma_2\}$,
  where $f(\eta_j)$ is the initial edge of the segment $\sigma\omega_j$. 
Proposition \ref{twogdisatd}  also  implies  that  $f$ is injective, which is clearly not true since
$ \{\sigma\sigma_1, \sigma\sigma_2\}$  contains only two elements  while $\{\eta_1, \eta_2, \eta_3\}$  contains three elements.

\end{proof}

\b{Le}\label{sec8la2}
{Suppose $R_1,  R_2\not=(2,3,8)$ are right  triangles.
Let $\xi_1\xi_2,  \eta_1\eta_2\subset \De^{(1)}_1$ \e{($\xi_1,  \xi_2,  \eta_1, \eta_2\in \p\De_1$)} 
 be two  geodesics locally contained in an apartment.  Let $v=\xi_1\xi_2\cap \eta_1\eta_2$.
  If  $m(v)=2$, then  $\xi'_1\xi'_2\cap \eta'_1\eta'_2\not=\phi$.  }

\end{Le}

\b{proof}
Pick an edge $e=vv'$ such that $v'\eta_i=v'v\cup v\eta_i$ ($i=1,2$). We 
extend the edge $vv'$ to obtain a ray $v\eta_3$, where $\eta_3\in \p\De_1$.  By construction,
 $\eta_i\eta_j=v\eta_i\cup v\eta_j\subset \De_1^{(1)}$.  It follows that $\eta'_i\eta'_j\subset \De_2^{(1)}$.
  Since $R_2\not=(2,3,8)$,   Lemma \ref{threeev}  implies  
  there is a vertex $w\in \De_2$ such that 
 $\eta'_i\eta'_j=w\eta'_i\cup w\eta'_j$.   On the other hand, since $m(v)=2$, 
$\xi_1\xi_2$ and $\eta_i\eta_j$  make a right angle and are locally contained in an apartment.  
Lemma \ref{ssufdas}
  implies that 
$\xi_1\xi_2$ and $\eta_i\eta_j$ are at different sides. Hence 
$\xi'_1\xi'_2$ and $\eta'_i\eta'_j$ are  also   at different sides.
Now Lemma \ref{sec8la2.5}
implies $w\in \xi'_1\xi'_2$.  Consequently, $w\in \xi'_1\xi'_2\cap \eta'_1\eta'_2$.

\end{proof}

\b{Le}\label{sec8la3}
{Suppose $R_1$ and $R_2$ are right  triangles.
   Let $\xi_1\xi_2, \eta_1\eta_2\subset \De_1^{(1)}$ be two geodesics 
  that have nonempty intersection. If  at least one of the angles that 
$\xi_1\xi_2$    and $ \eta_1\eta_2$  make   is  even, then 
  $\xi'_1\xi'_2$    and $ \eta'_1\eta'_2$  have the same type.}

\end{Le}

\b{proof}
Let $v\in \xi_1\xi_2\cap \eta_1\eta_2$ be  a vertex. We  may assume 
$\angle_{v}(\xi_1, \eta_1)$ is even. 
 Since $\Link(\De_1, v)$ is a thick spherical building,  
Proposition \ref{oppoti2}
 implies there is an edge $e$ incident to $v$ such that $e$ makes an angle $\pi$ with both 
$v\xi_1$ and $v\eta_1$. We extend $e$ to obtain a geodesic ray $v\omega$.  Then 
 $\omega\xi_1=v\omega \cup v\xi_1$   and $\omega\eta_1=v\omega \cup v\eta_1$. 
 Consider the geodesics $\xi'_1\xi'_2$,  $\eta'_1\eta'_2$, $\omega'\xi'_1$ and 
$\omega'\eta'_1$ in $\De_2$.  We note if two geodesics in $\De_2$ intersect in a ray,
  then they have the same type. Now $\xi'_1\xi'_2\cap \omega'\xi'_1$, $\omega'\xi'_1\cap \omega'\eta'_1$,
   $\omega'\eta'_1\cap  \eta'_1\eta'_2$ are all rays. So the 4 geodesics  have the same type. In particular,
  $\xi'_1\xi'_2$    and $ \eta'_1\eta'_2$  have the same type.

\end{proof}

\b{Cor}\label{corse8c1}
{Suppose $R_1,  R_2$ are right  triangles.
 Let $\xi_1\xi_2, \eta_1\eta_2, \zeta_1\zeta_2   \subset \De_1^{(1)}$ be  three geodesics 
 such that any two of  $\xi'_1\xi'_2, \eta'_1\eta'_2, \zeta'_1\zeta'_2$ 
 have nonempty intersection. If $\xi_1\xi_2, \eta_1\eta_2, \zeta_1\zeta_2$ pairwise  have different types, then 
 $\xi'_1\xi'_2\cap \eta'_1\eta'_2\cap  \zeta'_1\zeta'_2=\phi$.}

\end{Cor}

\b{proof}
 Note that  there are  at most  two types of  geodesics through a fixed vertex.
Suppose $\xi'_1\xi'_2\cap \eta'_1\eta'_2\cap  \zeta'_1\zeta'_2$ contains a vertex $w$. 
Then two of   $\xi'_1\xi'_2$,  $\eta'_1\eta'_2$, $\zeta'_1\zeta'_2$,
 say,    $\xi'_1\xi'_2$  and   $\eta'_1\eta'_2$ have the same type.  
It follows that  at least one of the angles that $\xi'_1\xi'_2$  and   $\eta'_1\eta'_2$ make
  is even.  
     Lemma \ref{sec8la3} then
  implies  $\xi_1\xi_2$  and  $\eta_1\eta_2$ have the same type,
  contradicting to the assumption.

\end{proof}

\b{Le}\label{sec8la3.9}
{Let $\De$ be a Fuchsian building with $R=(2,4,6)$ or $(2, 4,8)$.  Suppose 
 $\xi_1\xi_2, \eta_1\eta_2\subset \De^{(1)}$ \e{($\xi_1, \xi_2, \eta_1,  \eta_2\in \p\De$)}
    are  disjoint  and  are at different sides. Assume   $\eta_1\eta_2$ contains vertices 
 $v$ with $m(v)=2$.  Let $pq$ \e{($p\in \eta_1\eta_2$, $q\in \xi_1\xi_2$)} be the perpendicular between $\eta_1\eta_2$  and 
$\xi_1\xi_2$.    If $\xi_3\xi_4\subset \De^{(1)}$ \e{($\xi_3, \xi_4\in \p\De$)}
  is a geodesic such that $\xi_3\xi_4\cap \xi_1\xi_2\not=\phi$,   $\xi_3\xi_4\cap \eta_1\eta_2=\phi$, 
 and      that  $\xi_3\xi_4$ and $\eta_1\eta_2$ are at different sides, then
  $q\in \xi_3\xi_4$  and $pq$ is also the perpendicular between  $\eta_1\eta_2$  and  
   $\xi_3\xi_4$.}

\end{Le}

\b{proof} We shall only write down the proof for $R=(2,4,8)$,
  the case for $R=(2,4,6)$  is similar. 
Suppose $q\notin \xi_3\xi_4$.   We may assume $\xi_3\xi_4\cap \xi_1\xi_2=xy\subset \interior(q\xi_2)$ with 
 $x\in \interior(q\xi_2)$ and $y\in x\xi_2\cup \{\xi_2\}$.
   Let $p'q'$ ($p'\in \eta_1\eta_2$, $q'\in \xi_3\xi_4$)  be the perpendicular 
between   $\eta_1\eta_2$  and  $ \xi_3\xi_4$.   
By Proposition \ref{twogdisatd}  and  our assumption,     $m(p)=m(p')=2$  and $\angle_{q'}(p', \xi_j)=3\pi/4$ 
($j=3,4$). Note $\angle_{y_2}(\eta_1, q)=\angle_{x_2}(y_1, \xi_2)=7\pi/8$.  
    If  $p'\in \interior(y_2\eta_1)$,  then $Q=(q', p', y_2, x) $ has 
angle sum $\ge 3\pi/4+\pi/2+7\pi/8+\pi/8>2\pi$, impossible.
If $p'=p$, then $q'q=q'p'\cup pq$ since $m(p)=2$; hence  $T=(q', q, x)$ has angle sum 
  $\ge 3\pi/4+3\pi/4+\pi/8>\pi$, again impossible.
Now assume $p'\in \interior(y_1\eta_2)$. If  $x\in \interior(x_2\xi_2)$, then $(q', p', x_2, x)$ has angle sum $\ge 3\pi/4+\pi/2+7\pi/8+\pi/8>2\pi$,   impossible. If $x=x_2$, then $(q', p', x_2)$ has angle sum 
$\ge  3\pi/4+\pi/2+\pi/8>\pi$, also impossible.  Hence  $q\in \xi_3\xi_4$.  
Proposition \ref{twogdisatd}  implies $d(\eta_1\eta_2, \xi_3\xi_4)=d(\eta_1\eta_2, \xi_1\xi_2)=d(p,q)$.
  So  $pq$ realizes the distance between  $\eta_1\eta_2$  and $\xi_3\xi_4$. 
  Since $\De$ is $\CAT(-1)$,  $pq$  must be the unique perpendicular between
$\eta_1\eta_2$  and $\xi_3\xi_4$.


\end{proof}

\b{Le}\label{sec8la3.92}
{Let $\De$ be a Fuchsian building with $R=(2,4,6)$ or $(2, 4,8)$.  Suppose 
 $\xi_1\xi_2, \eta_1\eta_2\subset \De^{(1)}$ \e{($\xi_1, \xi_2, \eta_1,  \eta_2\in \p\De$)}  
    are  disjoint    and  are at different sides.  Assume  $\eta_1\eta_2$ contains vertices 
 $v$ with $m(v)=2$.  Let $pq$ \e{($p\in \eta_1\eta_2$, $q\in \xi_1\xi_2$)}  be the perpendicular between 
$\eta_1\eta_2$  and  
$\xi_1\xi_2$, and $p_1, p_2\in \eta_1\eta_2$  the two vertices on $\eta_1\eta_2$ 
  adjacent to $p$.
   If $\eta_3\eta_4\subset \De^{(1)}$ \e{($\eta_3, \eta_4\in \p\De$)}
  is a geodesic such that $\eta_3\eta_4\cap \eta_1\eta_2=v$  is  a  vertex, $\eta_3\eta_4\cap \xi_1\xi_2=\phi$, 
 and $\eta_3\eta_4$ and $\xi_1\xi_2$ are at different sides, then  $v=p_1$ or $p_2$.
  Furthermore,   if
  $p'q'$ \e{($p'\in \eta_3\eta_4$, $q'\in \xi_1\xi_2$)} is the perpendicular between
$\eta_3\eta_4$  and 
 $\xi_1\xi_2$,
  then $q'=q$  and $(v,  q, p')$ is the boundary of  a  chamber.}

\end{Le}

\b{proof}
An argument similar to the proof of Lemma \ref{sec8la3.9}  shows that $q'=q$  and 
 $v=p_1$  or  $p_2$.   
   Let  $Q=(p,q,p', v)$.    $Q$ has  right angles at $p$ and $p'$  and even angles 
 at  $q$ and $v$:   
  $\angle_q(p, p')$ is even because $m(q)=4$ and $\angle_q(p, \xi_1)=\angle_q(p', \xi_1)=3\pi/4$;
  the angle at $v$ is even   since  the assumptions and Proposition \ref{twogdisatd} 
  imply that $\eta_1\eta_2$ and $\eta_3\eta_4$ have the same type. 
It follows that  the angle sum $\ge 7\pi/4$ or $\ge 11\pi/6$ depending   on  
whether $R=(2,4,8)$ or $R=(2,4,6)$. 
Proposition \ref{quar11}  then implies  $n(Q)\le 2$.   Since  $Q$ has an even angle
  at  $q$,    we conclude that  $n(Q)=2$ 
and  $(v,  q, p')$ is the boundary of  a  chamber.

\end{proof}

\b{Prop}\label{sec8la4}
{Suppose $R_1,  R_2\not=(2,3,8)$ are right  triangles.
Let $A\subset\De_1$ be an  apartment and   $\xi_1\xi_2, \eta_1\eta_2\subset A$ two
 intersecting geodesics contained in the 1-skeleton of 
 $A$.    Then $\xi'_1\xi'_2\cap \eta'_1\eta'_2\not=\phi$. }


\end{Prop}

\b{proof}
  Denote  $v=\xi_1\xi_2\cap  \eta_1\eta_2$.
The proposition follows from Lemma \ref{sec8la2} if  $m(v)=2$.  From now on we  assume 
$m(v)=4, 6$ or $8$.
Suppose $\xi'_1\xi'_2\cap \eta'_1\eta'_2=\phi$. 
Proposition \ref{twogdisatd}
  implies $\xi'_1\xi'_2$ and   $\eta'_1\eta'_2$  have different types.
  By Lemma  \ref{sec8la3},  $\angle_v(\xi_1,  \eta_1)$  is an odd angle.
We may assume  $\eta'_1\eta'_2$ 
 contains vertices 
 $v'$ with $m(v')=2$, and  $\xi'_1\xi'_2$  does not.
Let $\xi_3\xi_4, \eta_3\eta_4 \subset A^{(1)}$ be  geodesics through $v$  such that 
$\angle_v(\xi_1, \xi_3)$ and $\angle_v(\eta_1,  \eta_3)$ are  nonzero  even angles. 
 Lemma \ref{sec8la3}  implies  $\xi'_1\xi'_2$ and $\xi'_3\xi'_4$   are of one  type,
    while    $\eta'_1\eta'_2$ and $\eta'_3\eta'_4$  are of    another  type.
Since we have assumed $\xi'_1\xi'_2\cap \eta'_1\eta'_2=\phi$, Lemma \ref{sec8la1} implies 
 $\xi'_3\xi'_4\cap \eta'_1\eta'_2=\phi$,  $\xi'_3\xi'_4\cap \eta'_3\eta'_4=\phi$
  and $\xi'_1\xi'_2\cap \eta'_3\eta'_4=\phi$.

\begin{figure}[h]
\centering\epsfig{file=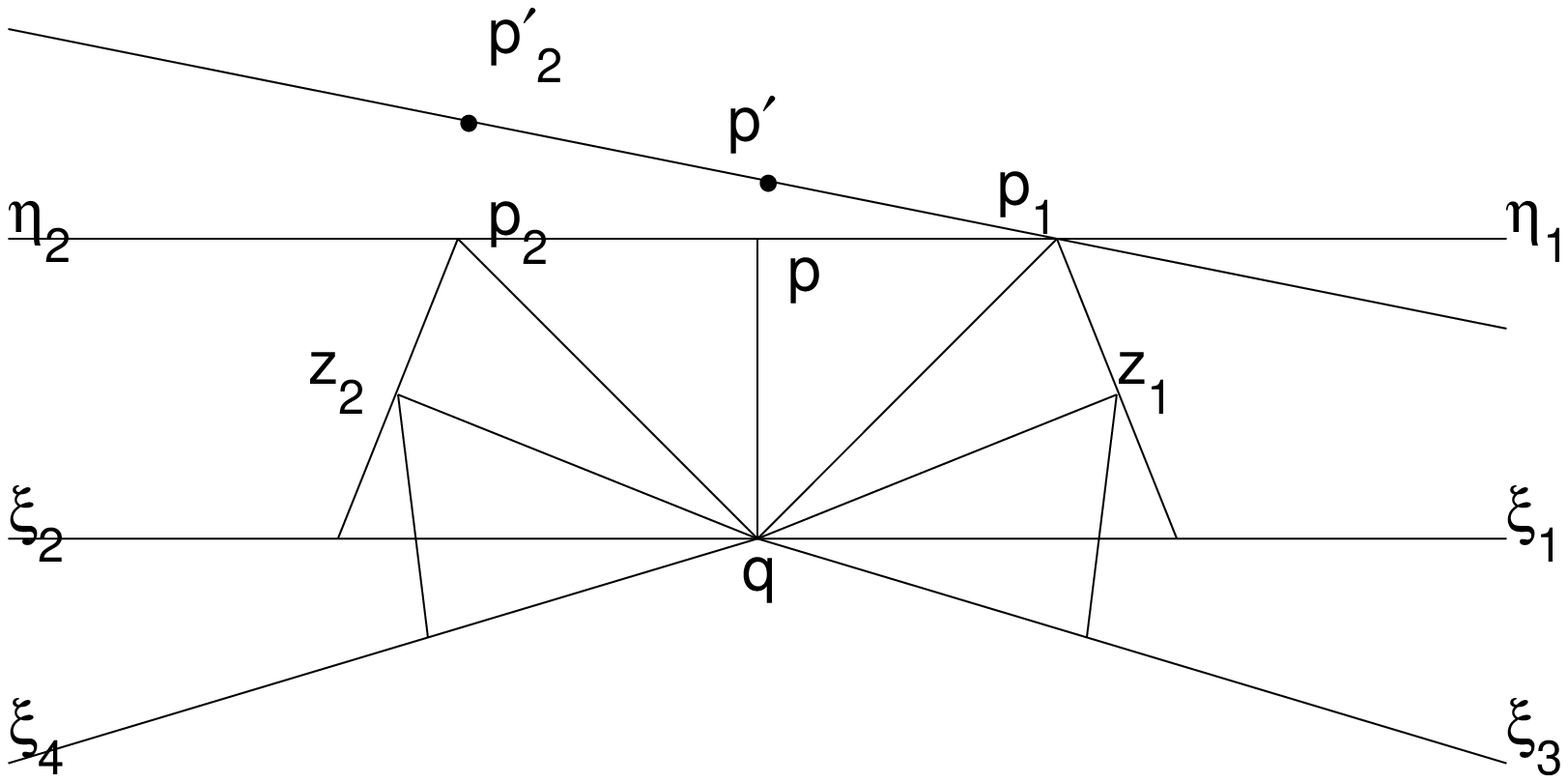,  height=4cm,  width=11cm}\caption{}\label{8p13}
\end{figure}

Let $pq$ ($p\in \eta'_1\eta'_2$, $q\in \xi'_1\xi'_2$) be the perpendicular between 
 $\eta'_1\eta'_2$  and  $ \xi'_1\xi'_2$. Note $m(p)=2,  m(q)=4$.   
Now Lemma \ref{sec8la3.9}    applied to $\eta'_1\eta'_2$ and $\xi'_1\xi'_2$, $\xi'_3\xi'_4$ 
implies $q\in \xi'_3\xi'_4$.
Let $p_i\in p\eta'_i$ ($i=1,2$) be the vertex adjacent to $p$. 
 Denote by $\sigma,   \tau_i$,    
  $\sigma_j\in  \Link(\De_2, q)$ 
     ($i=1,2$,   $j=1,2,3,4$)  the initial directions of $qp$, $qp_i$  and   $q\xi'_j$ respectively.   
 Proposition \ref{twogdisatd}  implies the     map  
 $f:\{\sigma_1, \sigma_2, \sigma_3, \sigma_4 \}\ra \{\sigma\tau_1, \sigma\tau_2\}$, 
  where $f(\sigma_j)$ is the initial edge of $\sigma\sigma_j\subset \Link(\De_2, q)$,
  is well-defined and satisfies the property: 
 $f(\sigma_1)\not=f(\sigma_2)$,  $f(\sigma_3)\not=f(\sigma_4)$.
  We may assume $f(\sigma_1)=f(\sigma_3)=\sigma\tau_1$  and 
$f(\sigma_2)=f(\sigma_4)=\sigma\tau_2$.
    Then  $\sigma\tau_1\subset \sigma\sigma_1\cap \sigma\sigma_3 $. 
The fact that $\xi'_1\xi'_2$ and $\xi'_3\xi'_4$ are at different sides implies that 
$0< \angle_q(\xi'_i, \xi'_j)<\pi$ for $1\le i\le 2$,   $3\le j\le 4$.
Since $m(q)=4$ and 
$\sigma\sigma_j$ has combinatorial length 3, we conclude that $\sigma\sigma_1\cap \sigma\sigma_3 $
  has combinatorial length 2. Similarly  $\sigma\sigma_2\cap \sigma\sigma_4 $  also  
  has combinatorial length 2.   Let $z_1, z_2$ be as shown in Figure \ref{8p13},  and $\omega_i\in \Link(\De_2, q)$ be the 
 initial direction of $qz_i$.  Then we have proved  $\sigma\sigma_1\cap \sigma\sigma_3=\sigma\omega_1 $,
 $\sigma\sigma_2\cap \sigma\sigma_4=\sigma\omega_2 $.  It follows that
  $\sigma_1\sigma_3=\sigma_1\omega_1\cup  \omega_1\sigma_3$  and 
 $\sigma_2\sigma_4=\sigma_2\omega_2\cup  \omega_2\sigma_4$.

Now Lemma \ref{sec8la3.92}
applied to $\xi'_1\xi'_2$   and 
$\eta'_1\eta'_2$, $\eta'_3\eta'_4$  implies that $p_i\in \eta'_3\eta'_4$ 
 for $i=1$ or $2$. We may assume $p_1\in \eta'_3\eta'_4$.  Let  $p'q'$ ($p'\in \eta'_3\eta'_4$, $q'\in \xi'_1\xi'_2$)
  be  the perpendicular between $\eta'_3\eta'_4$ and $\xi'_1\xi'_2$.   Then 
   $q'=q$ and  $T=(p_1,q, p')$ is the boundary of a chamber.  
 Since  $\eta'_1\eta'_2$, $\eta'_3\eta'_4$  are at different sides and 
$\angle_{p_1}(\eta'_1, z_1)=\pi$,  we have $p'\not=z_1$. 
Let $\sigma'\in \Link(\De_2, q)$ be the initial direction of $qp'$.    
Then $\sigma'\sigma_j=\sigma'\tau_1 \cup \tau_1\sigma_j$ for  $j=1,3$. Let $p'_2\in \eta'_3\eta'_4, p'_2\not=p_1$  be the 
  vertex that is adjacent to $p'$  but different from $p_1$,  and  
  $\tau'_2\in \Link(\De_2, q)$  the initial direction of $qp'_2$.  
Now the edge path
$\sigma_2\sigma\cup \sigma\tau_1\cup \tau_1\sigma'\cup \sigma'\sigma_2$ is a geodesic loop with length
  $2\pi$ in the   generalized polygon   $\Link(\De_2, q)$.  It follows that this loop is injective.  In particular, 
if $\omega'_2$   denotes   the only vertex 
 in the interior of $\sigma_2\tau'_2$,  then $\omega'_2\not=\omega_2$.  
Now the argument in the second paragraph applied to $\eta'_3\eta'_4$  and 
 $\xi'_1\xi'_2$,  $\xi'_3\xi'_4$    implies  that
  $\sigma_2\sigma_4=\sigma_2\omega'_2\cup \omega'_2\sigma_4$.    Consequently there are two 
 different   geodesic  segments with length $\pi/2$ from $\sigma_2$ to $\sigma_4$: 
   $\sigma_2\omega_2\cup  \omega_2\sigma_4$  and 
$\sigma_2\omega'_2\cup \omega'_2\sigma_4$. This contradicts to the fact that 
 $\Link(\De_2, q)$ is a $\CAT(1)$ space.

\end{proof}


For any $c\in {\mathcal{D}}_{A,v}$ we denote its image in ${\mathcal{D}}'_{A,v}$  by $c'$.

 \b{Le}\label{newl1}
{Proposition \ref{serk2po2}  holds for those $v$ with $m(v)\not= 4$.}

\end{Le}

\b{proof}
Note there are $m(v)$ geodesics in ${\mathcal{D}}'_{A,v}$.
  By Proposition \ref{sec8la4} any two geodesics in ${\mathcal{D}}'_{A,v}$ have nonempty 
 intersection,  in particular, Proposition \ref{serk2po2}  holds for those $v$ with $m(v)=2$. 
Suppose $m(v)=6$ or $8$.  Lemma \ref{sec8la3}
implies  that  if $\xi_1\xi_2, \eta_1\eta_2\in {\mathcal{D}}_{A,v}$ make  even angles, then 
  $\xi'_1\xi'_2$,   $\eta'_1\eta'_2$  also make  even angles. 
Fix three geodesics $c_1, c_2, c_3\in  {\mathcal{D}}_{A,v}$  that pairwise make  nonzero  even angles. 
 Since $R_2\not=(2,3,8)$,     Proposition \ref{tri2}  implies
  $c'_1, c'_2, c'_3$  intersect in a unique vertex $w\in \De_2$.
  Suppose there is a geodesic $c\in {\mathcal{D}}_{A,v} $ such  that  $w\notin c'$. 
Let $x_i=c'_i\cap c'$ ($i=1,2,3$). We may assume
 $x_2$ lies between $x_1$ and $x_3$.  
Consider the triangle $(w, x_1, x_3)$:  $x_2\in \interior(x_1x_3)$ and both
 $\angle_w(x_1, x_2)$ and $\angle_w(x_2, x_3) $  are even. 
Proposition \ref{tri2}
implies  $R_2=(2,3,8)$,  contradicting to our assumption.

\end{proof}

\begin{figure}[h]
\centering\epsfig{file=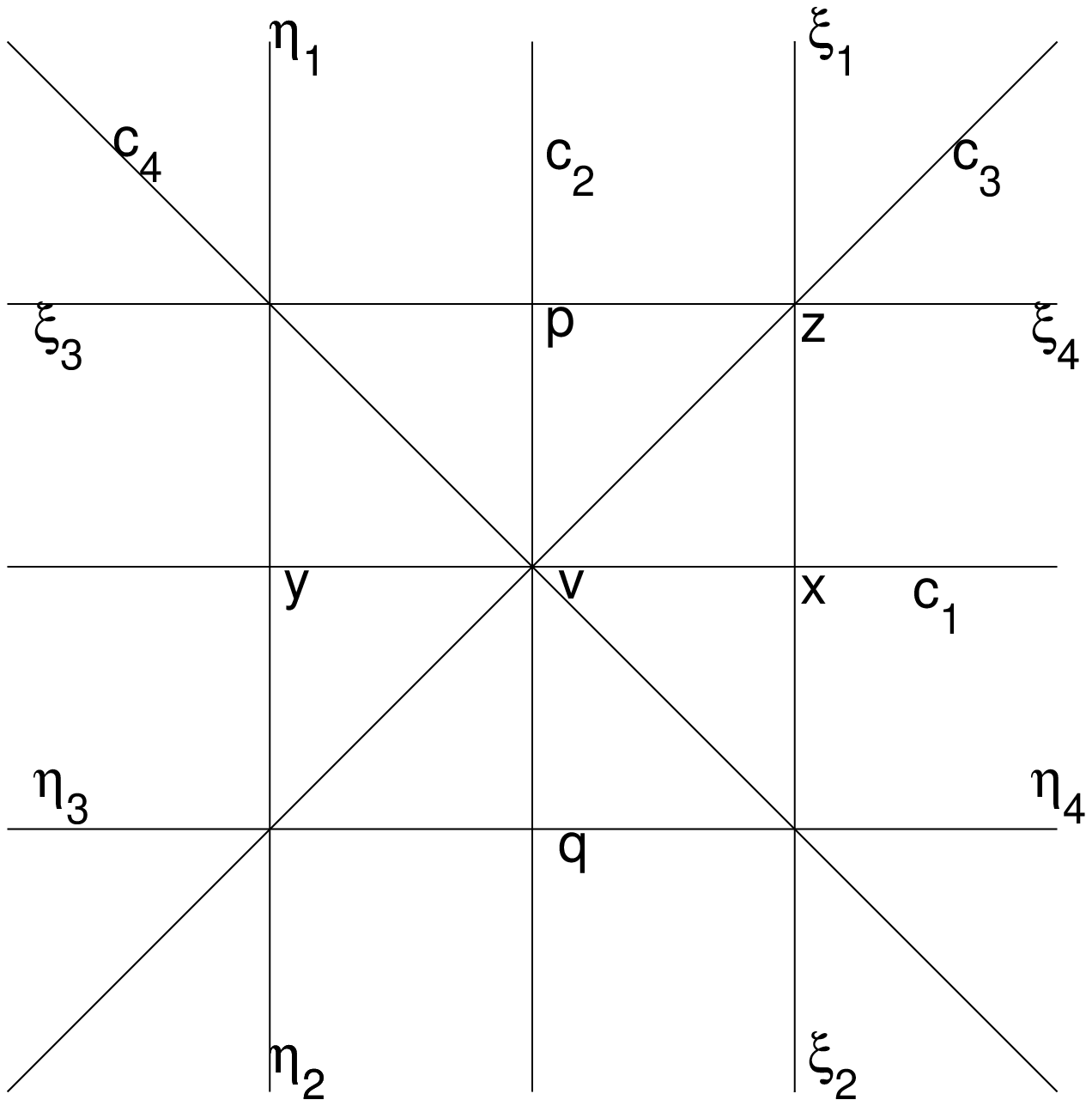,  height=6cm,  width=6cm}\caption{}\label{8l15}
\end{figure}

\b{Le}\label{newl2}
{Proposition \ref{serk2po2}  holds for those $v$ with $m(v)= 4$.}

\end{Le}

\b{proof}
Let $c_1, c_2\in {\mathcal{D}}_{A,v}$ be the two geodesics that contain vertices $v'$ with $m(v')=2$, 
 and $c_3$, $c_4$ the other two geodesics in   ${\mathcal{D}}_{A,v}$, see  Figure \ref{8l15}. Let $x, y\in c_1$ be the two  
   vertices on $c_1$ adjacent to 
  $v$,  and  $p, q\in c_2$  the two vertices on $c_2$ adjacent to $v$, as shown in Figure \ref{8l15}.   Let 
 $\xi_1\xi_2,\eta_1\eta_2\subset A^{(1)}$ be the two geodesics in $A$ perpendicular to $c_1$ and
 passing through $x$, $y$ respectively,  and  $\xi_3\xi_4, \eta_3\eta_4\subset A^{(1)}$  the two geodesics in $A$  
perpendicular to $c_2$ and
 passing through $p$ and $q$ respectively. 
Let $w=c'_1\cap c'_2$. Suppose $w\notin c'_3$.  Denote $w_i=c'_3\cap c'_i$ for $i=1,2,4$,
  and $w'_2=c'_2\cap c'_4$.  Since $c_1, c_2$ make  even angles, Lemma \ref{sec8la3}  implies that 
 $c'_1$, $c'_2$ make  even angles  at $w$.   Proposition \ref{tri2} applied to 
    $(w, w_1, w_2)$   shows  that 
$w_1w_2$ contains exactly one  vertex $w_0$  in the interior,  and $m(w_0)=2$ or $4$.   
Now 
  we  observe $w_4\in \{w_0,w_1, w_2\}$. Suppose   otherwise, say, $w_1\in \interior(w_2w_4)$ holds. 
 Then  $(w'_2, w_4, w_2)$ has 
  an even angle at $w_4$ and one side $w_2w_4$ contains two vertices $w_0$ and $w_1$
 in the interior, contradicting to Proposition \ref{tri2}.

We first consider the case $w_4\in \{w_1, w_2\}$, say, $w_4=w_1$. 
The triangle $(w'_2, w_1, w_2)$ has an even angle at $w_1$. Since the side 
 $w_1w_2$ contains $w_0$ in the interior,  $m(w_0)=2$ and $m(w_1)=6$ or $8$. 
Hence $\angle_{w_1}(w_2, w'_2)=\pi/3$ or $\pi/4$.  
The three geodesics $c_1, c_3$ and $\xi_1\xi_2$ have different types, 
   and any two of them have nonempty  intersection.   Corollary \ref{corse8c1}
implies $w_1\notin \xi'_1\xi'_2$.
 For the same reason, $w_1\notin \eta'_1\eta'_2$.  Let $z_i=c'_i\cap \xi'_1\xi'_2$ for $i=3,4$.  
 Then $(w_1, z_3, z_4)$ has an even angle at $w_1$. Since $m(w_1)=6$ or $8$,
  $\angle_{w_1}(z_3, z_4)=\pi/3$ or $\pi/4$.  It follows that either
 $\{w_2, w'_2\}=\{z_3, z_4\}$  or $\{w_2, w'_2\}\cap\{z_3, z_4\}=\phi$.
If $\{w_2, w'_2\}=\{z_3, z_4\}$,  then $\xi'_1$ is connected to an endpoint of $c'_2$ by a geodesic 
contained in  $\De_2^{(1)}$.
 It follows  that $\xi_1$ is connected to an endpoint of $c_2$ by 
a geodesic contained in  $\De_1^{(1)}$, which is not true.
 Hence  $\{w_2, w'_2\}\cap\{z_3, z_4\}=\phi$. Replacing $\xi'_1\xi'_2$ with $\eta'_1\eta'_2$ in the above argument,
  one sees $\eta'_1\eta'_2$  contains $z_3$ and $z_4$. Then $\eta'_1$ is connected to $\xi'_1$ or $\xi'_2$
 by  a geodesic in $\De_2^{(1)}$. It follows that $\eta_1$ is connected to $\xi_1$ or $\xi_2$
 by  a geodesic in $\De_1^{(1)}$, which is not true. 
  The contradiction shows  $w_4\not=w_1$.  Similarly   $w_4\not=w_2$.  Hence $w_4=w_0$. 
The triangle $(w_0, w_2, w'_2)$ has an even angle at $w_0$. It follows that $m(w_0)\not=2$ and hence $m(w_0)=4$. 
Proposition \ref{tri2} applied to $(w_0, w_2, w'_2)$ 
implies that $\angle_{w_2}(w_0, w'_2)=\pi/6$ or $\pi/8$.  Finally  we apply 
Proposition \ref{tri2} to   $(w, w_1, w_2)$  and conclude  that $w=w'_2$.

Suppose $w_0\notin \xi'_1\xi'_2$. 
 Then the three geodesics $c'_3$, $c'_4$
 and $\xi'_1\xi'_2$ form a triangle with  an even angle  at 
 $w_0$. 
It follows that either
 $w_1\in  \xi'_1\xi'_2$  or $w_2\in  \xi'_1\xi'_2$.
 On the other hand, the preceding paragraph shows 
  $w_1\notin \xi'_1\xi'_2$. Hence $w_2\in  \xi'_1\xi'_2$.
Let $w_5=\xi'_1\xi'_2\cap  c'_4$.   Note $w_5\not=w$, otherwise $\xi'_2$ is connected to one of the endpoints of $c'_2$ by a geodesic contained in $\De_2^{(1)}$, which is not true.
Let $z=c_3\cap \xi_1\xi_2$. Then $m(z)=6$ or $8$ and $z\in \xi_3\xi_4$.  Since $w_2=c'_3\cap \xi'_1\xi'_2$,  
Lemma \ref{newl1}
implies that $w_2\in \xi'_3\xi'_4$.  Let $w_6=c'_4\cap \xi'_3\xi'_4$. 
Notice $w_6\in ww_5$.   Suppose  otherwise, say $w_5\in \interior(ww_6)$  holds, then $(w_2, w_5, w_6)$ has an 
angle $\angle_{w_5}(w_2, w_6)=5\pi/6$ or $7\pi/8$, which is impossible.
There are three vertices on $ww_5$: $w$, $w_0$, $w_5$.  
Assume $w_6=w$. Then 
 $\xi'_3\xi'_4\cap c'_2$   is   a nontrivial interval, which implies $\xi'_3$ 
 is connected to an  endpoint of $c'_2$  by a geodesic  contained in $\De_2^{(1)}$.
It follows that $\xi_3$ is connected to an  endpoint of $c_2$  by a geodesic  contained in $\De_1^{(1)}$,
  which is not true. Similarly $w_6\not=w_0, w_5$.   The contradiction shows that 
$w_0\in \xi'_1\xi'_2$.  Similarly $w_0\in \eta'_1\eta'_2$.

 Since $c_3$ and $\xi_1\xi_2$ have different types,
 Lemma \ref{sec8la3}  implies $c'_3$ and $\xi'_1\xi'_2$ have different types. 
 Similarly  $c'_3$ and   $\eta'_1\eta'_2$   have different types.
  So $\xi'_1\xi'_2$  and $\eta'_1\eta'_2$  have the same type.  $\xi'_1\xi'_2$  and $\eta'_1\eta'_2$
  make 4 angles at $w_0$:  $\angle_{w_0}(\xi'_1, \eta'_1)$,  
 $\angle_{w_0}(\xi'_1, \eta'_2)$,  $\angle_{w_0}(\xi'_2, \eta'_1)$,  $\angle_{w_0}(\xi'_2, \eta'_2)$. 
 If any one of these angles is $0$ or $\pi$, then some $\xi'_i$ is connected to
some $\eta'_j$ by a geodesic contained in $\De_2^{(1)}$.   It follows that 
$\xi_i$ is connected to
some $\eta_j$ by a geodesic contained in $\De_1^{(1)}$,  which is not true.
Recall $m(w_0)=4$. It follows that 
 the 4 angles that $\xi'_1\xi'_2$  and $\eta'_1\eta'_2$
 make at $w_0$ are all $\pi/2$.  Hence   $\xi'_1\xi'_2$  and $\eta'_1\eta'_2$
are locally contained in an apartment and 
Lemma \ref{ritrdirigh}
  implies $\xi'_1\xi'_2$  and $\eta'_1\eta'_2$
are at different sides. Consequently, $\xi_1\xi_2$  and $\eta_1\eta_2$
are at different sides, which is not true.  The contradiction shows $w\in c'_3$. Similarly 
 $w\in c'_4$.

\end{proof}

\subsubsection{Geodesics through a fixed vertex }\label{containsno}

\b{Prop}\label{sec8.3p1}
{Let $\De_1$,  $\De_2$ be two Fuchsian buildings, and $h:\p\De_1\ra \p\De_2$  a homeomorphism that preserves the combinatorial cross ratio almost everywhere.
  Let  $v\in \De_1$ be  a vertex.
 If both $R_1$ and $R_2$ are right    triangles different from $(2,3,8)$, then 
the geodesics in ${\mathcal{D}}'_{v}$ intersect in a unique vertex of $ \De_2$.}

\end{Prop}

\b{proof}
Proposition \ref{serk2po2}  implies that 
  for  any vertex  $v\in \De_1$  and any apartment $A$ containing $v$, the geodesics in 
 ${\mathcal{D}}'_{A,v}$    intersect in  a  unique vertex  $v_A\in \De_2$. 
By Lemma \ref{redaprtowh}, we only need to show that $v_{A_1}=v_{A_2}$  holds for any two apartments 
 $A_1, A_2\subset \De_1$    containing  a vertex $v$   with  
 $\Link(A_1, v)\cap \Link(A_{2}, v)$     a  half apartment in $\Link(\De_1, v)$.
Let $B_i=\Link(A_i,v)$($i=1,2$), and $\omega_1$, $\omega_2$ the two endpoints of $B_1\cap B_2$.
Denote by $m$ the midpoint of $B_1\cap B_2$,  and  $m_i\in B_i$ the point in $B_i$ opposite to $m$.
Since $m(v)$ is even, $m, m_1$ and $m_2$ are all  vertices in $\Link(\De_1, v)$.    
Let $\xi_1, \xi_3\in \p A_1$ and $\xi_2, \xi_4\in \p A_2$ such that
 $v\xi_1$, $v\xi_2$ have initial direction $\omega_1$ and 
  $v\xi_3$, $v\xi_4$ have initial direction $\omega_2$. Similarly
let $\eta_1,\eta_3\in \p A_1$ and $\eta_2, \eta_4\in \p A_2$ such that 
$v\eta_3$, $v\eta_4$ have initial direction $m$,   and 
  $v\eta_i$ ($i=1,2$) has initial direction $m_i$.
Then we have $\eta_i\eta_j=v\eta_i\cup v\eta_j\subset \De_1^{(1)}$  if  $i\not=j$,   $\{i,j\}\not=\{3,4\}$.

Now consider the three geodesics $\eta'_1\eta'_2$,  $\eta'_1\eta'_3$,
$\eta'_3\eta'_2\subset \De_2^{(1)}$.  Since $R_2\not=(2,3,8)$, Lemma \ref{threeev} implies 
there is some vertex $w\in \De_2$ such that $\eta'_1\eta'_2=w\eta'_1\cup w\eta'_2$,
$\eta'_1\eta'_3=w\eta'_1\cup w\eta'_3$  and 
$\eta'_3\eta'_2=w\eta'_3\cup w\eta'_2$.   Notice that  $\eta_i\eta_j$ ($1\le i\not=j\le 3$)  and $\xi_1\xi_3$ 
 make a right angle and are locally contained in an
 apartment.  Lemma \ref{ritrdirigh}  implies $\eta_i\eta_j$ and $\xi_1\xi_3$
  are at different sides. It follows that $\eta'_i\eta'_j$ and $\xi'_1\xi'_3$
  are also  at different sides.   Now Lemma \ref{sec8la2.5}
implies  $w\in \xi'_1\xi'_3$.    
 Similarly $w\in \xi'_2\xi'_4$. 
Note $v_{A_1}=\xi'_1\xi'_3\cap \eta'_1\eta'_3=w$.

Now apply the above argument to $\eta_1, \eta_2, \eta_4$, instead of $\eta_1, \eta_2, \eta_3$,
 we see there is a   vertex $w'\in \De_2$ such that 
 $\eta'_1\eta'_2=w'\eta'_1\cup w'\eta'_2$,
$\eta'_1\eta'_4=w'\eta'_1\cup w'\eta'_4$  and 
$\eta'_4\eta'_2=w'\eta'_4\cup w'\eta'_2$,  and both $\xi'_1\xi'_3$  and  $\xi'_2\xi'_4$
 contain $w'$.   In particular, $w'=\xi'_1\xi'_3\cap \eta'_1\eta'_2=w$.  
    Therefore $v_{A_2}=\xi'_2\xi'_4\cap \eta'_2\eta'_4=w'=w=v_{A_1}$.

\end{proof}

\b{Prop}\label{onerionex}
{Let $\De_1$,  $\De_2$ be two Fuchsian buildings, and $h:\p\De_1\ra \p\De_2$  a homeomorphism that preserves the combinatorial cross ratio almost everywhere. 
If one of $R_1$,   $R_2$ is $(2,3,8)$, then so is the other.}
\end{Prop}

\b{proof} The proof is similar to that of Proposition \ref{onetoneq}.
Suppose $R_1=(2,3,8)$  and $R_2\not=(2,3,8)$.   Proposition \ref{onetoneq}
implies that $R_2$ is also a right  triangle. 
 Let $A$ be an apartment of $\De_1$. Let $\mathcal{G}$ be the  set of  geodesics contained 
 in $A^{(1)}$ that do not contain any vertex $v$ with $m(v)=3$.
  Let $c_1,c_2\in \mathcal{G}$
 with $c_1\cap c_2\not=\phi$.  Then $c_1$ and $c_2$ are at different sides and make  even angles  at 
 $c_1\cap c_2$. Lemma \ref{sec8la3}       implies 
  that their images $c'_1$, $c'_2$ in $\De_2$ have the same type. 
Note $c'_1$, $c'_2$  are  also  at different sides.
Since $R_2\not=(2,3,8)$,  Proposition \ref{twogdisatd}
implies  $c'_1\cap c'_2\not=\phi$. It follows that $c'_1$, $c'_2$
make  even angles   at $c'_1\cap c'_2$.

The geodesics in $\mathcal{G}$ divide $A$ into triangles, which shall be called \lq\lq chambers".
  Let $D$ be such a \lq\lq chamber", and $c_1$, $c_2$, $c_3\in{\mathcal{G}}$ the three geodesic containing the three
 sides of $D$.  Since $R_2\not=(2,3,8)$,  Proposition  \ref{tri2}  
 implies the images $c'_i\subset \De_2$ ($i=1,2,3$) have a common
vertex, which we shall denote by $w_{D}$.  Similarly, 
for  a vertex $v\in A$ with $m(v)=8$   and  $c_i$ ($i=1,2,3,4$) the 4  geodesics in $\mathcal{G}$ through  $v$, 
  Proposition \ref{tri2} 
 implies the images $c'_i\subset \De_2$ ($i=1,2,3,4$) have a common
vertex, which is denoted by $w_v$.

Now let $D_1$, $D_2$ be two adjacent \lq\lq chambers" in $A$. Then there is a vertex
 $v\in A$ and 5 geodesics $c_i$, $1\le i\le 5$ such that 
 $c_1$, $c_2$, $c_3$ pass through $v$, $c_1$, $c_2$, $c_4$ contain the sides of $D_1$,
  and $c_2$, $c_3$, $c_5$ contain the sides of $D_2$.     The above paragraph shows
 $w_D=w_v=w_{D'}$. Notice  $A$ is  gallery connected: given any two 
\lq\lq chamber"s $D$, $D'$, there is a finite sequence 
$D_0=D,  D_1, \cdots,  D_k=D'$ such   that   $D_i$ and $D_{i+1}$ share a side. 
It follows that 
 there is a vertex $w\in \De_2$ such that the images of the geodesics in $\mathcal{G}$
 all contain $w$.  We get a contradiction as in the proof of Proposition \ref{onetoneq}.

\end{proof}

\subsection{The exceptional case $(2,3,8)$} \label{secsingular}

In this section we prove Theorem \ref{main} when both $R_1$  and $R_2$ are $(2,3,8)$.  
As usual  Proposition \ref{tri2}  is  used in the manner as
indicated in Remark \ref{prop5.2e}.

\subsubsection{Buildings with chamber $(2,3,8)$} \label{subsecnos}

Let $\De$ be a Fuchsian building with chamber $(2,3,8)$.  
 The vertex links of $\De$ are generalized polygons.  
By Section \ref{2hyb} there are two integers $p,q\ge 2$,
  $p\not=q$ such that the following holds:  an edge $e$ is contained in exactly $p+1$ chambers if $e$ is incident 
 to a vertex $v$ with $m(v)=3$, and is contained in exactly $q+1$ chambers  otherwise.
  There are only two types of complete geodesics in $\De$:\newline
(1) Type I:  geodesics that  do not contain  any vertex   $v$ with  $m(v)=3$;\newline
(2)   Type II: geodesics that   contain   vertices    $v$ with  $m(v)=3$.

We shall use the following terminologies.  We say an edge is \e{numbered} by $i$ ($i=p$ or $q$) if it is contained in exactly
 $i+1$ chambers.  All edges in 
 Type I geodesics are numbered by $q$, and  all edges in Type II geodesics are numbered by $p$. 
We say a  vertex $v\in \De$ has  \e{index}  $i$ ($i=2,3$ or $8$) if $m(v)=i$;  
  we   also  say $v$ is \e{indexed } by $i$.
 If we record the indexes of vertices on a Type II geodesic  in linear order, then they are:  $\cdots, 8,3,2,3, \cdots$, with period 4. 
 We say the vertices on a Type II  geodesic are \e{periodically indexed}  by $8,3,2,3$.
Similarly  if  we record the indexes of vertices on a Type I geodesic in linear order, 
    then they are:  $\cdots, 8,2 \cdots$, with period 2. 
  We say  the vertices on a Type I geodesic are periodically indexed by $8,2$.  

The following lemma follows  directly from Proposition \ref{tri2}.

\b{Le}\label{sec9.1l1}
{Let $T\subset \De^{(1)}$ be a triangle  homeomorphic to  a  circle 
such that all the sides of  $T$  are contained  in  Type II geodesics.
Then $S(T)$ is the union of two chambers  as shown in 
  Figure  \ref{238}\e{(f)}.}

\end{Le}


Recall the chamber $(2,3,8)$ has area $A_0=\pi/24$.

 \b{Le}\label{sec9.1l2}
{There is no quadrilateral  $Q\subset \De^{(1)}$ with the following properties: 
all sides   of  $Q$ are contained in Type II geodesics, and there is some side $xy$ of $Q$ such that 
 the angles at $x$ and $y$ are $2\pi/3$.}

\end{Le}

\b{proof}
Suppose there is such a quadrilateral. Notice each angle of $Q$ is $\ge \pi/4$. 
  Proposition  \ref{quar11}
  implies   that  $n(Q)\le 4$.  But there are at least two chambers  in $S(Q)$ 
  incident  to each of
 $x$ and $y$. It follows that $S(Q)$ is the union of these 4 chambers. 
But the boundary of  this  union is not a quadrilateral, contradiction.  

\end{proof}

 \b{Le}\label{sec9.1l3}
{Let $Q=(x,y,z,w)\subset \De^{(1)}$ be a  quadrilateral  homeomorphic to  a  circle 
  with the following properties: 
all sides of  $Q$ are contained in Type II geodesics,    
 the angles at $x$ and $y$ are even and at least 
  two vertices  in the interior of $xy$  are not indexed by $2$.
  Then  $m(v)\not=8$ for every   vertex $v\in \interior(xy)$  and there is a vertex 
 $v'\in zw$  with $m(v')=8$.}

\end{Le}

\b{proof}
If $m(x)=2$, then $Q$ is actually a triangle and we obtain a contradiction to Lemma \ref{sec9.1l1}.
 Hence $m(x)\not=2$, and similarly $m(y)\not=2$.
At least one of $m(x), m(y)$ is $\not=3$, otherwise we have a contradiction to Lemma \ref{sec9.1l2}.

First assume exactly one of $m(x), m(y)$ is $3$, say, $m(x)=3$.
Since the angle at $x$ is even and each angle is at least $\pi/4$, Proposition  \ref{quar11}
  implies    $n(Q)\le 14$. 
Recall the vertices in a  Type II geodesic are periodically indexed by $8,3,2,3$.
 Since $m(x)=3$, $m(y)=8$ and at least two vertices in   $\interior(xy)$ 
  are not indexed by $2$,  there are 
    vertices  $v_1, v_2,  v_3, v_4 \in \interior(xy)$ 
 such that $x, v_1, v_2, v_3, v_4, y$ are in linear order and 
$m(v_1)=8, m(v_2)=3, m(v_3)=2$,  $m(v_4)=3$. Recall the angles at $x$ and $y$ are even. 
We count the chambers in $S(Q)$: 
there are at least 2 chambers incident to $x$,  7 incident to $v_1$ but not $x$, 
 2 incident to $v_2$ but not $v_1$,
 1  incident to $v_3$ but not $v_2$, 2 incident to $v_4$ but not $v_3$, 
and  1 incident to $y$ but not $v_4$, for a total of 15. This contradicts to 
the above conclusion that $n(Q)\le 14$.

Now assume  $m(x)=m(y)=8$. Suppose there is  a  vertex $v\in \interior(xy)$  with $m(v)=8$. 
 Let $x_1\in xy $ be the vertex on $xy$ closest to $x$ with  $m(x_1)=2$,
 and $y_1\in xy$ the vertex on $xy$ closest to $y$ with  $m(y_1)=2$. Let 
$C_1, C_2$ be chambers in $S(Q)$ incident to $x_1$, $y_1$ respectively, and $v_3$, $v_4$ their
 vertices  with $m(v_3)=m(v_4)=8$ respectively. 
Note none of $v_3, v_4$  lies on the sides of $Q$ incident to $x$ or  $y$
 because $\angle_{x}(v_3, y), \angle_{y}(v_4, x)=\pi/8$ and the angles of $Q$ at $x$ and $y$ are even. 
It follows that  $\{v_3, v_4\}\cap \{z,w\}=\phi$.  
Hence there are at least 8 chambers of $S(Q)$ incident to each of $v_3$, $v_4$, $v$. There are also at least two  
chambers incident to each of $x$ and $y$. Therefore there are at least 28 chambers in $S(Q)$.
 On the other hand, since each angle  of $Q$ is at least $\pi/4$, Proposition \ref{quar11}  implies
$n(Q)\le 24$, a contradiction. 
Hence $m(v)\not=8$ for all vertices   $v\in \interior(xy)$.

Since 
 $m(v)\not=8$ for all vertices   $v\in \interior(xy)$,  
 there is exactly one vertex $v_0\in xy$ with $m(v_0)=2$. 
Let $C$ be a chamber in $S(Q)$ incident to $v_0$, and $v'$ the vertex of $C$ with $m(v')=8$. 
If $v'\in Q$, then the above argument shows $v'\in \interior(zw)$ and the lemma holds. 
So we assume $v'$ lies in the interior of $S(Q)$.
Then there are at least 16 chambers in $S(Q)$ incident to $v'$. There are also (at least) two chambers 
 incident to each of $x$, $y$. Hence there are at least 20 chambers in $S(Q)$.
Suppose  there is no vertex $v'\in zw$ with $m(v')=8$.  Then the angles at $z$ and $w$ are $\ge \pi/3$. 
Since the angles at $x$ and $y$ are $\ge \pi/4$, Proposition \ref{quar11}  implies $n(Q)\le 20$. 
 So $S(Q)$ must be the union of the 20 chambers exhibited above. However, the boundary of this union is not 
 a quadrilateral, contradiction.

\end{proof}





\b{Le}\label{sec9.1l12}
{Let $\xi_1\xi_2, \eta_1\eta_2\subset \De^{(1)}$ \e{($\xi_1, \xi_2, \eta_1, \eta_2\in \p\De$)}
be two  geodesics that are locally contained in an apartment. If 
  $\xi_1\xi_2$, $\eta_1\eta_2$   have the same type  and make a right angle, then they are at different sides.}

\end{Le}

\b{proof} Let $x_0=\xi_1\xi_2\cap \eta_1\eta_2$.  The assumption implies $m(x_0)=8$.
If there is some $x\in \eta_1\eta_2$, $x\not=x_0$
such that $x\xi_1\cap \xi_1\xi_2=x'\xi_1$ ($x'\in\xi_1\xi_2$) is a ray and $x\xi_1\cap \eta_1\eta_2=\{x\}$,
   then $(x', x_0, x)\subset \De^{(1)}$ has an even angle at $x'$ and a right angle at $x_0$. Proposition \ref{tri2} implies 
  that $m(x_0)=2$, contradicting to the above observation. Hence  
   $x\xi_1\cap \xi_1\xi_2=\phi$  for any $x\in \eta_1\eta_2$, $x\not=x_0$.
Similarly  $x\xi_2\cap \xi_1\xi_2=\phi$  for any $x\in \eta_1\eta_2$, $x\not=x_0$
and  
$y\eta_i\cap \eta_1\eta_2=\phi$ ($i=1,2$) for any $y\in \xi_1\xi_2$, $y\not=x_0$.  
Then one argues as in  Lemma \ref{noritri}
 that $\xi_1\xi_2$, $\eta_1\eta_2$ are at different sides.

\end{proof}

\subsubsection{Two geodesics at different sides} \label{subsecspecial1}

Let $\De$ be a Fuchsian building with chamber $(2,3,8)$.  There are two integers
$p,q\ge 2$,
  $p\not=q$ 
such that all the edges in Type I geodesics  are  numbered by $q$,  
 and all the edges in Type II geodesics are  numbered by $p$.
The goal of this section is to show the following:

\b{Prop}\label{sec9p1}
{Let $\De$ be a Fuchsian building with chamber $(2,3,8)$.
If two disjoint geodesics $c_1, c_2\subset \De^{(1)}$ are at different sides, 
  then they have different types.}

\end{Prop}

The proof of Proposition   \ref{sec9p1}
is divided into two propositions (Propositions \ref{sec9p2} and \ref{sec9p3}).

\b{Prop}\label{sec9p2}
{Let $\De$ be a Fuchsian building with chamber $(2,3,8)$.
If two Type I  geodesics $\xi_1\xi_2, \eta_1\eta_2\subset \De^{(1)}$ are at different sides, 
  then  $\xi_1\xi_2  \cap  \eta_1\eta_2\not=\phi$.}

\end{Prop}

Recall the vertices on a Type I geodesic are periodically indexed by $8,2$. 

\b{Le}\label{sec9.3l1}
{Let $\xi_1\xi_2$, $\eta_1\eta_2$ be two disjoint Type I geodesics.
 If they are at different sides, then there are no vertices $y'\in \xi_1\xi_2$,
   $y\in \eta_1\eta_2$ such that $\angle_{y'}(y, \xi_1)=\angle_{y'}(y, \xi_2)=\pi$.}

\end{Le}

\b{proof}
Suppose there are vertices $y'\in \xi_1\xi_2$,
   $y\in \eta_1\eta_2$ such that $\angle_{y'}(y, \xi_1)=\angle_{y'}(y, \xi_2)=\pi$.
The assumption and Lemma \ref{onedif}  implies  that there exists 
 $x_i\in \xi_1\xi_2$ ($i=1,2$)  such that $x_i\eta_i\cap \eta_1\eta_2=x'_i\eta_i$ ($x'_i\in \eta_1\eta_2$) 
  is a  ray.  We may assume $x_i\eta_i\cap \xi_1\xi_2=\{x_i\}$.  
We claim $\angle_{y}(y', \eta_i)=\pi/4$ for $i=1,2$. The 
triangle inequality  then implies $\angle_y(\eta_1, \eta_2)\le \pi/2$,
 contradicting to the fact that $y\in \eta_1\eta_2$. Next we prove the claim.

The fact that $\xi_1\xi_2$, $\eta_1\eta_2$ 
are at different sides  implies  that  $\angle_y(y', \eta_i)<\pi$ for $i=1,2$. 
It follows that $x'_i\not=y$.  
The uniqueness of geodesics implies that $x'_i\in \interior(y\eta_i)$. 
Notice $x_iy=x_iy'\cup y'y$ and $(x_i, y, x'_i)$ is a triangle. 
Since $\xi_1\xi_2$, $\eta_1\eta_2$  are  Type I geodesics, all  the angles of 
$(x_i, y, x'_i)$  are  even.  Proposition \ref{tri2} then   implies  
$\angle_{y}(y', \eta_i)=\angle_y(y', x'_i)=\pi/4$.

\end{proof}

\b{Le}\label{sec9.3l2}
{Let $\xi_1\xi_2$, $\eta_1\eta_2$ be two disjoint Type I geodesics.
 Suppose $\xi_1\xi_2$, $\eta_1\eta_2$  are at different sides. 
Let $y_i\in \eta_1\eta_2$ \e{($i=1,2$)} be  a vertex such that $y_i\xi_i\cap \xi_1\xi_2=y'_i\xi_i$ is a ray
  and  $y_i\xi_i\cap \eta_1\eta_2=\{y_i\}$.
Then $y_1\not=y_2$.}

\end{Le}

\b{proof}
We suppose $y_1=y_2$.  Lemma \ref{sec9.3l1}
implies $y'_1\not=y'_2$.   The uniqueness of geodesic implies $y'_1\in \interior(y'_2\xi_1)$. 
The triangle $(y'_1, y'_2, y_1)$ has three even angles and 
  Proposition \ref{tri2}  implies $\angle_{y_1}(y'_1,y'_2)=\pi/4$ and $y_1y'_j$ is the union of two edges. 
Let $x_i\in \xi_1\xi_2$ be such that $x_i\eta_i\cap \eta_1\eta_2=x'_i\eta_i$ ($i=1,2$) is a ray
 and  $x_i\eta_i\cap \xi_1\xi_2=\{x_i\}$.

First suppose $x_i\notin y'_1y'_2$ for  $i=1$ or $2$. Then $x_i\in \interior(y'_j\xi_j)$ for some $j$. 
   The  uniqueness of geodesic implies $x'_i\in \interior(y_1\eta_i)$. 
Then $(y_1, x_i, x'_i)$ 
 has three even angles and has a   side 
$y_1x_i=y_1y'_j\cup y'_jx_i$ consisting of  at least three edges,  contradicting to Proposition \ref{tri2}. 
Hence $x_i\in y'_1y'_2$ for $i=1, 2$. Note $x_i\eta_i$ is  of Type I since it 
  shares a  ray with a Type I geodesic.  
 It follows that    $x_ix'_i$ makes even angles with $\xi_1\xi_2$. 
  Hence  
$m(x_i)\not=2$, otherwise $\angle_{x_i}(\xi_1, x'_i)=\angle_{x_i}(\xi_2, x'_i)=\pi$, contradicting to 
Lemma \ref{sec9.3l1}.  On the other hand, $y'_1y'_2$ contains only one vertex $v$ in the interior and 
it satisfies $m(v)=2$. Hence $\{x_1, x_2\}\subset \{y'_1,  y'_2\}$. In particular, there is some $j$ with 
 $x_1=y'_j$.
The uniqueness of geodesic implies 
$x'_1\in y_1\eta_1$.  The equality $x'_1=y_1$ does not hold since it would imply 
$\xi_j$ and  $\eta_1$ are connected by a geodesic contained in the 1-skeleton, contradicting to the assumption that 
    $\xi_1\xi_2$, $\eta_1\eta_2$  are at different sides.
Now $(x_1, y_1, x'_1)$ has three even angles and so $y_1x'_1$ is the union of two edges
 and $\angle_{y_1}(y'_j, \eta_1)=\pi/4$.  Similarly 
$x_2=y'_l$ for some $l=1,2$, $x'_2\in \interior(y_1\eta_2)$ and $\angle_{y_1}(y'_l, \eta_2)=\pi/4$.
 If $l=j$, then triangle inequality implies $\angle_{y_1}(\eta_1,\eta_2)\le \pi/2$,
contradiction to the fact that $y_1\in \eta_1\eta_2$.
If $l\not=j$, then 
$\angle_{y_1}(\eta_1, \eta_2)\le \angle_{y_1}(y'_j, \eta_1)+\angle_{y_1}(y'_1,y'_2)+\angle_{y_1}(y'_l, \eta_2)\le 3\pi/4$,
again a contradiction.

\end{proof}

The index $i+1$  is taken  mod 2  in the proofs of  Lemma \ref{sec9.3l3},   
 Propositions \ref{sec9p2} and \ref{sec9p3}.

\b{Le}\label{sec9.3l3}
{Let $\xi_1\xi_2$, $\eta_1\eta_2$ be two disjoint Type I geodesics.
 Suppose $\xi_1\xi_2$, $\eta_1\eta_2$  are at different sides. 
Let $y_i\in \eta_1\eta_2$ \e{($i=1,2$)} be such that $y_i\xi_i\cap \xi_1\xi_2=y'_i\xi_i$ is a ray.
Then $y'_1\not=y'_2$.}

\end{Le}

\b{proof}
We suppose  $y'_1=y'_2$.  We may assume $y_i\xi_i\cap \eta_1\eta_2=\{y_i\}$.   
The  triangle  $(y'_1, y_1, y_2)$ has three even angles,
 and hence each side is the union of two edges  and $\angle_{y_i}(y'_1, y_{i+1})=\pi/4$. 
Let $x_1\in \xi_1\xi_2$ be such 
that $x_1\eta_1\cap \eta_1\eta_2=x'_1\eta_1$ is  a ray   and  $x_1\xi_1\cap \xi_1\xi_2=\{x_1\}$. 
First assume $x_1\in \interior(y'_1\xi_1)$. 
The uniqueness of geodesic shows 
$x'_1\in \interior(y_1\eta_1)$. 
 Then $(x'_1, y_1, x_1)$ has three even angles  and its side  $x_1y_1=x_1y'_1\cup y'_1y_1$
 contains at least three edges, contradicting to Proposition \ref{tri2}. 
   Similarly $x_1\in \interior(y'_1\xi_2)$  cannot hold, and  we must have $x_1=y'_1$. 
Uniqueness of geodesic implies $x'_1\notin \interior(y_i\eta_2)$  for $i=1,2$. 
The fact that $\xi_1\xi_2$, $\eta_1\eta_2$  are at different sides  implies $x'_1\not=y'_1, y'_2$.
There is some $i=1,2$ such that $x'_1\in \interior(y_i\eta_1)$ and $y_{i+1}\notin y_i\eta_1$.
 Then $(x_1, x'_1, y_i)$ has three even angles and hence we have $\angle_{y_i}(x'_1, x_1)=\pi/4$.
 It follows that $\angle_{y_i}(\eta_1, \eta_2)\le \pi/2$, contradicting to the fact that
 $y_i\in \eta_1\eta_2$.

\end{proof}

Lemmas   \ref{sec9.3l2}, \ref{sec9.3l3} and uniqueness of geodesic  imply the following:

\b{Cor}\label{sec9.3c1}
{Let $\xi_1\xi_2$, $\eta_1\eta_2$ be two disjoint Type I geodesics.
 Suppose $\xi_1\xi_2$, $\eta_1\eta_2$  are at different sides. 
Let $y_i\in \eta_1\eta_2$ \e{($i=1,2$)} be such that $y_i\xi_i\cap \xi_1\xi_2=y'_i\xi_i$ is a ray
  and  $y_i\xi_i\cap  \eta_1\eta_2=\{y_i\}$.
Then $y_1\not=y_2$  and  $y'_1\in \interior(y'_2\xi_1)$.}

\end{Cor}

\b{Le}\label{sec9.3l4}
{Let $\xi_1\xi_2$, $\eta_1\eta_2$ be two disjoint Type I geodesics.
 Suppose $\xi_1\xi_2$, $\eta_1\eta_2$  are at different sides.
Let $y_1\in \eta_1\eta_2$,  $x_1\in \xi_1\xi_2$ be such that 
 $y_1\xi_1\cap \xi_1\xi_2=y'_1\xi_1$,  $x_1\eta_1\cap \eta_1\eta_2=x'_1\eta_1$ are rays
  and  $y_1\xi_1\cap  \eta_1\eta_2=\{y_1\}$,  $x_1\eta_1\cap \xi_1\xi_2=\{x_1\}$. 
Then $y'_1\in \interior(x_1\xi_1)$ and $x'_1\in \interior(y_1\eta_1)$.}

\end{Le}

\b{proof}
We show  (1) and (2) of Proposition \ref{tgeoddif}
cannot occur. 
First suppose (1) of Proposition \ref{tgeoddif}
occurs, that is, $y_1'\in\interior(x_1\xi_1)$  and $y_1\in \interior(x_1'\eta_1)$. In this case,
 $x_1y_1=x_1x'_1\cup x'_1y_1$ and $(y'_1, y_1, x_1)$ has three even angles. It follows that 
 the only vertex $x'_1$ in the interior of $x_1y_1$ satisfies $m(x'_1)=2$.  Since $\angle_{x'_1}(\eta_1, x_1)=\pi$,
  we also have $\angle_{x'_1}(\eta_2, x_1)=\pi$,  contradicting to Lemma \ref{sec9.3l1}.

Now suppose (2) of Proposition \ref{tgeoddif}
occurs, that is, $y_1'\in\interior(x_1\xi_1)$  and  $y_1=x_1'$.
Let $x_2\in \xi_1\xi_2$,
 $y_2\in \eta_1\eta_2$ 
 be such that $x_2\eta_2\cap \eta_1\eta_2=x'_2\eta_2$,  $y_2\xi_2\cap \xi_1\xi_2=y'_2\xi_2$   are  rays
  and  $x_2\eta_2\cap  \xi_1\xi_2=\{x_2\}$,  $ y_2\xi_2\cap  \eta_1\eta_2=\{y_2\}$. 
 Corollary \ref{sec9.3c1}  implies   $x'_2\in \interior(x'_1\eta_2)$  and 
  $y'_2\in \interior(y'_1\xi_2)$.  Notice $m(x_i)=m(x'_i)=m(y_i)=m(y'_i)=8$ for $i=1,2$, otherwise 
if,  say,  $m(x_1)=2$, then $\angle_{x_1}(x'_1, \xi_i)=\pi$ for $i=1,2$,    contradicting  to 
Lemma \ref{sec9.3l1}. Since $(x_1, x'_1, y'_1)$ has three even angles, $x_1y'_1$ contains no vertex indexed by $8$ in the interior.
 So we have $y'_2\in x_1\xi_2$.

First consider the case $y'_2=x_1$.  Then $y_2\not=y_1$ by Lemma \ref{sec9.3l2}.  
 Since $\angle_{x'_1}(x_1,\eta_1)=\pi$, we have $y_2\in \interior(y_1\eta_2)$. 
 The triangle $(x_1, y_1, y_2)$ has three even angles and so the side $y_1y_2$ contains 
 no vertex indexed by $8$ in the interior. 
Consider $x_2x'_2$.  Corollary \ref{sec9.3c1}  implies $x_2\not=x_1$  and 
$x'_2\in \interior(x'_1\eta_2)$. 
 Since $m(x'_2)=8$, we have $x'_2\in y_2\eta_2$. By the first paragraph  $x_2\in \interior(x_1\xi_1)$.
Since $m(x_2)=8$ and $x_1y'_1$ contains no vertex  indexed by $8$ in the interior, $x_2\in y'_1\xi_1$.
 Now $x_2y_1=x_2y'_1\cup y'_1y_1$ and $(x_2, y_1, x'_2)$ has three even angles. 
 The three triangles $(x_2, y_1, x'_2)$, $(x'_1, y'_1, x_1)$ and $(x_1, y_2, y_1)$ give rise
 to a loop with length $3\pi/4$ in $\Link(\De, y_1)$, which is impossible.

Now consider the case $y'_2\in \interior(x_1\xi_2)$.  By considering  $(x_1, y_1, y'_1)$
 we see $\angle_{x_1}(y_1, y'_1)=\pi/4$ and so $\angle_{x_1}(y_1, \xi_2)\ge 3\pi/4$. 
It follows that $y_2\in \interior(y_1\eta_2)$. 
Consider $Q=(x_1, y_1, y_2, y'_2)$.  $Q$ has one angle $\angle_{x_1}(y_1, \xi_2)\ge 3\pi/4$
 and the other three angles $\ge \pi/4$.  Proposition \ref{quar11}  implies $n(Q)\le 12$.
  On the other hand, there are at least 6 chambers in $S(Q)$ incident to $x_1$, at least 2 chambers incident to 
 each of $y'_2$,  $y_2$, $x'_1$, for a total of 12. It follows that $S(Q)$ is the union of 
 these 12 chambers. However, the boundary of this union is not a quadrilateral, contradiction.

\end{proof}

\noindent
{\bf{Proof of Proposition \ref{sec9p2}}}.
Suppose $\xi_1\xi_2\cap \eta_1\eta_2=\phi$.  
Let $y_i\in \eta_1\eta_2$ and $x_i\in \xi_1\xi_2$ ($i=1,2$)
 be such that $y_i\xi_i\cap \xi_1\xi_2=y'_i\xi_i$,
$x_i\eta_i\cap \eta_1\eta_2=x'_i\eta_i$ are rays and  $y_i\xi_i\cap \eta_1\eta_2=\{y_i\}$,
   $x_i\eta_i\cap  \xi_1\xi_2=\{x_i\}$.
Corollary  \ref{sec9.3c1}
and Lemma \ref{sec9.3l4}  imply that there are two cases:\newline
(1) there is some $i=1,2$ such that $\xi_1, y'_1, x_i, x_{i+1}, y'_2, \xi_2$
  and $\eta_i, x'_i, y_1, y_2, x'_{i+1}, \eta_{i+1}$  are in linear order;\newline
(2)  there is some $i=1,2$ such that $\xi_1, y'_1, x_i, x_{i+1}, y'_2, \xi_2$
  and $\eta_i, x'_i, y_2, y_1, x'_{i+1}, \eta_{i+1}$  are in linear order.

 We prove (2) is impossible, one similarly proves that 
      (1) is  impossible.    Suppose (2) holds.
Recall  $m(v)=8$ if $v$ is one of 
the vertices $x_i, x'_i, y_i, y'_i$ ($i=1,2$). Consider $Q=(y'_1,y_1, y_2, y'_2)$.
 Since all the angles of $Q$ are even, Proposition \ref{quar11}  implies 
  $n(Q)\le 24$.  We count the chambers in
$S(Q)$ by counting the  chambers   incident  to vertices   indexed by $8$. 
The vertices $x_i, x_{i+1}$ lie in the interior of $y'_1y'_2$, and so there are 8 chambers incident  to each of them. 
 There are also at least two chambers incident  to each of $y_1, y'_1, y_2, y'_2$.
We have exhibited   24 chambers in $S(Q)$. 
    It follows that   $n(Q)=24$ and $\angle_{y_1}(y'_1, y_2)=\pi/4$.
 Now consider $Q'=(y'_1, y_1, x'_{i+1}, x_{i+1})$. 
 The quadrilateral $Q'$ has an angle $\angle_{y_1}(y'_1, x'_{i+1})\ge 3\pi/4$.
It follows that $m(Q')\le 12$. On the other hand, there are at least 8 chambers
 in $S(Q')$ incident   to $x_i$ and at least 2 chambers incident  to each of $y'_1, y_1, x'_{i+1}, x_{i+1}$,
 for a total of 16, contradiction.

\qed

We next consider Type II geodesics.

\b{Prop}\label{sec9p3}
{Let $\De$ be a Fuchsian building with chamber $(2,3,8)$.
If two Type II  geodesics $\xi_1\xi_2, \eta_1\eta_2\subset \De^{(1)}$ are at different sides, 
  then  $\xi_1\xi_2\cap \eta_1\eta_2\not=\phi$.}

\end{Prop}

We start with  some  lemmas. 

\b{Le}\label{sec9.3l5}
{Let $\xi_1\xi_2$, $\eta_1\eta_2$ be two disjoint Type II geodesics.
 Suppose $\xi_1\xi_2$, $\eta_1\eta_2$  are at different sides. 
Let $y_i\in \eta_1\eta_2$ \e{($i=1,2$)} be such that $y_i\xi_i\cap \xi_1\xi_2=y'_i\xi_i$ is a
   ray  and  $y_i\xi_i\cap  \eta_1\eta_2=\{y_i\}$.
 Then $y_1\not=y_2$.}

\end{Le}

\b{proof}
Suppose $y_1=y_2$.  If $y'_1\not=y'_2$, then $y'_1\in \interior(y'_2\xi_1)$.   $(y'_1, y'_2, y_1)$ has two even 
 angles and all its sides are contained in Type II geodesics, 
    contradicting to Lemma \ref{sec9.1l1}. So we have $y'_1=y'_2$. In this case,
 $\angle_{y'_1}(\xi_1, y_1)=\angle_{y'_1}(\xi_2, y_1)=\pi$. 
 Let $x_i\in \xi_1\xi_2$ ($i=1,2$) be such that $x_i\eta_i\cap \eta_1\eta_2=x'_i\eta_i$ is a ray
  and  $x_i\eta_i\cap  \xi_1\xi_2=\{x_i\}$. 
  The uniqueness of geodesic implies $x'_i\in y_1\eta_i$.   Since  
 $\xi_1\xi_2$, $\eta_1\eta_2$  are at different sides,  we  have    $x'_i\not=y_1$. 
 Lemma \ref{sec9.1l1} applied to  $(x_i, x'_i, y_1)$   shows 
$\angle_{y_1}(\eta_i, y'_1)=\angle_{y_1}(x'_i, y'_1)=\pi/4$ or $\pi/3$. 
 Triangle inequality   then  implies $\angle_{y_1}(\eta_1,\eta_2)\le 2\pi/3$,
 contradicting to the fact that $y_1\in \eta_1\eta_2$.

\end{proof}

\b{Le}\label{sec9.3l6}
{Let $\xi_1\xi_2$, $\eta_1\eta_2$ be two disjoint Type II geodesics.
 Suppose $\xi_1\xi_2$, $\eta_1\eta_2$  are at different sides.
Let $y_1\in \eta_1\eta_2$,  $x_1\in \xi_1\xi_2$ be such that 
 $y_1\xi_1\cap \xi_1\xi_2=y'_1\xi_1$, $x_1\eta_1\cap \eta_1\eta_2=x'_1\eta_1$ are rays
  and  $y_1\xi_1\cap  \eta_1\eta_2=\{y_1\}$,  $x_1\eta_1\cap   \xi_1\xi_2=\{x_1\}$. 
Then $y'_1\in \interior(x_1\xi_1)$ and $x'_1\in \interior(y_1\eta_1)$.}

\end{Le}

\b{proof}
We show  (1) and (2) of Proposition \ref{tgeoddif}
cannot occur. 
First suppose (1) of Proposition \ref{tgeoddif}
occurs, that is, $y_1'\in\interior(x_1\xi_1)$  and $y_1\in \interior(x_1'\eta_1)$.
 In this case, all sides of $(x_1, y_1, y'_1)$ are 
  contained in Type II geodesics and one of them   $x_1y_1=x_1x'_1\cup x'_1y_1$ consists  of (at least) two edges.
 Lemma \ref{sec9.1l1} implies $m(x'_1)=2$. It follows that 
$\angle_{x'_1}(x_1, \eta_1)=\angle_{x'_1}(x_1, \eta_2)=\pi$, contradicting to Lemma \ref{sec9.3l5}.

Now suppose (2) of Proposition \ref{tgeoddif}
occurs, that is, $y_1'\in\interior(x_1\xi_1)$  and  $y_1=x_1'$.
The triangle $(x_1, y_1, y'_1)$ has an even angle at $y'_1$. 
Lemma \ref{sec9.1l1} implies $y'_1y_1$, $y'_1x_1$  are  edges in $\De$, $m(y'_1)=8$,  $m(y_1)=3$   
  and $\angle_{x_1}(y_1, y'_1)=\pi/3$.
Let $x_2\in \xi_1\xi_2$ be such that $x_2\eta_2\cap \eta_1\eta_2=x'_2\eta_2$ is  a ray
  and  $x_2\eta_2\cap  \xi_1\xi_2=\{x_2\}$. 
 Lemma \ref{sec9.3l5} implies $x_1\not=x_2$.
  The uniqueness of geodesic  implies   $x'_2\in x'_1\eta_2$. 
 First assume $x'_2=x'_1$. 
Since $\angle_{y'_1}(\xi_1, y_1)=\pi$, the uniqueness of geodesic implies 
 $x_2\in y'_1\xi_2$.  The fact that $\xi_1\xi_2$, $\eta_1\eta_2$  are at different sides
 implies $x_2\not=y'_1$. Hence $x_2\in \interior(x_1\xi_2)$. In this case,
 $(x_1, x_2, y_1)$ has an angle $\angle_{x_1}(x_2, y_1)\ge 2\pi/3$, contradicting to Lemma \ref{sec9.1l1}.
    Therefore we must have $x'_2\in \interior(x'_1\eta_2)$. 
 If $x_2\in \interior(y'_1\xi_1)$, then  the side
 $x_2y_1$  of  
 $(x_2, y_1, x'_2)$  contains $y'_1$  in the  interior  with $m(y'_1)=8$, contradicting to Lemma \ref{sec9.1l1}. 
If $x_2=y'_1$, then $(x_2, y_1, x'_2)$ has an even angle at $x'_2$ 
and Lemma \ref{sec9.1l1}  implies  $y'_1y_1$  consists of two edges, contradicting to the above observation that 
$y'_1y_1$ is an edge. 
The only remaining possibility is $x_2\in \interior(x_1\xi_2)$. 
In this case, the two angles of $(x_2, x_1, y_1, x'_2)$  at $x_1$ and $y_1$ are $2\pi/3$,
contradicting to Lemma \ref{sec9.1l2}.

\end{proof}

\b{Le}\label{sec9.3l7}
{Let $\xi_1\xi_2$, $\eta_1\eta_2$ be two disjoint Type II geodesics.
 Suppose $\xi_1\xi_2$, $\eta_1\eta_2$  are at different sides. 
Let $y_i\in \eta_1\eta_2$ \e{($i=1,2$)} be such that $y_i\xi_i\cap \xi_1\xi_2=y'_i\xi_i$ is a ray
  and $y_i\xi_i\cap  \eta_1\eta_2=\{y_i\}$.
 Then $y_1\not=y_2$  and 
$y'_1\in \interior(y'_2\xi_1)$.}

\end{Le}

\b{proof}
The claim  $y_1\not=y_2$  is the content of Lemma \ref{sec9.3l5}.  
Let $x_1\in \xi_1\xi_2$ be such that $x_1\eta_1\cap \eta_1\eta_2=x'_1\eta_1$ is a ray 
 and $x_1\eta_1\cap  \xi_1\xi_2=\{x_1\}$.
 Lemma \ref{sec9.3l6} implies $y'_1\in \interior(x_1\xi_1)$.   The same lemma applied to
   $y_2\xi_2$ and $x_1\eta_1$ shows $y'_2\in \interior(x_1\xi_2)$. Hence 
$y'_1\in \interior(y'_2\xi_1)$.

\end{proof}

\noindent
{\bf{Proof of Proposition \ref{sec9p3}}}.
Suppose $\xi_1\xi_2\cap \eta_1\eta_2=\phi$.
Let $y_i\in \eta_1\eta_2$ and $x_i\in \xi_1\xi_2$ ($i=1,2$)
 be such that $y_i\xi_i\cap \xi_1\xi_2=y'_i\xi_i$,
$x_i\eta_i\cap \eta_1\eta_2=x'_i\eta_i$ are rays  and    $y_i\xi_i\cap \eta_1\eta_2=\{y_i\}$,
   $x_i\eta_i\cap  \xi_1\xi_2=\{x_i\}$.
 Lemma \ref{sec9.3l6}   and  Lemma  \ref{sec9.3l7}  imply that there are two possibilities:\newline
(1) there is some $i=1,2$ such that $\xi_1, y'_1, x_i, x_{i+1}, y'_2, \xi_2$
  and $\eta_i, x'_i, y_1, y_2, x'_{i+1}, \eta_{i+1}$  are in linear order;\newline
(2)  there is some $i=1,2$ such that $\xi_1, y'_1, x_i, x_{i+1}, y'_2, \xi_2$
  and $\eta_i, x'_i, y_2, y_1, x'_{i+1}, \eta_{i+1}$  are in linear order.\newline
 We prove (1) is impossible, the proof that (2) is   impossible is similar.

We assume (1) holds.
We first observe that $m(v)\not=2$ if $v$ is one of the $x_i, y_i, x'_i,y'_i$:
if, for example,  $m(x'_1)=2$, then $\angle_{x'_1}(\eta_1, x_1)=\angle_{x'_1}(\eta_2, x_1)=\pi$,
 contradicting to Lemma \ref{sec9.3l5}. Consider the quadrilateral $(x_i, x_{i+1}, x'_{i+1}, x'_i)$.
 The two angles at $x'_i$ and $x'_{i+1}$ are even, and the side $x'_ix'_{i+1}$ contains the
 two vertices  $y_1$, $y_2$  in the interior  with $m(y_1), m(y_2)\not=2$. Lemma \ref{sec9.1l3} implies 
$m(v)\not=8$ for every  vertex $v\in \interior(x'_ix'_{i+1})$ 
   and there is a vertex $v'\in x_ix_{i+1}$ with $m(v')=8$.
But the same lemma applied to $(y'_1, y'_2, y_2, y_1)$ implies 
$m(v)\not=8$ for every   vertex $v\in \interior(y'_1y'_2)$.   Here we have a contradiction since 
$x_ix_{i+1}\subset \interior(y'_1y'_2)$ and $v'\in x_ix_{i+1}$ with $m(v')=8$.

\qed

\subsubsection{The image of an apartment} \label{imofasecspecial1}

Let $\De_1$ and $\De_2$ be two Fuchsian buildings with chamber $(2,3,8)$, and $h:\p\De_1\ra \p\De_2$  a homeomorphism that preserves the combinatorial cross ratio almost everywhere.
  There are integers $p\not=q\ge 2$ such that 
all edges in $\De_1$ incident to  vertices   $v$ with $m(v)=3$  are contained in exactly
$p+1$ chambers, and all other edges are contained in exactly $q+1$ chambers.   By Lemma \ref{redaprtowh0}
  there are two possibilities: \newline
(1) A geodesic $c\subset \De_1^{(1)}$ is Type I if and only if its image $c'$ is Type I. Hence 
all edges in $\De_2$ incident to   vertices   $w$ with $m(w)=3$   are contained in exactly
$p+1$ chambers, and all other edges are contained in exactly $q+1$ chambers;\newline
(2)  A geodesic $c\subset \De_1^{(1)}$ is Type I if and only if its image $c'$ is Type II. Hence
 all edges in $\De_2$ incident to  vertices   $w$ with $m(w)=3$ 
     are contained in exactly
$q+1$ chambers, and all other edges are contained in exactly $p+1$ chambers.

\b{Le}\label{sec9.4l1}
{Let $A$ be an  apartment in $\De_1$, and $v\in A$ a vertex with $m(v)=8$.  
Then the images of the 4 Type I geodesics \e{(}in $A$\e{)} through $v$  intersect in
  a unique vertex  of  $\De_2$.}

\end{Le}

\b{proof}
Let $\gamma_1,  \gamma_2, \gamma_3\subset \De_2^{(1)}$ be three geodesics of the same type. 
Suppose any two of them intersect in a vertex and $\gamma_1\cap   \gamma_2\cap  \gamma_3=\phi$.
Denote $x_i=\gamma_{i-1}\cap \gamma_{i+1}$ ($i\mod 3$) and apply Proposition \ref{tri2}  to 
 $T=(x_1, x_2, x_3)$.  
If the $\gamma_i$'s are Type I,  then all the angles of $T$ are $\pi/4$; if the $\gamma_i$'s are Type II,
 then the angles of $T$ are $\pi/3$, $\pi/3$ and $\pi/4$. In any case,  the angles of $T$
 are $<\pi/2$.

Now  let  $c_i\subset A$ ($i=1,2,3,4$) be the 4 Type I geodesics  through $v$, and $c'_i$ their images.
  Then the  geodesics  $c'_i$, $1\le i\le 4$  have the same type.  
 Proposition \ref{sec9p1}
  implies   $c'_i\cap c'_j\not=\phi$ for $1\le i,j\le 4$. 
 First suppose $c'_1$, $c'_2$, $c'_3$ intersect in a point $w$ and $w\notin c'_4$.
Let $y_i=c'_4\cap c'_i$ for $1\le i\le 3$.  We may assume $y_2$ lies between $y_1$ and $y_3$.
Consider $(w,y_1,y_2)$  and  $(w,y_2,y_3)$.   The first paragraph implies 
$\angle_{y_2}(w, y_1)<\pi/2$  and $\angle_{y_2}(w, y_3)<\pi/2$. The triangle inequality implies that
$\angle_{y_2}(y_1, y_3)<\pi$, contradicting to $y_2\in \interior(y_1y_3)$.

Next we assume any three of $c'_i$, $1\le i\le 4$ have empty intersection. 
Let $y_i$, $i=1,2,3$ be as above. We may assume $y_2$ lies between $y_1$ and $y_3$.
Let  $w=c'_1\cap c'_3$  and  $w_1=c'_1\cap c'_2$. 
 Consider $(w_1, y_1, y_2)$. The first paragraph implies $\angle_{y_2}(w_1,  y_1)<\pi/2$.
In particular, $m(y_2)\not=2$.  Since $y_2\in \interior(y_1y_3)$, Proposition \ref{tri2}
 applied to $(w, y_1, y_3)$ implies $m(y_2)=2$, a contradiction.

\end{proof}

Let $A$ be an apartment in $\De_1$. Type I geodesics in $A$ divide $A$ into triangular regions, each of which
 is the union of 6 chambers.  We consider the  triangulation of $A$ whose 1-skeleton
is  the union of Type I geodesics.  
We denote by  $A^{[i]}$, $i=0, 1, 2$ the $i$-skeleton of this triangulation.   
 Notice $A^{[0]}$ is the set of vertices $v$ in $A$ with $m(v)=8$.

\b{Prop}\label{sec9.4p1}
{Let $\De_1$ and $\De_2$ be two Fuchsian buildings with chamber $(2,3,8)$, and $h:\p\De_1\ra \p\De_2$  a homeomorphism that preserves the combinatorial cross ratio almost everywhere.  
   Then a  geodesic $c\subset \De_1^{(1)}$ is Type I if and only if its image $c'$ is Type I.

}

\end{Prop}

\b{proof} 
We suppose the opposite holds, that is, a  geodesic $c\subset \De_1^{(1)}$ is Type I if and only if its image $c'$ is Type II.
Fix  an apartment  $A$  of $\De_1$.  
First suppose there is a  2-simplex  $D$ in $A^{[2]}$ such that $c'_1\cap c'_2\cap c'_3\not=\phi$,
  where $c_1, c_2, c_3\subset A$  are the three geodesics that contain the edges of $D$. 
Then Lemma \ref{sec9.4l1}
and the proof of Proposition \ref{onerionex}
imply  that   the images of all the Type I geodesics in $A$ intersect in a unique vertex of $\De_2$,
  which is impossible.  Hence for any   2-simplex   $D$ in $A^{[2]}$, if $c_1, c_2, c_3\subset A$  
  are the three geodesics that contain the edges of $D$, then $c'_1\cap c'_2\cap c'_3=\phi$. 
Since $c'_1, c'_2, c'_3$ are Type II geodesics,  Lemma \ref{sec9.1l1} implies they form a triangle in $\De_2$ which has two angles 
 $=\pi/3$ and one angle $=\pi/4$.

Let $D_1$,  $D_2$ be two  2-simplices in $A^{[2]}$ that share an edge. Let $v$ be one of the common vertices
   of $D_1$ and $D_2$,
 $c_2$ the geodesic in $A$ containing the common edge of $D_1$ and $D_2$, $c_1\subset A^{(1)}$ 
 the geodesic containing the other edge of $D_1$ incident to  $v$, 
$c_3\subset A^{(1)}$ 
 the geodesic containing the other edge of $D_2$ incident to  $v$, 
and $\gamma_i$ ($i=1,2$)  the geodesic containing the third edge of $D_i$.
 Denote $w=c'_1\cap c'_2$,  $x_i=c'_i\cap \gamma'_i$ ($i=1,2$),  $x_3=c'_3\cap \gamma'_2$.  By the last paragraph,
  $T_i=(w, x_i, x_{i+1})$ is a triangle  with two vertices indexed by $3$. We may assume $m(w)=3$.  
 Notice $S(T_1)\cap S(T_2)\subset \De_2$ is a convex  subcomplex containing $wx_2$.  We claim   $S(T_1)\cap S(T_2)=wx_2$.
By  Lemma \ref{sec9.1l1}
   $S(T_i)$ is the union of two chambers. If $S(T_1)\cap S(T_2)\not=wx_2$, then $S(T_1)\cap S(T_2)$ contains 
 a chamber incident to $w$, which implies $c'_1\cap c'_3$  contains an edge of $\De_2$. This is impossible since 
 $c'_1$ and $c'_3$ are at different sides, see Lemma \ref{diinet1}.

Let  $c_4\subset A$ be the remaining  Type  I geodesic    through 
$v$.    There is a 2-simplex   $D_i$ ($i=3,4$) of $A^{[2]}$
  incident to $v$ and  with  edges lying on $c_i$ and $c_{i+1}$ ($i\mod 4$).  Let $\gamma_i$ be the geodesic containing the third 
 edge of $D_i$, $i=3,  4$.  Denote $x_4=c'_4\cap \gamma'_4$ and  $x_5=c'_1\cap \gamma'_4$. 
   There are   4 triangles $T_i=(w, x_i, x_{i+1})$, $1\le i\le 4$.  
Notice $x_1, x_5\in c'_1$ and $d(w, x_1)=d(w, x_5)$.   Hence there are two possibilities: either $x_5=x_1$ or 
   $w$ is the midpoint of $x_1x_5$.  If  $x_5=x_1$, then the claim in the last paragraph implies that the 4 triangles  $T_i$
 give rise to a geodesic loop
 with length $\frac{\pi}{3}\times 4$ in the $\CAT(1)$ space  $\Link(\De_2, w)$, a contradiction.
 If $w$ is the midpoint of $x_1x_5$,   then  the same claim implies that the $T_i$'s give rise to an edge path with length 4
  between two opposite vertices of the generalized polygon $\Link(\De_2, w)$.  This is impossible because 
  $\Link(\De_2, w)$ is a generalized 3-gon and opposite vertices have different types.

\end{proof}

Lemma \ref{sec9.4l1} implies for any apartment $A\subset \De_1$, 
 there is a well-defined map $f_A: A^{[0]}\ra \De^{(0)}$, where for any $v\in A^{[0]}$, 
 $f_A(v)$ is the unique intersection point of the images of the 4 Type I geodesics in $A$ through $v$. 
We denote $f_A(v)$ by $v'$.

\b{Prop}\label{sec9.4p2}
{Let $A$ be an apartment in $\De_1$.  Then there is a  unique  apartment $A'$ in $\De_2$  
with the 
following properties:\newline
\e{(1)} $f_{A}(A^{[0]})={A'}^{[0]}$;\newline
\e{(2)}  the map $f_{A}: A^{[0]}\ra {A'}^{[0]}$  extends to an isomorphism $g_A: A\ra A'$.}

\end{Prop}

\b{proof}
Let $D$ be a 2-simplex in $A^{[2]}$. Denote by $c_1, c_2,  c_3\subset A^{(1)}$ 
the three geodesics   that contain the three edges of $D$. If $c'_1\cap c'_2\cap  c'_3\not=\phi$, then 
 Lemma \ref{sec9.4l1}  and the proof of Proposition \ref{onerionex}
imply  that  the images of all the Type I geodesics in $A$ have nonempty intersection, impossible.
Hence, $c'_1\cap c'_2\cap  c'_3=\phi$. 
Denote  $x_i=c_{i-1}\cap c_{i+1} $ ($i\mod 3$).
Since  $c'_i$ is    a Type I geodesic, 
 $T'=(x'_1, x'_2, x'_3)$  has three even angles  and Proposition \ref{tri2}
implies that $S(T')$ is isometric to $D$.   We  define $g_D$ to be   the unique isometry  extending the map 
 $f_{A}|_{\{x_1, x_2, x_3\}}$. It is clear that if  $D_1$, $D_2\subset A^{[2]}$ are two 
2-simplices  with $D_1\cap D_2\not=\phi$,  then 
 $g_{D_1}$ and $g_{D_2}$ agree on $D_1\cap D_2$.
     Hence we  have a map $g_A:  A\ra \De_2$, where $g_A|_{D}=g_D$ for each 2-simplex    $D\subset A^{[2]}$.

Consider two arbitrary 2-simplices $D_1$, $D_2\subset A^{[2]}$  that share an edge.
Let $vx_2$ ($v, x_2\in  A^{[0]}$) be the common edge of $D_1$ and $D_2$,
$x_1$ the third vertex of $D_1$  and $x_3$ the third vertex of $D_2$.   Let  
$c_i\subset A$ ($i=1,2,3$) be the geodesic containing $vx_i$, and $\gamma_j\subset A$ ($j=1,2$) 
   the geodesic containing $x_jx_{j+1}$.     Let $T'_i=(w, x'_i, x'_{i+1})$ ($i=1,2$).  Notice 
$\image(g_{D_i})=S(T'_i)$.  
We claim  $\image(g_{D_1})\cap \image(g_{D_2})=v'x'_2$.  
The claim implies  that $g_A$ is an isometry into   $\De_2$  and hence 
   $\image(g_A)$ is an apartment in $\De_2$.  
Next we prove the claim.

Suppose  the claim is false. 
Recall  $S(T'_i)\subset \De_2$ is a convex subcomplex and is the union of 6 chambers.   We observe that 
$S(T'_1)\cap  S(T'_2)$ is the union of the two chambers in  $S(T'_1)$  that  have edges lying on $v'x'_2$.
Otherwise $S(T'_1)\cap  S(T'_2)$  would contain a nontrivial segment of $c'_1$ or $\gamma'_1$, which implies 
 $c'_1\cap c'_3$ or $\gamma'_1\cap \gamma'_2$ contains a nontrivial interval, contradicting to the fact that 
$c'_1$, $c'_3$ are at different sides and $\gamma'_1$,  $\gamma'_2$ are also at different sides.
It follows that $c'_1$ and $c'_3$ make an angle $\pi/4$. 
 Let  $ c_4\subset A$ be the  fourth  Type I  geodesic through $v$.  
 Also let $D_i\subset A^{[2]}$ ($i=3,4$) be the 2-simplex  incident to $v$, lying between $c_i$ and $c_{i+1}$ ($i+1\mod 4$),
  and sharing an edge with $D_{i-1}$. 
 Let $\gamma_i\subset A$, $i=3,4$,  be the geodesic containing the third edge of $D_i$.
 Denote $x_4=c_4\cap \gamma_4$, $x_5=c_1\cap \gamma_4$.
Notice  $v', x'_1, x'_5\in c'_1$ and $d(v', x'_1)=d(v', x'_5)$.

Let $T'_i=(v', x'_i, x'_{i+1})$ ($i=3,4$).  Note  each $T'_i$ has angle $\pi/4$ at $v'$, and 
  $v'x'_{i+1}\subset S(T'_i)\cap S(T'_{i+1})$.  Since $c'_1$ and $c'_3$ make an angle $\pi/4$,  the 
triangle inequality implies $\angle_{v'}(x'_1, x'_5)\le 3\pi/4$.   Hence   $x'_1=x'_5$.
      If $S(T'_2)\cap S(T'_3)=v'x'_3$, then 
$$(S(T'_1)-S(T'_1)\cap S(T'_2))\cup (S(T'_2)-S(T'_1)\cap S(T'_2))\cup S(T'_3)$$
 gives rise to an injective path in $\Link(\De_2, v')$ 
with length $\pi/2$ from
 the initial direction of $v'x'_1$ to that of $v'x'_4$.   On the other hand,  $\angle_{v'}(x'_4, x'_5)=\pi/4$ and $x'_1=x'_5$,
  a contradiction.
 Hence $S(T'_2)\cap S(T'_3)\not=v'x'_3$.  We have proved the implication:
if $S(T'_1)\cap S(T'_2)$ is not  a segment, then $S(T'_2)\cap S(T'_3)$  is not  a segment.  
By working with the  2-simplices of $A^{[2]}$ around  $x_2$  instead of $v$ 
and using gallery connectedness one 
   concludes  that 
 for any two  2-simplices $E_1, E_2\subset A^{[2]}$ sharing an edge, 
the intersection $\image(g_{E_1})\cap \image(g_{E_2})$  is not   a  segment. 
The above argument also shows that for any Type I geodesic  $c\subset A$, and any three consecutive vertices
 $x,y,z\in c$ with $m(x)=m(y)=m(z)=8$,   we have $x'=z'$.   Consequently,  
 for any Type I geodesic $c\subset A$, there are infinitely many $v_i\in c\cap A^{[0]}$ ($i\ge 1$)
 such that the family $\{\gamma': \gamma \subset A\; \text{is Type I and } \gamma\owns v_i \;\; \text{for some}\; i\, \}$ 
  have nonempty intersection.  We get a contradiction as in the proof of Proposition \ref{onetoneq}.

\end{proof}

\b{Cor}\label{sec9.4c2}
{Let $A$ be an apartment in $\De_1$.  Then for any vertex $v\in A$, the  geodesics in ${\mathcal{D}}'_{A,v}$
  intersect  in a unique  vertex of $\De_2$.}

\end{Cor}

\b{proof} 
Since $g_A:A\ra A'$ is an isomorphism, for any vertex $v\in A$,  
the  $g_A$-images of the  geodesics in ${\mathcal{D}}_{A,v}$  
 intersect in a  unique  vertex.  Hence it suffices  to show that  for any $\xi_1\xi_2\subset A^{(1)}$, we have
 $g_A(\xi_i)=\xi'_i$ ($i=1,2$),  where  $g_A: \p A\ra \p A'$ also denotes the boundary map of 
   $g_A:A\ra A'$.  
By the definition of $g_A$   this claim  holds for all Type I geodesics in
  $A$.   
For an arbitrary geodesic $\xi_+\xi_-\subset A^{(1)}$,    one can find a sequence of  Type I geodesics 
 $\eta_i\omega_i\subset A^{(1)}$, $i\in \Z$,   such that $\xi_+\xi_-\cap \eta_i\omega_i\not=\phi$,
  $\eta_i\omega_i\ra \xi_+$  as $i\ra +\i$ and    $\eta_i\omega_i\ra \xi_-$  as $i\ra -\i$. 
It follows that $\eta_i\ra \xi_+$  as $i\ra +\i$  and $\eta_i\ra \xi_-$  as $i\ra -\i$.
 Since  $h:\p\De_1\ra \p\De_2$ is a homeomorphism,
 we have $\xi'_+=h(\xi_+)=h(\lim_{i\ra +\i}\eta_i)=\lim_{i\ra +\i}h(\eta_i)
=\lim_{i\ra +\i}\eta'_i
=\lim_{i\ra +\i}g_A(\eta_i)
=g_A(\xi_+)$.
Similarly $\xi'_-=g_A(\xi_-)$  and we are done.

\end{proof}

\subsubsection{Geodesics through a fixed vertex} \label{subsecspecial}

In this section we show that  for any vertex $v\in \De_1$,
the   geodesics in  ${\mathcal{D}}'_v$  intersect  in  a unique vertex.

By Corollary \ref{sec9.4c2}
for  any vertex  $v\in \De_1$  and any apartment $A$ containing $v$, 
the   geodesics in  ${\mathcal{D}}'_{A,v}$  intersect  in  a unique vertex  $v_A$ of $\De_2$. 
The proof of Corollary \ref{sec9.4c2}  actually shows $v_A=g_A(v)$, where $g_A: A\ra A'$ is the isomorphism
in Proposition \ref{sec9.4p2}.

\b{Prop}\label{sec9.5p1}
{Let $\De_1$ and $\De_2$ be two Fuchsian buildings with chamber $(2,3,8)$, and $h:\p\De_1\ra \p\De_2$  a homeomorphism that preserves the combinatorial cross ratio almost everywhere.  Then for any vertex $v\in \De_1$ with $m(v)=8$,
  the   geodesics in  ${\mathcal{D}}'_v$  intersect  in  a unique vertex of $\De_2$.}

\end{Prop}

\b{proof} 
By Lemma \ref{redaprtowh}, we only need to show that $v_{A_1}=v_{A_2}$  holds for any two apartments 
 $A_1, A_2\subset \De_1$    containing   $v$   with  
 $\Link(A_1, v)\cap \Link(A_{2}, v)$     a  half apartment in $\Link(\De_1, v)$.
Let $B_i=\Link(A_i,v)$ ($i=1,2$), and $\omega_1$, $\omega_2$ the two endpoints of $B_1\cap B_2$.
Denote by $m$ the midpoint of $B_1\cap B_2$,  and  $m_i\in B_i$ the point in $B_i$ opposite to $m$.
Since $m(v)=8$, $m_1$,   $m_2$ and $m$  are  vertices. Notice $m_1$ and $m_2$ are opposite to each other.
Let $\xi_1, \xi_3\in \p A_1$ and $\xi_2, \xi_4\in \p A_2$ such that
 $v\xi_1$, $v\xi_2$ have initial direction $\omega_1$ and 
  $v\xi_3$, $v\xi_4$ have initial direction $\omega_2$. Similarly
let $\eta_1,\eta_3\in \p A_1$ and $\eta_2, \eta_4\in \p A_2$ such that 
$v\eta_3$, $v\eta_4$ have initial direction $m$,   and 
  $v\eta_i$ ($i=1,2$) has initial direction $m_i$.
Then we have $\eta_i\eta_j=v\eta_i\cup v\eta_j$  for $i\not=j$, $\{i,j\}\not=\{3,4\}$.
Lemma \ref{sec9.1l12}  implies  $\xi_1\xi_3$ and $\eta_i\eta_j$ ($1\le i\not=j\le 4$, $\{i,j\}\not=\{3,4\}$) 
   are at different sides, and $\xi_2\xi_4$ and $\eta_i\eta_j$ are also  
    at different sides.

 We claim   the three geodesics $\eta'_1\eta'_2$,  $\eta'_1\eta'_3$,
$\eta'_3\eta'_2$  have nonempty intersection.  The rest of the proof is the same as in 
  Proposition \ref{sec8.3p1}.  Suppose the claim is not true. Then Lemma \ref{threeev}
  implies     there are  vertices $v_1, v_2,  v_3$ such that 
$\eta'_1\eta'_2\cap \eta'_1\eta'_3=v_1\eta'_1$, $\eta'_2\eta'_1\cap \eta'_2\eta'_3=v_2\eta'_2$, 
$\eta'_1\eta'_3\cap \eta'_2\eta'_3=v_3\eta'_3$; furthermore $T:=(v_1, v_2,v_3)$ has three even angles and $S(T)$
is the union of 6 chambers,   as shown in  Figure \ref{238}(h).  In this case,  
  $\xi'_1\xi'_3$ and $\eta'_i\eta'_j$ ($1\le i\not=j\le 3$) are   Type I   and 
   are at different sides.  Proposition \ref{sec9p1}  implies 
 that  $\xi'_1\xi'_3\cap \eta'_i\eta'_j$ is a vertex  for $1\le i\not=j\le 3$. 
Denote $w_i=\xi'_1\xi'_3\cap \eta'_{i+1}\eta'_{i+2}$ ($i\mod 3$).  
 Since $\xi'_1\xi'_3$ and $\eta'_i\eta'_j$  make  even angles
  and  are  at different sides, $m(w_i)\not=2$.  Hence  $w_i\in v_{i+1}\eta'_{i+1}$ or $w_i\in v_{i+2}\eta'_{i+2}$.
 We may assume $w_1\in v_3\eta'_3$.  If  $w_3\in v_1\eta'_1$, then 
 $\xi'_1\xi'_3\cap \eta'_1\eta'_3\supset w_1w_3$ is not a vertex,
 contradiction.  Similarly we obtain a contradiction if $w_3\in v_2\eta'_2$.


\end{proof}

\b{Prop}\label{sec9.5p2}
{Let $\De_1$ and $\De_2$ be two Fuchsian buildings with chamber $(2,3,8)$, and $h:\p\De_1\ra \p\De_2$  a homeomorphism that preserves the combinatorial cross ratio almost everywhere.
  Then for any vertex $v\in \De_1$ with $m(v)=3$,
the   geodesics in  ${\mathcal{D}}'_v$  intersect  in  a unique vertex of $\De_2$.}

\end{Prop}

\b{proof} Let $v\in \De_1$ be  a  vertex  with $m(v)=3$.
We only need to show that $g_{A_1}(v)=g_{A_2}(v)$  holds for any two apartments 
 $A_1, A_2\subset \De_1$    containing   $v$   with  
 $\Link(A_1, v)\cap \Link(A_{2}, v)$     a  half apartment in $\Link(\De_1, v)$.
Let $B_i=\Link(A_i,v)$  ($i=1,2$), and $a$, $b$ the two endpoints of $B_1\cap B_2$.
 There are two edges $e_1=vv_1$ and $e_2=vv_2$ that have initial directions  $a$ and $b$ respectively. 
Note $\{m(v_1), m(v_2)\}=\{2,8\}$.  
We may assume 
 $m(v_1)=8$ and $m(v_2)=2$. Let $c, d$ be the other two vertices on $B_1\cap B_2$ such that
  $c$ lies between $a$ and $d$.  Let $vv_3$, $vv_4$ be the edges that have initial directions 
$c$ and $d$ respectively.  Then $m(v_4)=8$, and $T_3=(v, v_1, v_3)$, $T_4=(v, v_3, v_4)$
 are (boundaries of) chambers. Let $X=S(T_3)\cup S(T_4)$.  Notice 
$X$ is a convex subcomplex of $\De_1$, and $X\subset A_1\cap A_2$.

Since  $g_{A_i}: A_i\ra A'_i$  is an isomorphism,      $g_{A_i}(X)\subset \De_2$ is  a convex subcomplex.
Proposition \ref{sec9.5p1}  and  the fact 
  $m(v_1)=m(v_4)=8$     imply that   $g_{A_1}(v_1)=g_{A_2}(v_1)$
 and $g_{A_1}(v_4)=g_{A_2}(v_4)$.  As $v_3$ is the midpoint of $v_1v_4$, we also have 
$g_{A_1}(v_3)=g_{A_2}(v_3)$.
Set $v'_1=g_{A_1}(v_1)$, $v'_3=g_{A_1}(v_3)$, 
$v'_4=g_{A_1}(v_4)$.
Hence the intersection of the two convex subcomplexes $g_{A_1}(X),  g_{A_2}(X)\subset \De_2$
 contains  $v'_1v'_4$. 
  It is easy to see that either $g_{A_1}(X)\cap  g_{A_2}(X)=v'_1v'_4$
  or  $g_{A_1}(X)=g_{A_2}(X)$.   If $g_{A_1}(X)=g_{A_2}(X)$, then $g_{A_1}(v)=g_{A_2}(v)$
  and  we are done. We shall show that $g_{A_1}(X)\cap  g_{A_2}(X)=v'_1v'_4$
does not hold.

Let $\xi_i\in \p A_i$ ($i=1,2$) be the point such that the initial direction of $v\xi_i$  is 
   opposite to $c$ (the initial direction of $vv_3$).
Then $\angle_v(\xi_1, \xi_2)=2\pi/3$. Notice $\xi_1\xi_2\not\subset \De_1^{(1)}$, otherwise
 there would be a triangle with three even angles and one of the angles $=2\pi/3$, 
 contradicting to  Proposition \ref{tri2}.

Suppose $g_{A_1}(X)\cap  g_{A_2}(X)=v'_1v'_4$.
  Since  $m(v'_3)=m(v_3)=2$,      $\angle_{v'_3}(g_{A_1}(v), g_{A_2}(v))=\pi$.
 Hence   $v'_3\xi'_1\cup v'_3\xi'_2=v'_3g_{A_1}(\xi_1)\cup v'_3g_{A_2}(\xi_2)\subset \De_2^{(1)}$  is 
 the geodesic from $\xi'_1$ to $\xi'_2$.
It follows that $\xi_1\xi_2\subset\De_1^{(1)} $, contradicting to the preceding paragraph.

\end{proof}

\b{Prop}\label{sec9.5p3}
{Let $\De_1$ and $\De_2$ be two Fuchsian buildings with chamber $(2,3,8)$, 
and $h:\p\De_1\ra \p\De_2$  a homeomorphism that preserves the combinatorial cross ratio almost everywhere.
  Then for any vertex $v\in \De_1$ with $m(v)=2$,
the   geodesics in  ${\mathcal{D}}'_v$  intersect  in  a unique vertex of $\De_2$.}

\end{Prop}

\b{proof}
Let $A_1, A_2\subset \De_1$   be two apartments  containing   $v$   with  
 $\Link(A_1, v)\cap \Link(A_{2}, v)$     a  half apartment in $\Link(\De_1, v)$.
Let $B_i=\Link(A_i,v)$($i=1,2$), and $a$, $b$ the two endpoints of $B_1\cap B_2$.
There are two edges $e_1=vv_1, e_2=vv_2\subset A_1\cap A_2$ that have initial directions  $a$ and $b$ at $v$
  respectively. 
Note $m(v_1)=m(v_2)=3$ or $8$. By Propositions \ref{sec9.5p1}, \ref{sec9.5p2},
$g_{A_1}(v_1)=g_{A_2}(v_1)$  and $g_{A_1}(v_2)=g_{A_2}(v_2)$. 
Since $g_{A_i}(v)$ is the only vertex in the interior of $g_{A_i}(v_1)g_{A_i}(v_2)$,
 we have $g_{A_1}(v)=g_{A_2}(v)$.

\end{proof}



\subsection{Triangles and quadrilaterals  }\label{triqua}

In this section we prove Propositions \ref{tri1}, \ref{tri2}    and    \ref{quar11}.
  We have defined 
 triangles and quadrilaterals    in Section  \ref{triqua1sta}.

\subsubsection{Support sets of   triangles and quadrilaterals} \label{supor}

Let $X$ be a locally finite $\CAT(0)$ $2$-complex, e.g.,  a   Fuchsian 
   building.
For each $x\in X$, let $\log_x: X-\{x\}\ra Link(X,x)$ be the map that sends each $y\not=x$
 to the initial direction of $xy$ at $x$.  

\begin{Def}\label{supp}
{Let  $c\subset X$ be a subset that is   homeomorphic to a circle. 
The \emph{support set} $\supp(c)$ of $c$ is the set of $x\in X-c$ such that  $\log_x(c)$  
   represents a  nontrivial    class   in 
$H_1(\Link(X, x))$.}
\end{Def}

The goal of this section is to show that a 
 quadrilateral homeomorphic to a circle   must  bound a compact surface in $X$.  
The triangle case was proved in \cite{X}.

\b{Prop}\label{trisupp}\e{(\cite{X})}
{Let $X$  be a   
  locally finite  $\CAT(0)$   $2$-complex  
 and  $c\subset X$ a triangle  homeomorphic to a circle.   
  Then $\overline{\supp(c)}$
  is a compact surface with boundary $c$   and  $\overline{\supp(c)}=\supp(c)\cup c$.}
\end{Prop}

Given a triangle or quadrilateral  $c$, we subdivide the 2-complex $X$ such that 
 $c$ becomes part of the  1-skeleton of $X$. Thus if $c$ is a triangle  
homeomorphic to a circle and $x,y,z$  its   corners, then the segment from
 $\log_x(y)$ to  $\log_x(z)$  is an  edge path in the link $\Link(X, x)$. 
Note $\angle_x(y,z)<\pi$ since $c$ is homeomorphic to a circle. 
 Each edge in $\log_x(y)\log_x(z)$  corresponds to a 2-cell incident to $x$. Let 
$U(c, x)$ be the union of these 2-cells.   The proof  in \cite{X}    shows that 
   $U(c, x)$
  is  a  neighborhood of  $x$ in  $\overline{\supp(c)}$.

\begin{Le}\label{s1}
{Let $c\subset X$ be a triangle or  quadrilateral  homeomorphic to  a  circle,  and  
 $x\in X-c$. Then $x\in \supp(c)$ if and only if  $\log_x(c)$  is homotopic to a simple loop  with 
 length $<4\pi$ 
in $\Link(X,x)$. }
\end{Le}

\b{proof} 
One direction is clear  since  $\Link(X,x)$ is a finite graph. For the other direction let $x\in \supp(c)$. 
Then $\log_x(c)$ is homotopically nontrivial in the finite $\CAT(1)$ graph    $\Link(X,x)$. 
Notice if $x_1x_2$ is a geodesic segment which does not contain $x$, then
$\log_x(x_1)\log_x(x_2)$  has length  $<\pi$
  and  $\log_x(x_1x_2)$   is  homotopic to  $\log_x(x_1)\log_x(x_2)$.
It follows that 
$\log_x(c)$ is homotopic  to  a  closed path with 
 length  $<4\pi$.  The unique  closed geodesic $\sigma\subset  Link(X,x)$
in the free homotopy class of 
$\log_x(c)$   has the shortest length in the class. It follows that 
$\length(\sigma)<  4\pi  $.  Since  a closed geodesic 
 in a  $\CAT(1)$ space   has length at least $2\pi$,   $\sigma$   must  be  a  simple 
 closed geodesic.

\end{proof}

Lemma \ref{s1}  and the arguments in \cite{X}  imply the following:

\b{Le}\label{s2}
{Let $c$ be a quadrilateral in $X$ that is homeomorphic to a circle.
  Then  $\supp(c)$ is a topological surface,  $\overline{\supp(c)}$  is compact 
  and  $\overline{\supp(c)}\subset c\cup \supp(c)$.}

\end{Le}

Now we are ready to prove 

\b{Prop}\label{quipp}
{Let $X$  be a   
  locally finite  $\CAT(0)$   $2$-complex  
 and  $c\subset X$ a quadrilateral 
 homeomorphic to a circle. 
  Then $\overline{\supp(c)}$  is  a  compact surface with boundary $c$
   and  
$\overline{\supp(c)}=\supp(c)\cup c$.}

\end{Prop}

\b{proof}We shall show that each point $p\in c$ has a neighborhood $U$ such that $U\cap \overline{\supp(c)}$
  is homeomorphic to   a neighborhood of the origin in the closed upper half plane. 
Let $x, y, z, w\in c$ be the  4   corners of $c$ in cyclic order. 
Let $c_1=xy\cup yz\cup zx$ and $c_2=wx\cup xz\cup zw$.  We first consider 
    the case when $p\in c$ is one of corners, say, $p=w$. 
We may assume $w\notin   c_1$, otherwise $c$ is a triangle and the proposition  follows from Proposition \ref{trisupp}.
We orient $c_1$ and $c_2$ so that they have opposite orientations   on $xz$.  We also orient $c$ so that 
 $c$ is homotopic to $c_1*c_2$.  First suppose $w\notin \supp(c_1)$.  By Lemma \ref{s2} and the remark before
Lemma \ref{s1}  there is a neighborhood $U$ of $w$ in $X$ such that $U\cap \overline{\supp(c_2)}=U(c_2, w)\cap U$
  and $ U\cap \overline{\supp(c_1)}=\phi$.  Since $c$ is homotopic to $c_1*c_2$,   it follows from   
 the definition of support  set   that $U\cap \overline{\supp(c)}=U(c_2, w)\cap U$.

Now suppose  $w\in \supp(c_1)$.  Recall (see \cite{X}) that 
  the closed path $\log_w(c_1)$ is homotopic to a simple loop ${c_1}_w\subset  \Link(X,w)$
 with length $<3\pi$. The simple  loop ${c_1}_w$
can be constructed as follows. The three segments $\log_w(x)\log_w(y)$, $\log_w(y)\log_w(z)$,
$\log_w(z)\log_w(x)$  all have length $<\pi$ and their  intersections 
are segments   (possibly   degenerate): 
  there are $a_0,  b_0, c_0\in  \Link(X,w)$
  with  
$$\log_w(x)\log_w(y)\cap \log_w(x)\log_w(z)=\log_w(x)a_0,$$
$$\log_w(y)\log_w(x)\cap \log_w(y)\log_w(z)=\log_w(y)b_0$$  and 
$$\log_w(z)\log_w(y)\cap \log_w(z)\log_w(x)=\log_w(z)c_0.$$
Then ${c_1}_w=a_0b_0\cup b_0c_0\cup c_0a_0$. Let $U(c_1,w)$ be the union of the 2-cells of $X$ that  correspond
 to the edges  in ${c_1}_w$.  Then $U(c_1,w)$  is  a neighborhood of $w$ in $\supp(c_1)$. Since 
$c$ is homotopic to $c_1*c_2$, for a small neighborhood $U$ of $w$ in $X$  we have 
$$\overline{\supp(c)}\cap U
\subset (U(c_1,w)\cup U(c_2,w))\cap U$$  and  
$$\overline{(U(c_2,w)-U(c_1,w))\cup (U(c_1,w)-U(c_2,w))}\cap U\subset \overline{\supp(c)}\cap U.$$   
We observe that the open 2-cells in 
$U(c_1,w)\cap U(c_2,w)$ are  disjoint from $\supp(c)$:  Let $p$ be a point in some open 2-cell of  
$U(c_1,w)\cap U(c_2,w)$; by the above remark   $\log_p(c_1)$ and $\log_p(c_2)$ are homotopic to the circle
$\Link(X,p)$ with length $2\pi$; since $c$ is homotopic to $c_1*c_2$,  the path $\log_p(c)$ is 
either null-homotopic or homotopic to 
twice  of   a  generate in $\pi_1(\Link(X,p))$;  the fact that 
$\log_p(c)$ has length $<4\pi$ implies $\log_p(c)$ is 
   null-homotopic.   It  follows that 
$\overline{\supp(c)}\cap U=\overline{(U(c_2,w)-U(c_1,w))\cup (U(c_1,w)-U(c_2,w))}\cap U$. 
 Note $\overline{(U(c_2,w)-U(c_1,w))\cup (U(c_1,w)-U(c_2,w))}$
 is a  closed disk and it is the union of the 2-cells corresponding to the  edges in the geodesic
  $\log_w(x)b_0\cup b_0 \log_w(z)\subset \Link(X,w)$. We also observe that 
  $\log_w(x)b_0\cup b_0 \log_w(z)$ has length  $<2\pi$.

We next consider the case when $p$ is  not a corner. 
  We may assume  $p\in \interior(zw)$.    Then one of the following occurs: \newline
(1) $pw\cap px=pq $ is a nontrivial segment;\newline
(2)  $pw\cap px=\{p\} $ and  $px\cap yz$   contains   a  point  $r$;\newline
(3)   $pw\cap px=\{p\} $ and  $px\cap yz=\phi$.\newline

\noindent 
  Case (1):  In this case $zx=zq\cup qx$.
  Then  $T_1=xw\cup wq\cup qx$
  and $T_2=xy\cup yz\cup zx$    are  two  triangles, and $c$ is homotopic to $T_1*T_2$ with suitable orientations. 
     Note   $p\notin \supp(T_1)$ since 
$\log_p(T_1)$  is homotopic to  a path with  length $<\pi$.   It follows that there is a small neighborhood $U$ of $p$ in $X$ such
  that $U\cap \overline{\supp(c)}=U\cap \overline{\supp(T_2)}$. \newline

\noindent
Case (2): Let $T_1=(z,p,r)$, $T_2=(r,y,x)$ and $T_3=(p,x,w)$.
  Then $c$ is homotopic to $T_1*T_2*T_3$ with suitable orientations. 
 Note $p\notin \supp(T_2)$ since $r\in px$   implies the path $\log_p(T_2)$ is 
homotopic to  a path with
 length $<2\pi$. 
The segments $\log_p(r)\log_p(z)$ and $\log_p(x)\log_p(w)$ in $\Link(X,p)$ intersect in a 
(possibly degenerate) segment $\log_p(r)a$ ($a\in \Link(X,p)$).  The path $\log_p(z)a\cup a \log_p(w)$
  is a geodesic in $\Link(X,p)$. Th remark before Lemma \ref{s1} implies 
$U(T_1,p)$ and $U(T_3,p)$ are neighborhoods of $p$ in $\overline{\supp(T_1)}$ and
 $\overline{\supp(T_3)}$  respectively.   The argument in the second paragraph now implies 
 that the closed  disk $\overline{(U(T_1,p)-U(T_3,p))}\cup \overline{(U(T_3,p)-U(T_1,p))}$
 is a neighborhood of $p$ in $\overline{\supp(c)}$, and  that the disk is the union of the 2-cells that 
 correspond to the  edges  in $\log_p(z)a\cup a \log_p(w)$.

\noindent
Case (3):
In this case $T:=(p,x, w)$  and $Q:=(z, y,x,p)$ are  homeomorphic to a circle.  
 The path $c$ is homotopic to $T*Q$ with suitable orientations.
The discussion in the second paragraph implies 
 $\overline{\supp(Q)}$ determines a  geodesic segment $Q_p\subset \Link(X,p)$ with length $<2\pi$,
  and  $U(Q,p)$ is a neighborhood of $p$ in $\overline{\supp(Q)}$, where 
$U(Q,p)$   is the union of the  2-cells that  correspond to the edges in $Q_p$. 
 On the other hand, $U(T,p)$ is a neighborhood of $p$ in $\overline{\supp(T)}$,  and the link of 
 $U(T,p)$  at $p$ is the segment $\log_p(x)\log_p(w)$ which has length $<\pi$. 
 Since $\Link(X,p)$ is a finite $\CAT(1)$ graph  and the distance in $\Link(X,p)$
 between $\log_p(z)$  and $\log_p(w)$  is at least $\pi$,  
 the intersection  $Q_p\cap \log_p(x)\log_p(w)=\log_p(x)a$ is a 
(possibly degenerate) segment.  Let $U(c,p)$ be the union of the 
  2-cells that correspond to the edges in
$(Q_p-\log_p(x)a)\cup a\log_p(w)$.   Then $U(c,p)$  is  a  disk.  
The argument in the second paragraph  shows that 
 $U(c,p)$   
is a neighborhood of $p$ in $\overline{\supp(c)}$.

\end{proof}




The following lemma is 
a special case of a general observation due to B. Kleiner.

\b{Le}\label{extd}
{Let $c\subset  X$ be a triangle or  quadrilateral  homeomorphic to a circle,  and 
$x\in \supp(c)$.  Then  any   nontrivial geodesic segment $yx$ in $X$  can be extended into $\supp(c)$, 
 that is, there is some point $z\in X-\{x\}$ such that $xz\subset \supp(c)$ and $x\in yz$.} 

\end{Le}

\b{proof}  Let $\xi\in \Link(X,x)$  be  the initial  direction of  $yx$  at $x$.
  Since $x\in \supp(c)$,   $\log_x(c)$   is  homotopic to     a  nontrivial loop  $c_x$
     in the $\CAT(1)$
 graph   $\Link(X,x)$   and   $\interior(U(c,x))\subset \supp(c)$, where   $U(c,x)$ is the union of 
 2-cells corresponding to the edges in $c_x$. 
It follows that there exists  some  $\eta\in c_x$ such that the distance in $\Link(X,x)$ from $\xi$ to $\eta$ 
 is $\pi$.   Now we choose a geodesic segment $xz\subset \interior(U(c,x))$  such that the initial direction of 
 $xz$ at $x$ is $\eta$.   Then $yz=yx\cup xz$ is an extention of $yx$ into $\supp(c)$.

\end{proof}






\subsubsection{Gauss-Bonnet formula for piecewise Riemannian 2-complexes}  \label{gaussb}

In this section we recall the Gauss-Bonnet formula for 
finite 2-complexes with piecewise Riemannian  metrics  (see for example \cite{BB}).

Let $X$ be  a finite 2-complex.     We suppose  each  2-cell    of $X$  has a Riemannian 
metric  such that the 1-cells   are geodesics, and for any two 2-cells $f_1$, $f_2$ with 
 $f_1\cap f_2\not=\phi$, the induced metrics on $f_1\cap f_2$ from $f_1$  and  $f_2$
  coincide. 

For a vertex $v$  let
  $$\chi(v)=\chi(\Link(X,v)),$$
 the Euler characteristic of the link   $\Link(X,v)$.  For a 2-cell $f$   incident   to $v$,  denote by 
 $\alpha(v,f)$ the interior angle of $f$ at $v$.  The complete angle at $v$ is
$$\alpha(v)=\Sigma \alpha(v,f),$$
 where $f$ varies over all 2-cells  incident  to $v$.   The  curvature measure of $v$ is then by definition 
$$\kappa(v)=(2-\chi(v))\pi-\alpha(v).$$

For a 2-cell  $f$ denote by $K$ the Gaussian curvature of $f$. Then the curvature of $f$ by
 definition is  
$$\kappa(f)=\int_f K.$$
Now the total curvature of $X$ is defined by:
$$\kappa(X)=\Sigma_s \kappa(s)$$
 where $s$ varies over  all the 0-cells   and 2-cells of $X$ (the curvatures of the 1-cells are 0 since they are geodesics).

The Gauss-Bonnet formula (see \cite{BB})  says
 $$\kappa(X)=2\pi \chi(X).$$
Note $\chi(X)=1$  when $X$ is  homeomorphic to  a disk.  If  a 2-cell $f$ has  constant 
Gaussian curvature  $-1$, then $\kappa(f)=-\Area(f)$.

Let $\De$  be  a   Fuchsian building.   For  a  finite sequence of pairwise 
 distinct points   $P=<x_1, \cdots, x_l>$ ($l\ge 3$)  
 in $\De$,  define $d(P)=(l-2)\pi-\Sigma_{i=1}^l \angle_{x_i}(x_{i-1}, x_{i+1})$ ($i\mod l$). 
  If $T=(x,y,z)$  is  a   triangle,  
     set   $d(T)=d(<x,y,z>)$. Similarly  we  define  $d(Q)$ 
    if  $Q$  is  a  quadrilateral.  
For  a      triangle or quadrilateral $c\subset \De^{(1)}$
homeomorphic  to  a 
 circle,   let  $n(c)$ be  the number of chambers 
in  $\overline{\supp(c)}$.   
  Recall  $\overline{\supp(c)}$  is  a  finite subcomplex.    
 We also denote by $A_0$ the area of a chamber.

\b{Cor}\label{ctq}
{Let $\De$ be a Fuchsian building, and $c\subset \De^{(1)}$ a triangle or quadrilateral homeomorphic  to  a 
 circle.   If  $\overline{\supp(c)}$  is  homeomorphic  to a  closed disk,
then  $d(c)\ge n(c)A_0$.}

\end{Cor}

\b{proof}
Let $X=\overline{\supp(c)}$.  
Since $\De$ is  $\CAT(-1)$,  $\kappa(v)\le 0$   if  $v\in X$ is a vertex  but not  a corner of $c$. 
  On the other hand, if $v$  is  a  corner of $c$, then the length of 
$\Link(X,v)$  is larger than or equal to the angle of $c$ at $v$. 
The corollary follows from these two observations and Gauss-Bonnet formula. 

\end{proof}

Let us make some observation about  $A_0$.   The Gauss-Bonnet formula 
  applied  to the chamber $R$  implies:
 $A_0=(k-2)\pi-\Sigma(\p R)$, where $R$ is  a $k$-gon  and $\Sigma(\p R)$  is the sum of angles at the vertices of $R$.

\noindent 
$\bullet$  if $k\ge 5$, then $A_0\ge \pi/2$, with equality precisely when $k=5$ and all the angles of $R$ equal to $\pi/2$;\newline
$\bullet$  if $k=4$, then $A_0\ge \pi/6$, with equality precisely when  
$R$ has angles $\frac{\pi}{2}, \frac{\pi}{2}, \frac{\pi}{2},  \frac{\pi}{3}$;\newline
$\bullet$  if $k=3$ and $R$ has no right angle, then $A_0\ge \pi/12$,  with equality precisely when  $R=(3,3,4)$;\newline
$\bullet$  if $R=(2,8,8)$, then $A_0=\pi/4$;\newline
$\bullet$  if $R=(2,6,8)$, then $A_0=5\pi/24$;\newline
$\bullet$  if $R=(2,6,6)$, then $A_0=\pi/6$;\newline
$\bullet$  if $R=(2,4,6)$, then $A_0=\pi/12$;\newline
$\bullet$  if $R=(2,4,8)$, then $A_0=\pi/8$;\newline
$\bullet$  if $R=(2,3,8)$, then $A_0=\pi/24$.

\subsubsection{Triangles  in the 1-skeleton} \label{tqdisk}

Let $\De$  be   a Fuchsian   building and $c\subset \De^{(1)}$  
a triangle or quadrilateral homeomorphic  to  a 
 circle.  Then $X:=\overline{\supp(c)}$  is  a finite subcomplex homeomorphic to a compact surface with boundary.
  Notice  a vertex $v\in X$ is   a   special point    if and  only  if    $\kappa(v)<0$.
The remark before Lemma  \ref{s1}
   implies that  the corners of a triangle are not special points.  
Proposition \ref{tri1} follows from the following result:

\b{Prop}\label{trcon}
{Let $\De$  be a Fuchsian building and $T\subset \De^{(1)}$ a triangle that is 
homeomorphic to a   circle.  Then $\overline{\supp(T)}$  contains  no special points.
 In particular,  $\overline{\supp(T)}$  is  homeomorphic  to  a  closed disk and is convex in $\De$.}

\end{Prop}

The second claim in the proposition follows easily from the first one:  assume 
$X=\overline{\supp(T)}$  contains  no special points;  then $X$ is locally convex in 
 the  $\CAT(-1)$ space $\De$; it follows that $X$  is actually convex in $\De$ and 
 hence is  contractible.

\b{remark}\label{trquds}
{I suspect that any triangle or quadrilateral homeomorphic to a circle in a  $\CAT(0)$ 2-complex 
 must bound a disk, at least for those $\CAT(0)$ 2-complexes that admit cocompact groups of isometries.}

\end{remark}

Let $T=(x,y,z)\subset \De^{(1)}$  be a triangle   homeomorphic  to  a  circle. 
 Recall  the angles of the chamber $R$ lie in $\{\pi/2, \pi/3, \pi/4, \pi/6, \pi/8\}$. In particular,
  they   are integral multiples of $\pi/24$.  It follows that $d(T)$ is an integral multiple of
  $\pi/24$.  
On the other hand, $d(T)>0$ because $\De$ is a $\CAT(-1)$ space. 
Therefore $\pi/24\le d(T)\le   5\pi/8$.  
We will prove Proposition \ref{trcon}  by inducting on $d(T)$, starting with triangles with 
 $d(T)=\pi/24$.    
We first make an observation that shall be used often   later.

The  following lemma follows easily from the  triangle inequality and the fact that $d(P)>0$
  for any finite sequence of pairwise distinct points  $P=<x_1, \cdots, x_l>$ ($l\ge 3$).

\b{Le}\label{divi}
{Let $P=<x_1, \cdots, x_l>$ \e{($l\ge 3$)}  be  a sequence of pairwise distinct points,
  and $y_i\in x_ix_{i+1}$,   $y_{j}\in  x_jx_{j+1}$  with  $i< j$.    
  If   $P'=<y_i, x_{i+1},  \cdots, x_j, y_{j}>$ is a sequence of 
pairwise distinct points,  then  $d(P')<d(P)$.} 
 
\end{Le}

 We   first study triangles with $d(T)=\pi/24$.

\b{Le}\label{tr23}
{Let $\De$  be a Fuchsian building and $T=(x,y,z)\subset \De^{(1)}$ a triangle that is 
homeomorphic to a   circle.  If $d(T)=\pi/24$, then $R=(2,3,8)$  and 
  $T$ is the boundary of   some chamber.}

\end{Le}

\b{proof}
We claim that none of the segments $xy$, $yz$ and $zx$ contains any vertex in the interior. 
Suppose, say, $xy$ contains   a vertex $p$ in the interior. Then there is an edge
 $pq\subset \overline{\supp(T)}$ such that $\interior(pq)\subset \supp(T)$. 
Since 
$\overline{\supp(T)}$  is a compact surface,
 Lemma \ref{extd} implies  $pq$ extends into $\supp(T)$ 
  and  the extension     hits 
 the boundary of $\overline{\supp(T)}$ (which is $T$) at some point $r$.  
 The segment $pr$ lies in $\De^{(1)}$.    Let $c'$ be the union of $pr$ with  one 
   of the two components of $T-\{p,r\}$.
  Then $c'$     is  either a triangle or a quadrilateral. Lemma \ref{divi} 
 implies  $d(c')< d(T)=\pi/24$, which is a contradiction since $d(c')>0$ is an 
integral multiple of $\pi/24$.
Hence $xy$, $yz$ and $zx$ are three edges in $\De$. Similarly one shows that 
 $\log_x(y)\log_x(z)$ is an edge in $\Link(\De, x)$. 
It follows that  $xy$ and $xz$ lie in the boundary of some chamber   $C$ and 
 we conclude that $C$ is actually a triangle and $T$ is its boundary. 
Since $d(T)=\pi/24$,  $A_0=\pi/24$. The observation about $A_0$ shows $R=(2,3,8)$.

\end{proof}

For a  triangle or quadrilateral $c\subset \De^{(1)}$  homeomorphic  to a circle  and 
 $x\in \overline{\supp(c)}$, let $c_x$  be the link $\Link(\overline{\supp(c)}, x)$.  
  By Propositions \ref{trisupp}
 and \ref{quipp}
$c_x$ is a circle if $x\in \supp(c)$ and is a segment if $x\in c$.

\b{Le}\label{bounosi}
{Let $T=(x,y,z)\subset \De^{(1)}$ be a  triangle 
  homeomorphic to a  circle  
with $d(T)=\frac{k\pi}{24}$,  $2\le k\le 15$.
 Assume   $\overline{\supp(T')}$  has   no special points for  every  
triangle $T'\subset \De^{(1)}$   homeomorphic  to a  circle with $d(T')\le \frac{(k-1)\pi}{24}$.
Then $T$ contains  no special points  of  $\overline{\supp(T)}$.}
 
\end{Le}

\b{proof}
Suppose the lemma  is false.  
Since the corners of $T$ are not special points, 
we may assume there is a special  point $p\in \interior(xy)$.  
  Then $c_p$ is  a segment with length $>\pi$. Note $c_p$ is an edge path in $\Link(\De, p)$.
   Let    $\eta\in c_p$   be the point  such that the subsegment of $c_p$ from $\log_p(x)$ to $\eta$ has length 
 $\pi$.    Then $\eta$ is a  vertex in $\Link(\De, p)$.  There is an edge   $pq\subset \overline{\supp(T)}$
   with initial direction $\eta$ and $\interior(pq)\subset\supp(T)$.   We extend  the geodesic $pq$  into 
  $\supp(T)$ until it hits a point $r\in T$.  The uniqueness of geodesic implies $r\in \interior(yz)$.
 Let $T_1=(x,r,z)$ and $T_2=(p,y,r)$.  
  Lemma \ref{divi} implies   $d(T_i)\le \frac{(k-1)\pi}{24}$  for $i=1,2$.
The   assumption   then implies that $\overline{\supp(T_i)}$  contains  no special points.
In particular $\overline{\supp(T_i)}$  is homeomorphic  to   a  closed  disk and we can apply 
Corollary \ref{ctq}
to  $T_i$.

Note $T_2$ has an even angle at $p$  and $A_0\ge \pi/24$. 
Since   $\angle_r(p, z)+\angle_r(p,y)\ge \pi$,  Corollary \ref{ctq}
applied to $T_2$ shows that $\angle_r(p,z)\ge \angle_p(r,y)+\angle_y(x,z)+  2 A_0\ge 
 2\times \frac{\pi}{m(p)}+\pi/8+2\times \pi/24$.  Since there are at least $m(p)$ chambers in 
$\overline{\supp(T_1)}$  incident to $p$, Corollary \ref{ctq}
applied to $T_1$   implies that  $\pi\ge \pi/8+\pi/8+\angle_r(p,z)+m(p)A_0$. 
Combining these two inequalities we obtain  $m(p)A_0+\frac{2\pi}{m(p)}\le \frac{13\pi}{24}$. 
  Since $A_0$ is an integral multiple of $\pi/24$ and $m(p)\in \{2,3, 4, 6, 8\}$,
 this inequality  never holds and we have a contradiction.

\end{proof}

\b{Le}\label{intounosi}
{Let $T=(x,y,z)\subset \De^{(1)}$ be a  triangle 
  homeomorphic to a  circle  
with $d(T)=\frac{k\pi}{24}$,  $2\le k\le 15$.
 Assume   $\overline{\supp(T')}$  has  no special points for  every  
triangle $T'\subset \De^{(1)}$   homeomorphic  to a  circle with $d(T')\le \frac{(k-1)\pi}{24}$.
 Then  there is  no special point in  $\supp(T)$.}

\end{Le}

\b{proof}
Suppose   there is a special  point  $p\in \supp(T)$.   Then $c_p$ has length $>2\pi$. 
 There is at least one edge $pq\subset \overline{\supp(T)}$  such that its extension in 
$\supp(T)$ eventually hits $T$ at  a  point  $p_1$  different from $x$, $y$ and $z$.
 We may  assume $p_1\in \interior(xy)$.   Denote $\eta_1=\log_p(p_1)$. 
 Let $\eta_2\not= \eta_3\in c_p$ be the two points such that their distance  in $c_p$ from 
$\eta_1$ is exactly $\pi$.  We extend $p_1p$ inside $\supp(T)$ in two ways, one in the direction
 $\eta_2$, the other in the direction $\eta_3$.  Suppose they eventually hit $T$ at 
 $p_2$ and $p_3$ respectively.  Note $p_2, p_3\in T-xy$.

We claim one of $p_2,  p_3$ lies in $xz$ and the other lies in $yz$. Suppose otherwise, say, $p_2, p_3\in yz$
  and $p_3\in \interior(p_2y)$. 
Let $T_1=(p_1, p_2, y)$.   Lemma \ref{divi} and the   assumption 
imply that  $\overline{\supp(T_1)}$  contains no special points. 
This contradicts to the facts that $p\in p_1p_2, p\in  p_1p_3$ and   that  the
 initial directions of $pp_2$ and $pp_3$ at $p$ are different.

We may assume $p_2\in xz$ and $p_3\in yz$. Let $T_2=(p_1, p_2, x)$,  $T_3=(p_1, p_3, y)$ and
 $Q=(p, p_2, z, p_3)$.  Also denote $m=m(p)$. Lemma \ref{divi} and the assumption imply that 
$\overline{\supp(T_i)}$ ($i=2,3$)  is homeomorphic to a closed disk and we can apply 
 Corollary \ref{ctq}  to $T_i$.  Since  $p\in \interior(p_1p_i)$,  
there are at least $m$ chambers in $\overline{\supp(T_i)}$.
Now Corollary \ref{ctq}  implies $\angle_{p_1}(p_2, x)+\angle_{p_2}(p_1, x)+\angle_x(y,z)+mA_0\le \pi$,
$\angle_{p_1}(p_3, y)+\angle_{p_3}(p_1, y)+\angle_y(x,z)+mA_0\le \pi$.
Since $\De$ is $\CAT(-1)$ and all the angles of $Q$ are integral multiples of $\pi/24$, we have
$\angle_{p_2}(p, z)+\angle_{p}(p_2, p_3)+\angle_{p_3}(p, z)+\angle_z(p_2, p_3)\le 2\pi-\pi/24$.
Notice $\angle_{p_2}(p,z)+\angle_{p_2}(p,x)\ge \pi$, 
$\angle_{p_3}(p,z)+\angle_{p_3}(p,y)\ge \pi$,
$\angle_{p_1}(p,y)+\angle_{p_1}(p,x)\ge \pi$.
Since the angles of $T$ are $\ge \pi/8$, the above inequalities imply 
$2\pi/m+2mA_0\le 7\pi/12$.  The second term here, $2mA_0\ge 4\times \pi/24=\pi/6$. It follows that
$2\pi/m\le 5\pi/12$, which implies   $m\ge 6$.  The first term $2\pi/m\ge \pi/4$.  It follows that $\pi/3\ge 2mA_0\ge 12A_0$,
or, $A_0\le \pi/36$,   contradicting to the fact that $A_0\ge \pi/24$.

\end{proof}

  The proof of  Proposition  \ref{trcon}  is  now complete. 
Next we     consider  Proposition \ref{tri2}.

\b{Le}\label{rnottr}
{Let $\De$ be a Fuchsian building with chamber $R$. 
 If $R$ is not a triangle, then  no triangle  in $\De^{(1)}$  is homeomorphic  to a  circle.}

\end{Le}

\b{proof}
Suppose $T\subset \De^{(1)}$  is  a triangle homeomorphic  to a circle.
Assume $R$ is a $k$-gon with $k\ge 4$ and let $\alpha_0$ be the smallest angle of $R$. 
Proposition \ref{trcon}  and Corollary \ref{ctq}  imply that $3\alpha_0+n(T)A_0\le \pi$.
 Since $T$ is  a triangle and $R$ is not, $n(T)\ge 2$. 
  Notice  $A_0\ge \pi/6$  for $k\ge 4$. It follows that $\alpha_0\le 2\pi/9$ and so 
    $\alpha_0=\pi/6$ or $\pi/8$.   By Gauss-Bonnet $A_0\ge \pi/3$.   Repeating   the above argument 
 we obtain $\alpha_0\le \pi/9$,  contradicting  to the 
 fact that  $\alpha_0\ge \pi/8$.

\end{proof}

There are two ways to prove the other claims in Proposition \ref{tri2}.  The first method   is to use
Corollary \ref{ctq}  and the observation about $A_0$. Basically Corollary \ref{ctq}  puts severe 
 restriction on the possible configurations for $\overline{\supp(T)}$
and one can list all the triangles in $\De^{(1)}$ that are homeomorphic to a circle.
  We omit  the details here. 
 On the other hand, 
Proposition \ref{trcon}  implies that $\overline{\supp(T)}$  is isomorphic to a convex  subcomplex 
in an apartment. 
 So the second method is 
to  find   all the triangles by examining   the tessellation of $\H^2$ by the chamber.
 Here we exhibit the  tessellations by $(2,6,8)$, $(2,4,8)$ and $(2,3,8)$.



\subsubsection{Quadrilaterals  must bound disks}  \label{nonexis}

In this section we  show that $\overline{\supp(Q)}$  is homeomorphic to a   closed 
  disk   for any   quadrilateral   $Q\subset \De^{(1)}$
homeomorphic  to  a  circle.   
Proposition \ref{quar11}  then follows  from   Corollary  \ref{ctq}.
The proof is based on the following lemma.

\b{Le}\label{lnosin}
{Let $Q=(x,y,z,w)\subset \De^{(1)}$ be   a quadrilateral that  is homeomorphic  to  a  circle. 
   Then $\overline{\supp(Q)}$  is homeomorphic to a   closed 
  disk
   if the following conditions are satisfied:\newline
\e{(1)} there is no special  point in $\supp(Q)$;\newline
\e{(2)}   $x, y, z, w$ are not   special   points;\newline
\e{(3)} for each $p\in Q$, the path $Q_p$ has length $\le 2\pi$;\newline
\e{(4)} there is at most one   special point in the interior  of each  side  of  $Q$.}

\end{Le}

\b{proof} 
  By Proposition \ref{quipp},  $\overline{\supp(Q)}$  is a compact surface with boundary
 $Q$.   It suffices to prove that  $\overline{\supp(Q)}$  is simply connected. 
Notice $\overline{\supp(Q)}$ has nonpositive curvature with respect to the path metric, which is denoted by   $d'$.
Suppose  $\overline{\supp(Q)}$ is  not  simply connected.  Then  
   there is  a simple closed geodesic (with respect to $d'$) $\gamma\subset \overline{\supp(Q)}$ (see p. 202 of \cite{BH}).  
  $\gamma$ must contain  special  points since otherwise $\gamma$ is a closed geodesic in the 
$\CAT(-1)$ space $\De$.  
One also sees that $\gamma$ contains at least three  special   points.
 The assumption implies   that  there are either 3 or  4 special  points on  $\gamma$,
  and  hence   $\gamma$  is  actually  a  triangle or quadrilateral in $\De$  with the special points as corners.
  For any special point $p\in \gamma\cap Q$, condition (3) implies that the angle (measured in $\De$)
that $\gamma$ makes at $p$ is at least $\pi/2$.   
Now we have  a contradiction since the angles at the corners 
 of the triangle or quadrilateral $\gamma$ 
are at least $\pi/2$ and $\De$  is a $\CAT(-1)$ space.

\end{proof}

We first show condition (3) of Lemma \ref{lnosin} holds:

\b{Le}\label{qststep1.1}
{Let $\De$  be a Fuchsian building and $Q=(x_1,x_2,x_3,x_4)\subset \De^{(1)}$ a   quadrilateral that is 
homeomorphic to a   circle.  Then $\length(Q_p)\le 2\pi$  for every $p\in Q$.}

\end{Le}

\b{proof}
Suppose there is some $p\in Q$ with $\length(Q_p)> 2\pi$.  Then $p$ is  vertex. 
The proof of Lemma \ref{quipp} shows that 
  $p$  is not a corner. We may assume $p\in interior(x_1x_2)$.  We claim $m(p)\not=2$.
Suppose $m(p)=2$. Then $\length{Q_p}\ge  3\pi$  and there are two edges $pq_1,  pq_2\subset \overline{\supp(Q)}$ such that 
 the angle between any two of $px_1, pq_1, pq_2$ and $px_2$ is $\pi$.  We extend $pq_1, pq_2$ inside 
$\supp(Q)$ until they hit $Q$ at $r_1$ and $r_2$ respectively. The uniqueness of geodesic implies $r_1, r_2\in x_3x_4$. 
 But then $r_1r_2=r_1p\cup pr_2$ is not contained in $x_3x_4$, contradicting to  the uniqueness of geodesic.
Hence $m(p)\ge 3$.   Since $\length{Q_p}> 2\pi$, there are edges $pq_1,  \cdots,  pq_l\subset \overline{\supp(Q)}$,
 with $l\ge 6$ and  $\interior(pq_i)\subset \supp(Q)$.  We extend each $pq_i$  until it  hits $Q$  at  
  some  point $r_i$.  We orient $x_1x_4\cup x_4x_3\cup x_3x_2$ from $x_1$ to $x_2$  and
  relabel the points $q_1, \cdots, q_l$ such that $r_1,  \cdots, r_l$ are in linear order in $x_1x_4\cup x_4x_3\cup x_3x_2$
  when one travels from $x_1$ to $x_2$.  We set $r_0=x_1$ and $r_{l+1}=x_2$. 
Since $l\ge 6$, there exist at least 3 points $r_i, r_{i+1},  r_{i+2}$ that lie on the same side of $Q$. 
 In particular, the triangle $(p, r_i, r_{i+2})$ has a side $r_ir_{i+2}$ containing at least two edges. 
 Proposition \ref{tri2} implies $R$ is a right triangle.

We claim $R=(2,3,8)$. Suppose $R\not=(2,3,8)$.  Then $m(p)\ge 4$ and $l\ge 8$.  There are 
 4 points $r_i, r_{i+1}, r_{i+2},  r_{i+3}$  that lie on the same side of $Q$. 
Proposition \ref{tri2} applied to $(p, r_i, r_{i+1})$,   $(p,  r_{i+1}, r_{i+2})$
  and  $(p, r_{i+1},  r_{i+3})$   shows that $\angle_{r_{i+1}}(r_i,  p)$, $\angle_{r_{i+1}}(r_{i+2},  p)$,
  $ \angle_{r_{i+2}}(r_{i+1},  p)$,  $\angle_{r_{i+2}}(r_{i+3},  p)\le \pi/2$. The triangle inequality 
   then implies all these 4 angles are equal to $\pi/2$.  In particular, $(p, r_{i+1}, r_{i+2})$ has two right angles, 
which is impossible.  Hence $R=(2,3,8)$ and $m(p)=3$ or 8.

Suppose $m(p)=8$. Since the initial direction of $px_1$ and $px_2$ are of the same type,  $l\ge 17$. 
 There are 7 points $r_i, \cdots, r_{i+6}$ that lie on the same side of $Q$. The argument in the preceding paragraph shows 
 at least one of the 4 angles $\angle_{r_{i+2}}(p, r_{i+1}),$  $\angle_ {r_{i+2}}(p, r_{i+3}),$ 
$\angle_{r_{i+3}}(p, r_{i+2}),$  $\angle_ {r_{i+3}}(p, r_{i+4}),$  is $>\pi/2$.  
Assume, say, $\angle_{r_{i+2}}(p, r_{i+1})>\pi/2$. Let $T_1=(p, r_i, r_{i+1})$ and $T_2=(p, r_{i+1},  r_{i+2})$.
  Proposition \ref{tri2}  applied to $T_2$ shows $\angle_{r_{i+1}}(p, r_{i+2})=\pi/8$.
  It follows that $\angle_{r_{i+1}}(p, r_i)=7\pi/8$  and $d(T_1)<0$, which is impossible.

Finally assume $m(p)=3$.  Then $l\ge 6$.  Notice Proposition \ref{tri2}
 implies if $(y_1, y_2, y_3)\subset \De^{(1)}$ is homeomorphic to a circle and $m(y_1)=3$, then 
 $\angle_{y_2}(y_1, y_3)\le \pi/2$ and equality holds only when $(y_1, y_2, y_3)$ is the boundary of  a chamber.
  It follows that if $r_i, r_{i+1}, r_{i+2}$ lie on the same side of $Q$, then $m(r_{i+1})=2$ and 
 $(p, r_i, r_{i+1})$, $(p,  r_{i+1},  r_{i+2})$  are boundaries of chambers.
This implies that each side of $Q$  contains at most 3 of the $r_i$'s.  We claim there is no $i$ such that
 $\{r_i, r_{i+1}, r_{i+2}\}\subset x_1x_4-\{x_4\}$ and $\{r_{i+3},  r_{i+4}, r_{i+5}\}\subset x_3x_4-\{x_4\}$;
 and similar claim also holds about the two sides of $Q$ incident to $x_3$.  
 Suppose there is  such an $i$.  Let $Q_1=(p, r_{i+2}, x_4, r_{i+3})$.  Then the observation about angles implies that 
  $d(Q_1)<0$, which is impossible.  
Recall there are at least 8 $r_j$'s. 
If $x_3=r_i$ for some $i$,  then $x_2x_3$ contains exactly 3 
 $r_j$'s,   and each of $x_1x_4-\{x_4\}$,  $x_3x_4-\{x_4\}$  contains exactly 3 $r_j$'s,  contradicting to the above claim. 
Hence $x_3$ is distinct from all the $r_j$'s. Similarly $x_4$ is also distinct from all the $r_j$'s. 
The above claim implies  $x_3x_4$ contains exactly two $r_j$'s and they lie in $\interior(x_3x_4)$, and 
 each of $x_1x_4-\{x_4\}$,  $x_2x_3-\{x_3\}$ also contains exactly 3 $r_j$'s. So we have 
$r_1, r_2\in \interior(x_1x_4)$,  $r_3, r_4\in \interior(x_3x_4)$   and $r_5, r_6\in \interior(x_2x_3)$. 
 Consider $T=(p, r_3, r_4)$.  Since $m(p)=3$,  Proposition \ref{tri2} implies  at least one of 
$\angle_{r_3}(p, r_4)$,  $\angle_{r_4}(p, r_3)$ is $\le \pi/4$. We may assume $\angle_{r_4}(p, r_3)\le \pi/4$.
It follows that $\angle_{r_4}(p, x_3)\ge 3\pi/4$ and $d(Q_2)<0$, which is impossible, where
 $Q_2=(p, r_4, x_3, r_5)$.   The contradiction proves the lemma.

\end{proof}

The proof of Proposition \ref{quar11}
 is similar to that of Proposition \ref{trcon}:
 we induct on $d(Q)=\frac{k\pi}{24}$, starting with $k=1,2$.

\b{Le}\label{qststep}
{Let $\De$  be a Fuchsian building and $Q=(x_1,x_2,x_3,x_4)\subset \De^{(1)}$ a   quadrilateral that is 
homeomorphic to a   circle.  Then $d(Q)\ge \frac{2\pi}{24}$.  Moreover, if 
 $d(Q)=\frac{2\pi}{24}$, then  $R=(2,3,8)$ and $\overline{\supp(Q)}$ can be obtained by gluing two chambers along an edge.}

\end{Le}

\b{proof}
Note $d(Q)>0$  is  an integral multiple of $\pi/24$.  Hence $d(Q)\ge \frac{\pi}{24}$.
Suppose $d(Q)= \frac{\pi}{24}$.    Then the proof of Lemma \ref{tr23} shows that 
$x_ix_{i+1}$ ($i \mod 4$) is an  edge  and $\log_{x_i}(x_{i-1})\log_{x_i}(x_{i+1})$
    is an edge in $\Link(\De, x_i)$.  From this it is not  hard to see that 
 $Q$ is the boundary of a chamber. In particular the chamber is a   4-gon.  
In this case $d(Q)=A_0\ge \pi/6$,
 contradicting to the assumption that $d(Q)= \frac{\pi}{24}$.
Hence $d(Q)\ge \frac{2\pi}{24}$.

Now suppose $d(Q)=\frac{2\pi}{24}$.   
Assume  some side of $Q$, say $x_1x_2$, contains a vertex  $p$ in the interior.  
For an edge $pp'\subset \overline{\supp(Q)}$ with $\interior(pp')\subset \supp(Q)$, we extend 
  $pp'$ inside $\supp(Q)$ until it hits $Q$ at some point $q$.   Since $d(Q')\ge \frac{2\pi}{24}$
 for  all quadrilaterals, Lemma \ref{divi} implies $q\in \interior(x_1x_4)$, or $q\in \interior(x_2x_3)$.
We may assume $q\in \interior(x_1x_4)$.  Denote  $T=(p,q, x_1)$ and  $P=<p, x_2, x_3, x_4, q>$. 
 Lemma \ref{divi} implies $d(T)=d(P)=\pi/24$. By Lemma \ref{tr23},
 $R=(2,3,8)$ and $T$ is the boundary of some chamber. At least one of $m(p)$, $m(q)$  is  $\not=2$.
 We may assume $m(p)\not=2$.  Then there is an  edge 
$pp''\subset \overline{\supp(Q)}$ different from $pp'$ 
with $\interior(pp'')\subset \supp(Q)$. We  extend $pp''$    inside $\supp(Q)$ until 
 it hits $Q$ at some 
 point $r\not=q$. Then $r\in \interior(x_1x_4)$, or $r\in \interior(x_2x_3)$.
 Note $r\in \interior(x_1x_4)$ cannot hold since otherwise both $x_1q$ and $x_1r$ are edges. 
 Now Lemma \ref{divi} applied to $P$ and $T'=(p,r, x_2)$ implies  $d(T')<\pi/24$, a  contradiction.
 Hence   each $x_ix_{i+1}$ is an edge.

The first paragraph shows there is some $i$ such that $\log_{x_i}(x_{i-1})\log_{x_i}(x_{i+1})$
 contains at least two edges. It implies that there is some 
edge $x_ip'\subset \overline{\supp(Q)}$ with $\interior(x_ip')\subset \supp(Q)$. 
We extend $x_ip'$ inside $\supp(Q)$ until it hits some $q\in Q$. 
Since $d(Q)=2\pi/24$  and $d(Q')\ge 2\pi/24$  for all quadrilaterals $Q'$, Lemma \ref{divi} implies 
 $q=x_{i+2}$.  Let $T_1=(x_i, x_{i+1}, x_{i+2})$ and $T_2=(x_i, x_{i-1}, x_{i+2})$.  
  By   Lemma \ref{divi},        $d(T_1)=d(T_2)=\pi/24$.
  Lemma \ref{tr23} then implies that  $R=(2,3,8)$  and 
both $T_1$ and $T_2$ are the boundaries of chambers.

\end{proof}

We first consider the case when   one  of   the corners of   $Q$  is a special point. 
 Recall for each $x\in \overline{\supp(Q)}$, $Q_x$ denotes the link of $\overline{\supp(Q)}$ at $x$.

\begin{figure}[h]
\centering\epsfig{file=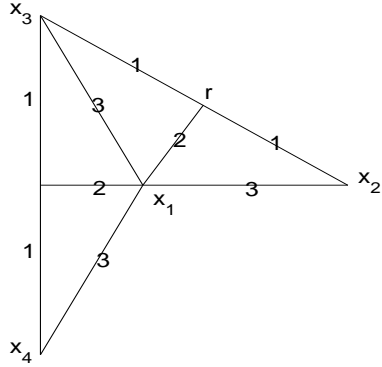,  height=5cm,  width=5cm}\caption{Corner $x_1$  is  a special point}\label{10l18}
\end{figure}

\b{Le}\label{versing}
{Let $\De$  be a Fuchsian building and $Q=(x_1,x_2,x_3,x_4)\subset \De^{(1)}$ a   quadrilateral that is 
homeomorphic to a   circle.  If $Q_{x_1}$ has length $>\pi$, then $R=(2,3,8)$ and $\overline{\supp(Q)}$ 
  must be as shown in   Figure \ref{10l18}.}

\end{Le}

\b{proof}
Let $\eta\in Q_{x_1}$ be the point such that the subsegment of $Q_{x_1}$ 
from $\log_{x_1}(x_4)$ to $\eta$  has length $\pi$. Let $x_1q\subset \overline{\supp(Q)}$ be the edge 
 with initial direction $\eta$ at $x_1$. We extend $x_1q$ inside  $\supp(Q)$   
   until it hits $Q$ at some point $r$. 
 The uniqueness of geodesic implies $r\in  \interior(x_2x_3)$.
 Let $T_1=(x_4, r, x_3)$ and $T_2=(x_1, x_2, r)$. 
Since     the  side  $x_4r=x_4x_1\cup x_1r$ of  $T_1$  contains  at least two edges,    
Proposition \ref{tri2} implies $R$ is a right  triangle  and 
 $\angle_r(x_1, x_3)\le  \pi/2$.  
We  claim $\angle_r(x_1, x_2)\le \pi/2$.  Suppose $\angle_r(x_1, x_2)> \pi/2$.
 Then  Proposition \ref{tri2}  implies that $R=(2,3,8)$  and $\angle_r(x_1, x_2)= 2\pi/3$,
   $\angle_{x_1}(r,  x_2)=\pi/8$. Now consider $T_1$:  the side $x_4r$ contains the  vertex $x_1$ 
in the interior  with $m(x_1)=8$ 
 and  one endpoint $r$ of $x_4r$ satisfies $m(r)=3$.  However, by Proposition \ref{tri2}
  no such triangle  exists.

Now  $\angle_r(x_1, x_3)\le  \pi/2$,    $\angle_r(x_1, x_2)\le \pi/2$  and 
 the triangle  inequality imply that 
$\angle_r(x_1, x_2)=\angle_r(x_1, x_3)=\pi/2$.
Consider $T_2$. Since $\angle_r(x_1, x_2)=\pi/2$,  $\angle_{x_1}(r, x_2)<\pi/2$.
   In particular,  $m(x_1)>2$.  Now consider $T_1$:  the angle at  $r$  is  $\pi/2$, the side $x_4r$ contains  the
  vertex  $x_1$ in the interior with $m(x_1)>2$.  Proposition \ref{tri2} implies $R=(2, 3, 8)$,
    $m(r)=2$, $m(x_1)=3$ and $\overline{\supp(T_1)}$  is as shown in    Figure \ref{238}(d).
It follows that $\angle_{x_1}(r, x_2)=\pi/3$ and by Proposition \ref{tri2} $T_2$ is the boundary of  a chamber. 
 Now the lemma follows.

\end{proof}

From now on we  assume that none of the 4   corners of $Q$ is a special point.

\b{Le}\label{twosigonsside}
{Let $Q=(x_1, x_2, x_3, x_4)\subset \De^{(1)}$ be a  quadrilateral that is homeomorphic to 
  a circle with $d(Q)=\frac{k\pi}{24}$, $k\ge 3$.  
Suppose $\overline{\supp(Q')}$ is  homeomorphic to  a closed disk for every 
 quadrilateral   $Q'\subset  \De^{(1)}$ homeomorphic to a  circle
 with $d(Q')\le \frac{(k-1)\pi}{24}$.     
  Then    for  each $i$,  $\interior(x_ix_{i+1})$ \e{(}$i\mod 4$\e{)}
contains at most one special point. }

\end{Le}

\b{proof}
Suppose the lemma is not true.  We may assume $\interior(x_1x_2)$ contains two 
 special points $y_1, y_2$, and $y_1\in x_1y_2$.  Let $\eta_i\in Q_{y_i}$  ($i=1,2$) 
 be the point in $Q_{y_i}$ such that the subsegment from $\log_{y_i}(x_{3-i})$  to  $\eta_i$ has length $\pi$.
  We extend  the geodesic $x_{3-i}y_i$ into $\supp(Q)$ in the direction of $\eta_i$  
until it hits a point $z_i$ on $Q$. 
 Note $y_iz_i\subset \De^{(1)}$.  The uniqueness of geodesic implies 
$z_1\in x_1x_4\cup x_4x_3$  and  
$z_2\in x_2x_3\cup x_3x_4$.  Note $z_1z_2=z_1y_1\cup  y_1y_2\cup  y_2z_2$.  
The uniqueness of geodesic implies  at most one of $z_1$, $z_2$ lies  on $x_3x_4$. We claim none of 
$z_1$, $z_2$ lies  on $x_3x_4$. Suppose the contrary holds, say, $z_2\in x_3x_4$. 
Then $z_1\in \interior(x_1x_4)$. Consider $T=(x_1,  z_2, x_4)$:  $z_1\in \interior(x_1x_4)$,
   $y_2\in \interior(x_1z_2)$ and $y_2z_1\cap   x_1z_2=y_2y_1$ is a nontrivial segment.
 This contradicts to Proposition  \ref{tri1}.    Hence  $z_1\in x_1x_4$ and $z_2\in   x_2x_3$.
Notice $\angle_{y_i}(x_i, z_i)$ ($i=1,2$) is an even angle.
Proposition \ref{tri2}   applied to  $T_i=(x_i, y_i,  z_i)$ 
    implies  $R$ is  a right triangle.
Let $Q'=(z_1, z_2,  x_3, x_4)$.  Lemma \ref{divi}   and the assumption imply  that  $\overline{\supp(Q')}$
  is  homeomorphic  to a   closed disk  and we can apply Corollary \ref{ctq}  to  $Q'$.

We claim $R=(2,3,8)$. 
   Suppose $R\not=(2,3,8)$.  Then $A_0\ge  \pi/12$. 
Since $\angle_{y_i}(x_i, z_i)$ ($i=1,2$) is an even angle,
 Proposition \ref{tri2}  implies   $m(y_i)\ge 4$  and   $\angle_{z_i}(x_i, y_i)\le  \pi/4$.   
It follows that $\angle_{z_1}(y_1, x_4), \angle_{z_2}(x_3, y_2)\ge  3\pi/4$  and $d(Q')\le \pi/4$.
Corollary \ref{ctq}   implies  that  $n(Q')\le 3$.  On the other hand, since 
 $y_1\in \interior(z_1z_2)$ and   $m(y_1)\ge 4$,  there are at least 4 chambers in 
$\overline{\supp(Q')}$ incident to $y_1$, a  contradiction.

 We next  claim $m(y_1)=m(y_2)=8$.  Suppose the claim is not true, say, $m(y_1)\not=8$. 
 Since $T_1=(x_1, y_1, z_1)$ has an even angle at $y_1$,  we  must have $m(y_1)=3$. 
Proposition \ref{tri2} applied to $T_1$  implies $\angle_{z_1}(x_1, y_1)=\pi/8$.
Proposition \ref{tri2} applied to $T_2$  also  implies
  $\angle_{z_2}(x_2, y_2)\le \pi/2$.  
It follows that  $\angle_{z_1}(x_4,  y_1)=7\pi/8$,  $\angle_{z_2}(x_3,  y_2)\ge \pi/2$
  and $d(Q')\le 3\pi/8$.  Corollary \ref{ctq}  then implies $n(Q')\le 9$. 
Now we count the chambers in $\overline{\supp(Q')}$  that intersect $z_1z_2$: 
 at least 7 chambers are incident to $z_1$, at least  3 are incident to $y_1$ and at most one of them is incident to both
 $z_1$ and $y_1$, at least 3 are incident to $y_2$ and at most one of them is incident to both $y_1$  and $y_2$. Hence there are at least 11 chambers in  $\overline{\supp(Q')}$, contradicting to $n(Q')\le 9$.

\begin{figure}[h]
\centering\epsfig{file=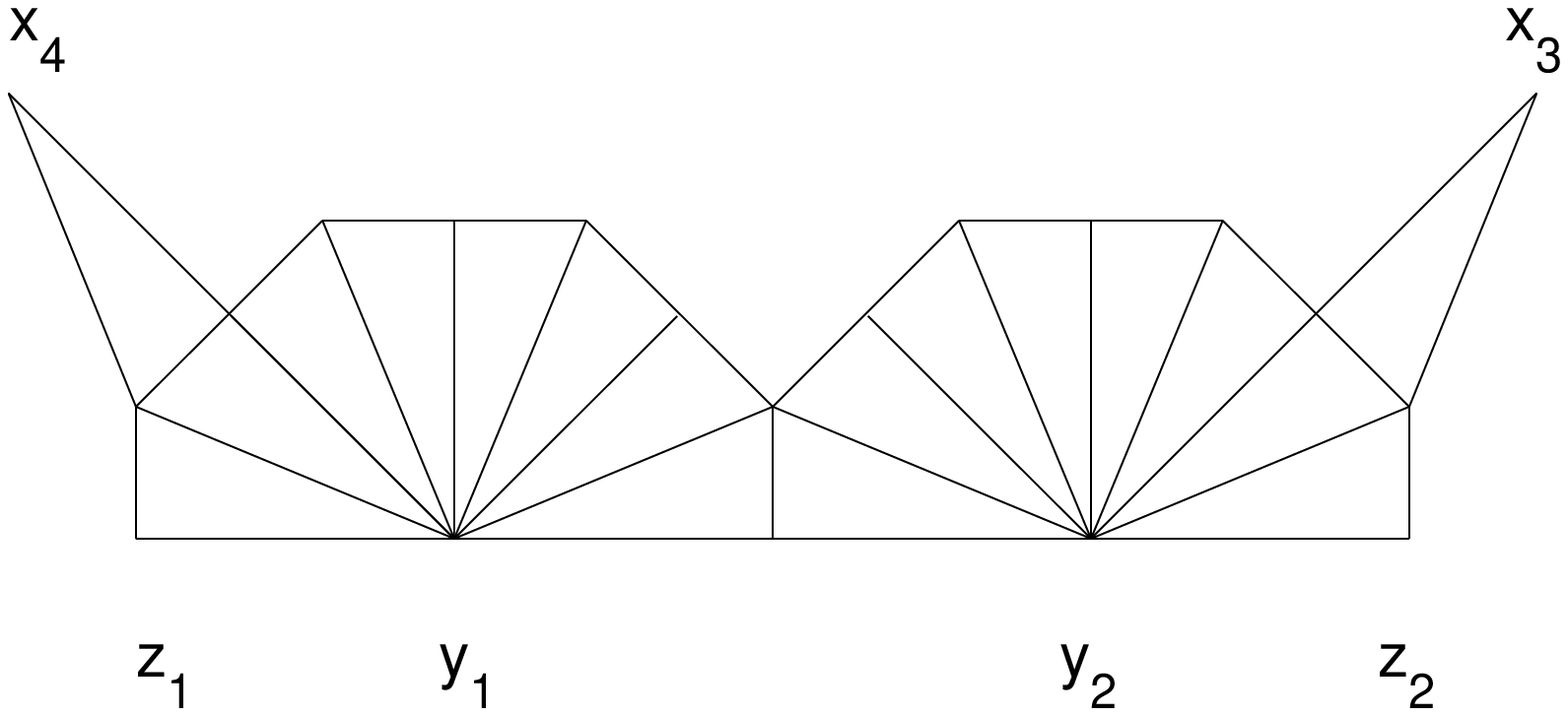,  height=4cm,  width=12cm}\caption{}\label{10l19}
\end{figure}

  Notice 
Proposition \ref{tri2}  implies  $\angle_{z_i}(x_i, y_i)\le  \pi/2$
  because    $\angle_{y_i}(x_i, z_i)$ ($i=1,2$) is an even angle.
 It follows that $\angle_{z_1}(y_1, x_4), \angle_{z_2}(x_3, y_2)\ge  \pi/2$.   
On the other hand, since  $y_1, y_2\in \interior(z_1z_2)$ and   $m(y_1)=m(y_2)=8$,  we have  $n(Q')\ge 16$.  
Corollary \ref{ctq}  then  implies that  
$\angle_{x_4}(x_3, z_1)=\angle_{x_3}(x_4,  z_2)=\pi/8$  and $\angle_{z_1}(y_1, x_4)=\angle_{z_2}(x_3,  y_2)= \pi/2$.
In particular,  $m(x_3)=m(x_4)=8$. 
      Corollary \ref{ctq} applied to $Q'$ again shows that
  $n(Q')=18$.  We count the chambers in $\overline{\supp(Q')}$  that are incident to vertices indexed by 8:
 at least 8 chambers are incident to each of $y_1, y_2$, and  at least one is incident to each of $x_3$, $x_4$, for a total of 18. 
 It follows that $\overline{\supp(Q')}$  is the union of these 18 chambers. 
Since  $\angle_{z_1}(y_1, x_4)=\angle_{z_2}(x_3,  y_2)= \pi/2$   and 
$\angle_{y_i}(x_i, z_i)$  is an even angle,  
 Proposition \ref{tri2}  applied to $T_i$  shows that  $m(z_i)=2$ and $z_iy_i$  is  an edge.
It follows that  the vertices on $z_1z_2$ are periodically indexed by  2 and 8,
  and  there is exactly one vertex in $\interior(y_1y_2)$ and it is indexed by $2$.
 One also observes that  there is exactly one vertex in the interior of each of the segments $z_1x_4$,  $z_2x_3$ and it
 is  indexed by 3.   From these  observations   one  sees that the union of the above 18 chambers is as shown in 
  Figure  \ref{10l19}. This union is not  a  quadrilateral, a contradiction.

\end{proof}

\begin{figure}[h]
\centering\epsfig{file=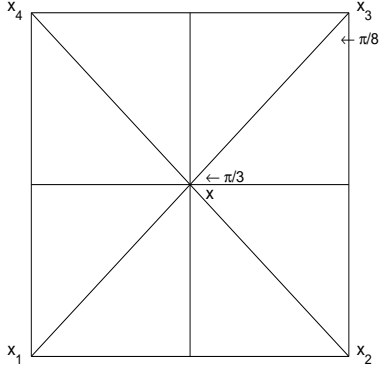,  height=5cm,  width=5cm}\caption{$x\in \supp(Q)$ is  a special point of  $\overline{\supp(Q)}$}\label{10l20}
\end{figure}

\b{Le}\label{interising}
{Let $Q=(x_1, x_2, x_3, x_4)\subset \De^{(1)}$ be a  quadrilateral that is homeomorphic to 
  a circle with $d(Q)=\frac{k\pi}{24}$, $k\ge 3$.  
Suppose $\overline{\supp(Q')}$ is  homeomorphic to  a closed disk for every 
 quadrilateral   $Q'\subset  \De^{(1)}$ homeomorphic to a  circle
 with $d(Q')\le \frac{(k-1)\pi}{24}$. 
If there is a special point $x\in \supp(Q)$, then 
 $R=(2,3,8)$  and $\overline{\supp(Q)}$ must be as shown in   Figure \ref{10l20}.}

\end{Le}


Since $x\in \supp(Q)$  is  a special point,  $Q_x$  is a circle with length $>2\pi$.
 It follows  that  there are at least 5 edges 
that are contained in $\overline{\supp(Q)}$ and incident to $x$. We extend these geodesics 
  inside  $\supp(Q)$  until they  hit $Q$. 
 Geodesic uniqueness implies these 5 points on $Q$ are  pairwise distinct.    
 At least one of these  5 points  is different from  all the $x_i's$. 
 We denote this point by $r_1$ and we may assume $r_1\in \interior(x_1x_2)$. 
 Let $\xi\in Q_x$ be the    initial   direction of $xr_1$ at $x$.  
Let $\eta_2, \eta_3\in Q_x$ be the two points in $Q_x$
that have distance $\pi$ from  $\xi$.  
We extend $r_1x$ inside  $\supp(Q)$  in the directions  $\eta_2$ and $\eta_3$  until  the extensions  
 hit $Q$ at two points $r_2$ and $r_3$ respectively. 
Note $r_2, r_3\in Q-x_1x_2$.
 Then  exactly one of the following occurs:\newline
(1) one of $r_2, r_3$ lies in $x_2x_3$  and  the other lies in $x_1x_4$, but $\{r_2, r_3\}\not=\{x_3, x_4\}$;\newline
(2) one of $r_2, r_3$  lies in $\interior(x_3x_4)$, the other lies in $\interior(x_2x_3)$ or in $\interior(x_1x_4)$;\newline
(3) $r_2, r_3$ both lie in $x_2x_3$ or both lie in $x_1x_4$;\newline
(4) $r_2, r_3\in\interior(x_3x_4)$;\newline
(5) one of $r_2, r_3$ lies in $\interior(x_3x_4)$ and the other lies in $\{x_3, x_4\}$;\newline
(6) $\{r_2, r_3\}=\{x_3, x_4\}$.

\b{Le}\label{interising1}
{Case \e{(1)} cannot occur.}

\end{Le}

\b{proof}
Suppose case (1)  occurs. We may assume $r_2\in \interior(x_2x_3)$ and $r_3\in x_1x_4$. 
Let $T_2=(r_2, x_2, r_1)$,  $T_3=(r_3, x_1, r_1)$.   
Notice the side $r_1r_i$ of $T_i$ contains at least two edges.  By Proposition \ref{tri2} 
  $R$ is  a  right triangle  and 
$\angle_{r_1}(x_2, r_2),  \angle_{r_1}(x_1, r_3)\le \pi/2$.  Using 
 the triangle inequality, one further concludes  that $\angle_{r_1}(x_2, r_2)=\angle_{r_1}(x_1, r_3)=\frac{\pi}{2}$.  
Applying Proposition \ref{tri2} again, one sees that $R=(2,3,8)$ and there are two possibilities for $T_i$:\newline
 (a)  Both   $T_2$ and  $T_3$  are as shown in    Figure \ref{238}(b), in which case $m(x)=2$ and 
 $T_i$ has angle $\pi/8$ at $r_i$;\newline
(b) Both   $T_2$ and  $T_3$  are as shown in   Figure \ref{238}(d), in which case $m(x)=3$ and 
 $T_i$ has angle $\pi/8$ at $r_i$.

   Suppose (a) occurs.  Then $r_2r_3=r_2x\cup xr_3$,  and $c=(r_2, r_3, x_4, x_3)$ 
is a quadrilateral if   $r_3\not=x_4$ and is a triangle if $r_3=x_4$.
  Since $T_i$ has angle $\pi/8$ at $r_i$ ($i=2,3$),    we have  
$d(c)\le 0$, contradicting to the fact that $\De$ is  a $\CAT(-1)$ space.

Suppose (b) holds.  First 
assume  $r_3=x_4$.  Then $Q'=(x_4, x, r_2, x_3)$ is a quadrilateral with $d(Q')<d(Q)$.  The assumption implies that 
 $\overline{\supp(Q')}$  is  a closed disk and we can apply Corollary \ref{ctq}.
   Note $\angle_x(r_2, r_3)=2\pi/3$ and $\angle_{r_2}(x, x_3)=7\pi/8$.
  Hence  $d(Q')\le 5\pi/24$ and $n(Q')\le 5$. On the other hand, there are at least 7 chambers in $\overline{\supp(Q')}$
  that  are incident to $r_2$,  a contradiction. 

Now assume $r_3\not=x_4$.    Proposition \ref{tri1} implies   $T_{2x}=\Link(\overline{\supp(T_2)}, x)$  has length $\pi$.
 The segment  $\eta_2\eta_3\subset \Link(\De, x)$  has length $2\pi/3$.
 Since $\angle_x(r_3, r_1)=\angle_x(r_2, r_1)=\pi$,   $\eta_2\eta_3\cap T_{2x}$ 
 either equals $\{\eta_2\}$ or is a segment with length $\pi/3$. 
 First suppose $\eta_2\eta_3\cap T_{2x}=\{\eta_2\}$. 
Let $y$ be the only vertex in $\interior(r_2x_2)$. Then  $m(y)=2$  and $r_3y=r_3x\cup xy$. 
Let $Q_1=(r_3, y, x_3, x_4)$. Then $d(Q_1)< d(Q)$ and $\overline{\supp(Q_1)}$  is  a closed disk. 
 Note $d(Q_1)\le 3\pi/8$.  Hence $n(Q_1)\le 9$.  Let us count the chambers in $\overline{\supp(Q_1)}$:
  since  $m(r_2)=m(r_3)=8$ and  $r_2\in  \interior(yx_3)$,   there are 
 at least 7 chambers  incident to $r_3$ and at least 8  incident to $r_2$.
 Hence $n(Q_1)\ge 15$, contradiction.

Now suppose $\eta_2\eta_3\cap T_{2x}$ is a segment with length $\pi/3$. 
 In this case $(x, y, r_3)$ is the boundary of a chamber. 
 In particular,   $r_3y$ is  an edge.   The fact $m(y)=2$ implies  $\angle_y(r_3, x_3)=\pi$  and 
  $T=(r_3, x_3, x_4)\subset  \De^{(1)}$.
Since $\angle_{r_3}(x, x_1)=\angle_{r_3}(x, y)=\pi/8$,  we conclude $\angle_{r_3}(x_4, x_3)\ge 6\pi/8$ and $d(T)\le 0$, 
 which is impossible.

\end{proof}

\b{Le}\label{interising2}
{Case \e{(2)} cannot occur.}

\end{Le}

\b{proof}
Suppose case (2)  occurs. We may assume
$r_2\in \interior(x_2x_3)$ and $r_3\in \interior(x_3x_4)$.
Let $T=(r_1, r_2, x_2)$,  $Q_1=(x_1, r_1, r_3, x_4)$ and $Q_2=(x, r_2, x_3, r_3)$.
Then $d(Q_i)<d(Q)$ ($i=1,2$) and $\overline{\supp(Q_i)}$ is  a  closed   disk. 
Since $x\in \interior(r_1r_2)$,   Proposition \ref{tri2} applied to $T$ 
  implies $R$ is a right  triangle.
Note $m(x)\not=2$: otherwise  $\angle_x(r_2, r_3)=\pi$,  
$xr_1$ and $xr_3$ are extensions of $r_2x$ and we have case (1), which does not
 occur.  Proposition \ref{tri2} then implies $R=(2,4,8), (2,4,6)$,  or  $ (2,3,8)$. 
 Assume $R=(2,4,8)$. Then $A_0=\pi/8$, $m(x)=4$  and $\angle_{r_2}(x, x_2)=\pi/8$.
  It follows that $\angle_{r_2}(x, x_3)=7\pi/8$. Note 
$\angle_x(r_2, r_3)\ge \pi/2$ since $m(x)=4$ and $\angle_x(r_2, r_3)$  is a  nonzero even angle.
Now $d(Q_2)\le 3\pi/8$ and $n(Q_2)\le 3$,  contradicting to the fact that there are at least 7 chambers in
 $\overline{\supp(Q_2)}$ incident to $r_2$.
Similarly $R\not=(2,4,6)$. Hence $R=(2,3,8)$.

We claim $m(v)\not=8$ for  every  vertex $v\in \interior(r_1r_2)$;  in  particular,  $m(x)=3$. 
 Suppose there is  a vertex  $v\in \interior(r_1r_2)$  with   $m(v)=8$. 
 Then  Proposition \ref{tri2} applied to $T$  implies $\angle_{r_2}(x, x_2)=\angle_{r_1}(x, x_2)=\pi/8$
  and $m(v)=2$ or $8$ for all vertices $v$ on $r_1r_2$. 
 It follows  that $\angle_{r_1}(x, x_1)=\angle_{r_2}(x, x_3)=7\pi/8$  and  $m(x)=8$ (since $m(x)\not=2$). 
We count the chambers in $\overline{\supp(Q_1)}$: at least 7 chambers are incident to $r_1$ and at least 8 are 
 incident to $x$. Hence $n(Q_1)\ge 15$. Corollary \ref{ctq} applied to $Q_1$ implies $\angle_{r_3}(x, x_4)\le \pi/4$.
 It follows that $\angle_{r_3}(x, x_3)\ge 3\pi/4$. Since $\angle_x(r_2, r_3)\ge \pi/4$ and 
$\angle_{r_2}(x, x_3)=7\pi/8$,  we conclude $d(Q_2)\le 0$, impossible.

Since $m(x)=3$, the two vertices in $r_1r_2$ that are adjacent to  $x$
  are  indexed by  $2$ and $8$ respectively.  It follows from the last paragraph that one of the following occurs:\newline
(a) $m(r_2)=8$ and $xr_2$ is an edge; \newline
(b) $m(r_1)=8$ and $xr_1$ is  an edge. \newline
First assume (a) holds. Then Proposition  \ref{tri2} implies $\angle_{r_2}(x, x_2)=\pi/8$ or $2\pi/8$. It follows that
$\angle_{r_2}(x, x_3)=7\pi/8$ or $6\pi/8$.
 If $\angle_{r_2}(x, x_3)=7\pi/8$, then $n(Q_2)\le 5$, contradicting to the fact that there are  at least 
   7 chambers in
 $\overline{\supp(Q_2)}$ incident to $r_2$. If $\angle_{r_2}(x, x_3)=6\pi/8$,  then $n(Q_2)\le 8$.
  A similar  argument shows   $\angle_{x_3}(r_3, r_2)<\pi/3$, in particular, $m(x_3)\not=3$. 
 It follows that there is a  vertex $v\in   \interior(r_2x_3)$ with $m(v)=3$.
 Now there are 6 chambers in $\overline{\supp(Q_2)}$ incident to $r_2$ and 2 chambers incident to $v$ but not to 
 $r_2$. Hence $\overline{\supp(Q_2)}$ is the union of these 8 chambers. However, this union is 
 not a  quadrilateral, a contradiction.

Now we assume (b) holds.  
Then $\angle_{r_1}(x, x_1)=7\pi/8$ or $6\pi/8$.  On the other hand, either 
$\angle_{r_3}(x, x_4)\ge \pi/2$ or $\angle_{r_3}(x, x_3)\ge \pi/2$.
Assume  $\angle_{r_3}(x, x_4)\ge \pi/2$   and $\angle_{r_1}(x, x_1)=7\pi/8$.  Then $n(Q_1)\le 9$.  
There are 7 chambers in $\overline{\supp(Q_1)}$ incident to $r_1$ and 2 chambers incident to $x$ but not to $r_1$.
 It follows that $\overline{\supp(Q_1)}$ is the union of these 9 chambers. However, 
this union is not a quadrilateral, a contradiction. Now assume $\angle_{r_3}(x, x_4)\ge \pi/2$  
 and $\angle_{r_1}(x, x_1)=6\pi/8$.  Then $\angle_{r_1}(x, x_2)=2\pi/8$. 
  Proposition \ref{tri2}
  implies that $\angle_{r_2}(x, x_2)=\pi/8$. Hence $\angle_{r_2}(x, x_3)=7\pi/8$. 
Also note $\angle_{x}(r_2, r_3)=2\pi/3$. It follows that $n(Q_2)\le 5$, contradicting to 
 $\angle_{r_2}(x, x_3)=7\pi/8$.
 Therefore we have $\angle_{r_3}(x, x_4)< \pi/2$. 
 It follows that $\angle_{r_3}(x, x_3)\ge 5\pi/8$.  Since $\angle_{r_2}(x, x_2)\le \pi/2$, we have $\angle_{r_2}(x, x_3)\ge \pi/2$.
  Consequently  $d(Q_2)\le \pi/12$ and $n(Q_2)\le 2$. It follows that $\overline{\supp(Q_1)}$ is the union of 
  the  two chambers incident to 
  $x$. In particular, $m(r_3)=2$, contradicting to $\angle_{r_3}(x, x_3)\ge 5\pi/8$.

\end{proof}

\b{Le}\label{interising3}
{Case \e{(3)} cannot occur.}

\end{Le}

\b{proof}
Suppose case (3)  occurs.
We may assume $r_2, r_3\in x_2x_3$  and $r_2\in \interior(r_3x_2)$.
 Let $T=(r_1,  r_3, x_2)$. Since $r_2\in \interior(r_3x_2)$ and $x\in  \interior(r_3r_1)$,
   Proposition  \ref{tri1} implies $\angle_x(r_1, r_2)<\pi$,  
 contradicting to the definition of $r_2$.

\end{proof}

\b{Le}\label{interising4}
{Case \e{(4)} cannot occur.}

\end{Le}

\b{proof}
Suppose case (4)  occurs.
We may assume $r_3\in \interior(x_4r_2)$.
Let $T=(x, r_2, r_3)$,  $Q_1=(x_1, r_1, r_3, x_4)$ and $Q_2=(x_2, r_1, r_2, x_3)$.
  Then $d(Q_i)<d(Q)$ ($i=1,2$)  and $\overline{\supp(Q_i)}$ is  a  closed  disk. 
Note $T$ has an even angle at $x$.   Proposition \ref{tri2} implies  $R$ is a right  triangle.
We claim $R=(2,3,8)$, that is, $R\not=(2,8,8)$,
$(2,6,8)$,  $(2,6,6)$,  $(2,4,8)$,
  $(2,4,6)$. We only show $R\not=(2,4,8)$,  the other cases can be handled similarly.
Assume $R=(2,4,8)$.   Then $A_0=\pi/8$.   At least one of $\angle_{r_1}(x, x_1), \angle_{r_1}(x, x_2)\ge \pi/2$.
 We may assume $\angle_{r_1}(x, x_2)\ge \pi/2$.
Proposition \ref{tri2} applied to $T$ shows  $\angle_{r_2}(x, r_3)\le \pi/4$, which 
implies $\angle_{r_2}(x, x_3)\ge 3\pi/4$.  It follows that $d(Q_2)\le \pi/2$  and $n(Q_2)\le 4$.  
On the other hand, since $m(x)\ge 4$ and $\angle_{r_2}(x, x_3)\ge 3\pi/4$, there are at least 4 chambers in 
$\overline{\supp(Q_2)}$ incident to $x$,  and at least 2 chambers 
 incident to $r_2$ but not to $x$. Hence $n(Q_2)\ge 6$, contradicting to $n(Q_2)\le 4$.

We claim $m(x)=8$. Assume $m(x)\not=8$.  Since $T$ has an even angle at $x$, $m(x)=3$.
Proposition \ref{tri2}  
implies  $\angle_{r_2}(r_3, x)=\angle_{r_3}(r_2, x)=\pi/8$. Then $\angle_{r_2}(x_3, x)=\angle_{r_3}(x_4, x)=7\pi/8$.
At least one of 
$\angle_{r_1}(x, x_1)$, $\angle_{r_1}(x_2, x)$  is $\ge \pi/2$.  We may assume 
$\angle_{r_1}(x_2, x)\ge \pi/2$.  Then $d(Q_2)\le 3\pi/8$ and $n(Q_2)\le 9$.  
We count the chambers in  $\overline{\supp(Q_2)}$:  
at least 7 chambers are 
 incident to $r_2$ and at least 2 chambers are incident to
 $x$ but not to $r_2$. Hence $\overline{\supp(Q_2)}$   must be the union of these 9 chambers. 
 But this   union  is  not a quadrilateral.   Contradiction.

Assume $T$ is as shown in   Figure   \ref{238}(d). We may assume $m(r_3)=2$ and $\angle_{r_2}(x, r_3)=\pi/8$.
  Since $m(x)=8$, by counting the chambers in $\overline{\supp(Q_2)}$ incident to $x$ and $r_2$, we see
 $n(Q_2)\ge 15$.  Corollary \ref{ctq} then implies $\angle_{r_1}(x, x_2)\le \pi/4$.  
 It follows that $m(r_1)=8$  and $\angle_{r_1}(x, x_1)\ge 3\pi/4$.  Note $\angle_{r_3}(x_4, x)=\pi/2$.  
Corollary \ref{ctq}  applied to $Q_1$ implies  $n(Q_1)\le 12$. On the other hand,
 by counting the chambers in  $\overline{\supp(Q_1)}$ incident to $x$ and $r_1$ we conclude $n(Q_1)\ge 14$,
  a contradiction.

Assume $T$ is as shown in   Figure  \ref{238}(f).  
  Then  $\angle_{r_2}(r_3, x)=\angle_{r_3}(r_2, x)=\pi/3$.
 At least one of 
$\angle_{r_1}(x, x_1)$, $\angle_{r_1}(x_2, x)$  is $\ge \pi/2$.  We may assume 
$\angle_{r_1}(x_2, x)\ge \pi/2$.  Corollary \ref{ctq}  implies  $n(Q_2)\le 14$.  
It follows   that  $m(v)\not=8$ for all vertices  $v\in \interior(xr_1)$.
Notice the vertices on $r_1r_2$ are periodically indexed by $8,3,2,3$.
It implies that  there are 
 at most three  vertices in $\interior(xr_1)$.   Assume there is no or exactly two vertices in 
$\interior(xr_1)$.  
Then 
$\angle_{r_1}(x_2, x)=2\pi/3$.    Corollary \ref{ctq}  implies    $n(Q_2)\le 10$.  By counting the chambers in 
 $\overline{\supp(Q_2)}$ intersecting $r_1r_2$ we see there cannot be exactly two vertices in 
$\interior(xr_1)$  and hence   $xr_1$ is  an edge.  In this case, there are at least 8 chambers in
$\overline{\supp(Q_2)}$ incident to $x$, 1 incident to $r_2$ but not to $x$ and 1 incident to $r_1$ but not to $x$.
 It follows that $\overline{\supp(Q_2)}$
is the union of these 10 chambers. However, this union is not a quadrilateral, a contradiction. 
The counting argument also shows that it is impossible to have $m(r_1)=8$ (with $\angle_{r_1}(x, x_2)\ge \pi/2$). 
 Hence there is exactly one vertex $v\in \interior(xr_1)$ and $m(v)=3$, $m(r_1)=2$.
Since $n(Q_2)\le 14$ and $x\in \interior(r_1r_2)$ with $m(x)=8$,  $m(z)\not=8$ for every  vertex 
$z\in \overline{\supp(Q_2)}$ different from $x$  and  the corners. 
   Let  $C$  be  the chamber    in $\overline{\supp(Q_2)}$ containing $vr_1$.  $C$ has a vertex $v'\not=x$ with $m(v')=8$.
 Notice $v'\in r_1x_2$.  It follows that  $v'=x_2$. 
Consider the three chambers in $\overline{\supp(Q_2)}$  that are incident to $v$.  
Two of them are incident to $x_2$.  It follows that 
    $\angle_{x_2}(r_1, x_3)\ge \pi/4$.   Corollary \ref{ctq} then  implies   $n(Q_2)\le 11$. 
One observes  that there are at least 11 chambers in $\overline{\supp(Q_2)}$
 intersecting $r_1r_2$, but their union is not a quadrilateral, a  contradiction.

For the remaining cases, $\angle_{r_2}(r_3, x)$  and     $\angle_{r_3}(r_2, x)\le \pi/4$  hold.
 In particular, $m(r_2)=m(r_3)=8$  and    Corollary \ref{ctq}   implies 
    $n(Q_1)\le 12$ or $n(Q_2)\le 12$.  We may assume $n(Q_2)\le 12$.  
 On the other hand, $n(Q_2)\ge 14$ since 
 there are at least 8 chambers in $\overline{\supp(Q_2)}$ incident to $x$ and at least 6 incident to $r_2$.
 A contradiction.

\end{proof}

\b{Le}\label{interising5}
{Case \e{(5)} cannot occur.}

\end{Le}

\b{proof}
Suppose case (5) does occur. 
 We may assume $r_2=x_3$ and $r_3\in \interior(x_3x_4)$. 
Let $T_1=(r_1, x_2, x_3)$,  $T_2=(x, x_3, r_3)$ and  $Q'=(x_1, r_1,r_3, x_4)$.
  Then $d(Q')<d(Q)$ and $\overline{\supp(Q')}$ is  a  closed  disk.  
Note $T_2$ has an even angle at $x$.   It follows   that  
  $R$ is a right  triangle and 
 $m(x)\not=2$.  We claim $R=(2,3,8)$.
Suppose $R\not=(2,3,8)$. Then $A_0\le \pi/12$. 
 Proposition \ref{tri2} applied to $T_1$ and $T_2$ implies 
$\angle_{r_3}(x, x_3), \angle_{r_1}(x, x_2)\le  \pi/4$. It follows 
that $\angle_{r_1}(x, x_1), \angle_{r_3}(x, x_4)\ge 3\pi/4$  and  $d(Q')\le\pi/4$. 
Corollary \ref{ctq}  then implies  $n(Q')\le 3$.
On the other hand, since $R\not=(2,3,8)$ is a right triangle and  $m(x)\not=2$, we have
 $m(x)\ge 4$. It follows that $n(Q')\ge 4$, contradicting to $n(Q')\le 3$.

Assume $m(x)=3$. Recall  $T_2$ has an even angle at $x$.
Proposition \ref{tri2} applied to $T_2$ and $T_1$ implies 
  $\angle_{r_3}(x, x_3)=\pi/8$,  
  $\angle_{r_1}(x, x_2)\le \pi/2$. 
 It follows that $\angle_{r_3}(x_4, x)=7\pi/8$ and $\angle_{r_1}(x_1,x)\ge \pi/2$.
 Consequently,  $d(Q')\le 3\pi/8$ and $n(Q')\le 9$.  There are at least 7 chambers in $\overline{\supp(Q')}$
incident to $r_3$ and at least 2 incident to $x$ but not to  $r_3$.  Hence 
 $\overline{\supp(Q')}$ must be the union of these 9 chambers.  However, this union is not a quadrilateral, 
 a contradiction.   Therefore  $m(x)=8$.  
Then Proposition \ref{tri2}
 implies $\angle_{r_1}(x, x_2)=\pi/8$,  $\angle_{r_3}(x, x_3)\le \pi/2$. 
In particular,  $m(r_1)=8$.
 It follows that $\angle_{r_1}(x_1, x)=7\pi/8$,  $d(Q')\le 3\pi/8$ and $n(Q')\le 9$. 
 On the other   hand, there are 8 chambers in $\overline{\supp(Q')}$ incident to $x$ and 7 incident to $r_1$.  Hence  
$n(Q')\ge 15$, 
  a contradiction.

\end{proof}

\b{Le}\label{interising6}
{Suppose case \e{(6)} occurs. Then $R=(2,3,8)$  and $\overline{\supp(Q)}$ must be as shown in Figure \ref{10l20}.}

\end{Le}

\b{proof}
We may assume $r_2=x_3$ and $r_3=x_4$.  Let $T_1=(x_1, r_1, x_4)$,  $T_2=(x_2, r_1, x_3)$
  and $T_3=(x, x_3, x_4)$.  Note $T_3$ has an even angle at $x$.
It follows   that $R$ is a right  triangle and 
 $m(x)\ge 3$.  We claim $R=(2,3,8)$. Suppose $R\not=(2,3,8)$. Then Proposition \ref{tri2} applied to $T_i$ ($i=1,2$) 
 implies $\angle_{r_1}(x, x_i)\le \pi/6$. By the  triangle inequality $\angle_{r_1}(x_1, x_2)\le \pi/3$, 
contradicting to $r_1\in x_1x_2$.

We claim  $m(v)\not=8$ for every vertex $v$
in the interior of $r_1x_3$ or  $r_1x_4$.
  Suppose   $m(v)=8$ for some vertex  $v\in \interior(r_1x_3)$.
Then Proposition \ref{tri2}    applied to $T_2$  implies  
 $\angle_{r_1}(x_2, x)=\pi/8$. It follows that $\angle_{r_1}(x, x_1)\ge 7\pi/8$ and 
 $d(T_1)<0$, which is impossible.  In particular, $m(x)\not=8$ and hence $m(x)=3$.
Since $T_3$ has an even angle at $x$,  Proposition \ref{tri2} implies 
$\overline{\supp(T_3)}$ is as shown in   Figure  \ref{238}(e). On the other hand,  Proposition \ref{tri2} applied to 
$T_1$ and $T_2$ shows $\angle_{r_1}(x, x_1)\le \pi/2$, $\angle_{r_1}(x, x_2)\le \pi/2$. 
Triangle inequality further implies $\angle_{r_1}(x, x_1)=\angle_{r_1}(x, x_2)= \pi/2$.
Now Proposition \ref{tri2} applied to 
$T_1$ and $T_2$  again implies that $\overline{\supp(T_i)}$ ($i=1,2$) must  be 
 as shown in   Figure   \ref{238}(d), with 
 $m(x)=3$.  The lemma is proved.

\end{proof}

The proof of 
Proposition \ref{quar11} is now complete.




\begin{thebibliography}{99}






\bibitem[BB]{BB} W. Ballmann, S. Buyalo,
\e{Nonpositively curved metrics on $2$-polyhedra,}
Math. Z.  {\bf{222}}  (1996),  no. 1, 97--134.



\bibitem[Bea]{Bea} A. Beardon,
\emph{The geometry of discrete groups,}   
   Graduate Texts in Mathematics, {\bf{91}}. Springer-Verlag, New York, 1995.

\bibitem[Ben]{Ben} N. Benakli,
\emph{Polyedres hyperboliques, passage du local au global,} Thesis, Universite Paris Sud, 1992.



\bibitem[Bo1]{Bo1}M. Bourdon,
\emph{Immeubles hyperboliques, dimension conforme et rigidité de Mostow,}
Geom. Funct. Anal.  {\bf{7}}  (1997),  no. 2, 245--268.


\bibitem[Bo2]{Bo2}M. Bourdon,
\emph{Sur les immeubles fuchsiens et leur type de quasi-isométrie,}
Ergodic Theory Dynam. Systems  {\bf{20}}  (2000),  no. 2, 343--364.


\bibitem[BP1]{BP1}M. Bourdon, H. Pajot,
 \emph{Poincaré inequalities and quasiconformal structure on the boundary of some hyperbolic buildings,}
Proc. Amer. Math. Soc.  {\bf{127}}  (1999),  no. 8, 2315--2324.


  



\bibitem[BP2]{BP2}M. Bourdon, H. Pajot,
 \emph{Rigidity of quasi-isometries for some hyperbolic buildings,}
  Comment. Math. Helv.  {\bf{75}}  (2000),  no. 4, 701--736.


\bibitem[BP3]{BP3}M. Bourdon, H. Pajot,
\e{Quasi-conformal geometry and hyperbolic geometry,} 
in  \lq\lq Rigidity in dynamics and geometry", (Cambridge, 2000),  1--17, Springer, Berlin, 2002.


 

	












    

\bibitem[BH]{BH}  M. Bridson,   A. Haefliger,
\emph{Metric spaces of nonpositive curvature,}   Grundlehren {\bf{319}},
  Springer-Verlag, Berlin(1999).
















\bibitem[F]{F}H. Federer,
\e{Geometric measure theory,}
  Grundlehren  {\bf{153}}, 
 Springer-Verlag, New York, 1969. 


	
\bibitem[FH]{FH}W. Feit, G. Higman,
\emph{The nonexistence of certain generalized polygons,  }  
J. Algebra  {\bf{1}} (1964) 114--131.
     

 



\bibitem[G]{G} M. Gromov,
\emph{Hyperbolic Groups, }  
   in   \lq\lq Essays in Group Theory" (ed. S.Gersten) M.S.R.I. Publications No.{\bf{8}},  
     Springer-Verlag(1987) 75--263.


\bibitem[GP]{GP} D. Gaboriau, F. Paulin,
\e{Sur les immeubles hyperboliques,}   Geom. Dedicata  {\bf{88}}  (2001),  no. 1-3, 153--197. 


\bibitem[H]{H} J. Heinonen,
\e{Lectures on analysis on metric spaces,}
 Universitext. Springer-Verlag, New York, 2001.

 


\bibitem[HK]{HK} J. Heinonen,  P. Koskela,
\e{Quasiconformal maps in metric spaces with controlled geometry,}  
Acta Math.  {\bf{181}}  (1998),  no. 1, 1--61.




\bibitem[KK]{KK} B. Kleiner, M. Kapovich,
\e{Hyperbolic groups with low-dimensional boundary,}  
Ann. Sci. École Norm. Sup. (4)  {\bf{33}}  (2000),  no. 5, 647--669.



\bibitem[KL]{KL} B. Kleiner,  B. Leeb,
\emph{Rigidity of quasi-isometries for symmetric spaces and Euclidean buildings,}
 C. R. Acad. Sci. Paris Sér. I Math.  {\bf{324}}  (1997),  no. 6, 639--643.




\bibitem [M]{M}B. Maskit,  
\emph{Kleinian Groups}, Springer, 1987.







\bibitem [P]{P}P. Pansu,
\emph{Métriques de Carnot-Carathéodory et quasiisométries des espaces symétriques de rang un,}  
  Ann. of Math. (2)  {\bf{129}}  (1989),  no. 1, 1--60.






\bibitem [R]{R}M. Ronan,
\emph{Lectures on buildings,}  Perspectives in Mathematics {\bf{7}}. Academic Press, Inc., Boston, MA, 1989.




\bibitem [RT]{RT}M. Ronan,J. Tits,
\e{Building buildings,}   Math. Ann.  {\bf{278}}  (1987),  no. 1-4, 291--306.

\bibitem [S1]{S1} D. Sullivan, 
\e{The density at infinity of a discrete group of hyperbolic motions,}
Inst. Hautes Études Sci. Publ. Math. No. {\bf{50}} (1979), 171--202.




\bibitem [S2]{S2} D. Sullivan, 
\e{Discrete conformal groups and measurable dynamics,}
      Bull. Amer. Math. Soc. (N.S.)  {\bf{6}}  (1982), no. 1, 57--73.




\bibitem [T]{T}J. Tits,
\emph{Buildings of spherical type and finite BN-pairs,}  Lecture Notes 
in Mathematics, Vol. {\bf{386}}. Springer-Verlag, Berlin-New York, 1974.


\bibitem [X]{X}X. Xie,
\emph{The Tits boundary of a $\CAT(0)$  2-complex,}  
 to appear in Transactions of AMS.



	  



 
 
\end{thebibliography}
\end{document}